\documentclass[preprint,12pt]{elsarticle}
\usepackage{amsmath}
\usepackage{amsfonts}
\usepackage{amssymb}
\usepackage{array,longtable}
\usepackage{amscd}
\usepackage[cp1251]{inputenc}
\usepackage[english]{babel}
\usepackage[usenames]{color}
\usepackage{graphicx}
\graphicspath{{pictures/}}
\DeclareGraphicsExtensions{.pdf,.png,.jpg}
\usepackage[matrix,arrow,curve]{xy}

\usepackage{amssymb}
 \usepackage{amsthm}

\journal{Advances in Mathematics}

\begin{document}
\newtheorem{definition}{Definition}
\newtheorem{theorem}{Theorem}
\newtheorem{lemma}{Lemma}
\newtheorem{proposition}{Proposition}
\newtheorem{observation}{Observation}
\newtheorem{corollary}{Corollary}
\newcommand{\diag}{{\rm diag}}
\begin{frontmatter}



\title{Geometric properties of LMI regions}


\author{Olga Y. Kushel}

\address{Shanghai University, \\ Department of Mathematics, \\ Shangda Road 99, \\ 200444 Shanghai, China \\
kushel@mail.ru}

\begin{abstract}
LMI (Linear Matrix Inequalities) regions is an important class of convex subsets of $\mathbb C$ arising in control theory. An LMI region $\mathfrak D$ is defined by its matrix-valued characteristic function $f_{\mathfrak D}(z) = {\mathbf L} + z{\mathbf M}+\overline{z}{\mathbf M}^T$ as follows: ${\mathfrak D} := \{z \in {\mathbb C}: f_{\mathfrak D}(z)\prec 0\}$. In this paper, we study LMI regions from the point of view of convex geometry, describing their boundaries, recession cones, lineality spaces and other characteristic in terms of the properties of matrices $\mathbf M$ and $\mathbf L$. Conversely, we study the link between the properties of matrices $\mathbf M$ and $\mathbf L$, e.g. normality, positive and negative definiteness, and the corresponding properties of an LMI region $\mathfrak D$. We provide the conditions, when an LMI region coincides with the intersection of elementary regions such as halfplanes, stripes, conic sectors and sides of hyperbolas. We also analyze the following problem, connected to pole placement: for a given LMI region $\mathfrak D$, defined by $f_{\mathfrak D}$, how to find a closed disk $D(x_0, r)$ centered at the real axis, such that $D(x_0, r) \subseteq {\mathfrak D}$?

\end{abstract}



\end{frontmatter}


\section{Introduction}
It is well-known (see \cite{GUT2}, \cite{GUJU}, \cite{MAZ}), that transient properties of a linear dynamical system are defined by its eigenvalues localization inside some given region of the complex plane. However, if the corresponding region is of a polynomial nature, such properties are difficult to analyze. A prominent idea to define an intersection of polynomial regions by a linear matrix inequality was proposed by Chilali and Gahinet in \cite{CHG}, where the following kind of regions was introduced.

{\bf Definition 1}. Let ${\mathcal M}^{n \times n}$ denote the set of all real $n \times n$ matrices. A subset ${\mathfrak D} \subset {\mathbb C}$ that can be defined as
\begin{equation}\label{LMI} {\mathfrak D} = \{z \in {\mathbb C}: \ {\mathbf L} + {\mathbf M}z+{\mathbf M}^T\overline{z} \prec 0\},\end{equation}
where ${\mathbf L}, {\mathbf M} \in {\mathcal M}^{n \times n}$, ${\mathbf L}^T = {\mathbf L}$, is called an {\it LMI region} with the {\it characteristic function} $f_{\mathfrak D}(z) = {\mathbf L} + z{\mathbf M}+\overline{z}{\mathbf M}^T$ (see \cite{CHG}, \cite{CGA}).

 Well-known examples of LMI regions are the left-hand side of the complex plane $${\mathbb C}^- = \{\lambda \in {\mathbb C} : {\rm Re}(\lambda) < 0\},$$
 with the characteristic function $$f_{{\mathbb C}^-}(z) = z + \overline{z},$$
  and the unit disk
 $$D(0,1) =  \{\lambda \in {\mathbb C} : |\lambda| < 1\},$$
 with the characteristic function  $$f_{D(0,1)}(z) = \begin{pmatrix} -1 & 0 \\  0 & -1 \\ \end{pmatrix} + \begin{pmatrix} 0 & 1 \\  0 & 0 \\ \end{pmatrix}z + \begin{pmatrix} 0 & 0 \\  1 & 0 \\ \end{pmatrix}\overline{z}.$$

 Since their introduction in \cite{CHG}, LMI regions have received enormous attention in systems and control theory (see, e.g. \cite{LIY}, \cite{OST}, \cite{SCH}, \cite{ZZX} and many others). The problem of locating all the closed loop poles of a controlled system inside a specific region ${\mathfrak D} \subset {\mathbb C}$, also known as ${\mathfrak D}$-pole placement problem and the related problem of matrix $\mathfrak D$-stablity with respect to a given LMI region ${\mathfrak D}$ (a matrix $\mathbf A \in {\mathcal M}^{n \times n}$ is called {\it stable with respect to $\mathfrak D$} or simply {\it $\mathfrak D$-stable} if its spectrum $\sigma({\mathbf A})$ belongs to ${\mathfrak D}$ (see \cite{CHG})) have appeared in various applications (see \cite{CHG}, \cite{CGA}, \cite{LZH}, \cite{MAO}, \cite{MEME}, \cite{ZHA}). The particular case of $\mathfrak D$-stability and $\mathfrak D$-stabilization problem with respect to a disk $D(x_0, r)$, centered at the point $x_0 \in {\mathbb R}$ of radius $r$ is widely studied (see \cite{FK}, \cite{CHCH}, \cite{HHP}, \cite{HON}, \cite{LLK}, \cite{MAO}, \cite{SKK}, \cite{XSX} and many others). Thus a natural question arise: given an LMI region $\mathfrak D$, defined by \eqref{LMI}, how to find a closed disk $D(x_0, r)$ such that $D(x_0, r) \subseteq {\mathfrak D}$? In this paper, we provide sufficient conditions for $x_0$ and $r$ in terms of the spectral characteristics of matrices $\mathbf L$ and $\mathbf M$.

A more general question arises, how the properties of the generating matrices $\mathbf L$ and $\mathbf M$ in Formula \eqref{LMI} are connected to the properties of an LMI region $\mathfrak D$? Furthermore, when studying robust $\mathfrak D$-stability problems, we are interested in certain characteristics of an LMI region $\mathfrak D$, such as its recession cone or lineality space. Finally, if we consider some perturbations (e.g. congruence transformation) or impose some additional properties on the matrices $\mathbf L$ and $\mathbf M$, how do we change the region? We are going to study all these questions.

The outline of the paper is as follows. Section 2 collects preliminary results from the matrix theory, focusing on the properties of normal and positive definite matrices and the techniques of simultaneous reduction of matrices to some diagonal (quasi-diagonal) forms.
Section 3 deals with topological properties of LMI regions. In this section, we describe the boundary, closure and completion of a given LMI region, and represent it as an intersection of polynomial regions. We also obtain some results on the localization of an LMI region inside an intersection of certain elementary regions. Section 4 studies LMI regions as convex sets. In this section, we provide a criterion of an LMI region to be a cone, describe recession cones and lineality spaces of LMI regions, give the criterion of an LMI region to be bounded. The main result of this section is the description of the recession cone ${\mathfrak D}_{rc}$ of an LMI region $\mathfrak D$ (Theorem \ref{resco}).
 Section 5 describes an LMI region $\mathfrak D$ using the canonical forms of its generating matrices $\mathbf M$ and $\mathbf L$. In some special cases (e.g. when $\mathbf M$ is normal and commute with $\mathbf L$), we conclude that an LMI region $\mathfrak D$ is an intersection of a certain family of halfplanes, cones, horizontal stripes and hyperbolas (see Theorems \ref{comp} and \ref{comp1}). Section 6 deals with the results potentially applicable to the study of robust problems, namely, inclusion relations between LMI regions, shifts, reflections and contractions of LMI regions, estimates of the angle of their recession cones. The main result of this section deals with the circle placement problem (see Theorem \ref{placement}).
 In Section 7, we consider the characteristics of several the most studied LMI regions with their applications to the theory of dynamical systems.
 \subsection{Example}
 Here, we consider an example of a problem we can solve with the help of the techniques, developed in this paper. Given an LMI region $\mathfrak D$, defined by its characteristic function
 $$f_{\mathfrak D} = \begin{pmatrix} -1 & 0 & 0 \\ 0 & -1 & 0 \\ 0 & 0 & -1 \\ \end{pmatrix} + \begin{pmatrix} 0 & 0 & 1 \\ 0 & 1 & -1 \\ 1 & 1 & 0 \\ \end{pmatrix}z + \begin{pmatrix} 0 & 1 & 1 \\ -1 & 1 & 0 \\ 1 & 0 & 0 \\ \end{pmatrix}\overline{z}, $$
 which represents the intersection of a parabola and a cubic (see Figure 1).
 \begin{figure}[h]
\center{\includegraphics[scale=0.5]{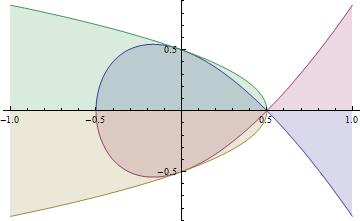}}
\caption{LMI region $\mathfrak D$}
\end{figure}
We need to find a radius $r$ such that an open disk $D(0,r)$, centered at the origin, is contained in $\mathfrak D$.

Step 1. Calculating ${\rm Sym}(\mathbf M)$ and its eigenvalues, we obtain
 $${\rm Sym}(\mathbf M) = \begin{pmatrix} 0 & 0 & 1 \\ 0 & 1 & 0 \\ 1 & 0 & 0 \\ \end{pmatrix}$$
$$\sigma({\rm Sym}(\mathbf M)) = \{-1, 1, 1\},$$
By Corollary \ref{coest}, we obtain, that the intersection ${\mathfrak D}\cap {\mathbb R} = (-\frac{1}{2}, \frac{1}{2})$.

Step 2. Calculating ${\rm Skew}(\mathbf M)$ and its eigenvalues, we obtain
 $${\rm Skew}(\mathbf M) = \begin{pmatrix} 0 & 0 & 0 \\ 0 & 0 & -1 \\ 0 & 1 & 0 \\ \end{pmatrix}$$
$$\sigma({\rm Skew}(\mathbf M)) = \{-i, i, 0\}.$$
By Theorem \ref{imag} the intersection ${\mathfrak D}\cap {\mathbb I} = (-\frac{i}{2}, \frac{i}{2})$.

Step 3. By Formula \eqref{radius}, we obtain $r({\mathfrak D}, 0) = \frac{1}{2\sqrt{2}}$ (see Figure 2).
\begin{figure}[h]
\center{\includegraphics[scale=0.5]{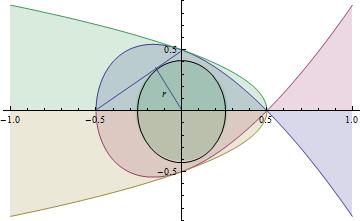}}
\caption{$D(0,r) \subset {\mathfrak D}$}
\end{figure}

\section{Preliminary results and techniques}
\subsection{Basic facts about matrices} Here, we mostly consider the matrices with real entries, denoting ${\mathbf A} \in {\mathcal M}^{n \times n}$. If a matrix $\mathbf A$ is supposed to be complex, we specify ${\mathbf A} \in {\mathcal M}^{n \times n}({\mathbb C})$.

A matrix ${\mathbf A} \in {\mathcal M}^{n \times n}$ is called
\begin{enumerate}
\item[\rm 1.] {\it symmetric} if ${\mathbf A} = {\mathbf A}^T$;
\item[\rm 2.] {\it skew-symmetric} if ${\mathbf A} = -{\mathbf A}^T$;
\item[\rm 3.] {\it orthogonal} if ${\mathbf A}{\mathbf A}^T = {\mathbf I}$;
\item[\rm 4.] {\it normal} if ${\mathbf A}{\mathbf A}^T = {\mathbf I}$.
\end{enumerate}

For an arbitrary ${\mathbf A} \in {\mathcal M}^{n \times n}$, we introduce the notations ${\rm Sym}({\mathbf A}) = \frac{{\mathbf A} + {\mathbf A}^T}{2}$ for its symmetric part and ${\rm Skew}({\mathbf A}) = \frac{{\mathbf A} - {\mathbf A}^T}{2}$ for its skew-symmetric part. Thus ${\mathbf A} = {\rm Sym}({\mathbf A}) + {\rm Skew}({\mathbf A})$, ${\mathbf A}^T = {\rm Sym}({\mathbf A}) - {\rm Skew}({\mathbf A})$.

An matrix ${\mathbf A} \in {\mathcal M}^{n \times n}({\mathbb C})$ is called
\begin{enumerate}
\item[\rm 1.] {\it Hermitian} if ${\mathbf A} = {\mathbf A}^*$, where $(\cdot)^*$ means conjugate transpose;
\item[\rm 2.] {\it unitary} if ${\mathbf A}{\mathbf A}^* = {\mathbf I}$.
\end{enumerate}

 Let $\mathcal H$ denote the set of all Hermitian matrices from ${\mathcal M}^{n \times n}({\mathbb C})$ and $\mathcal S$ denote the set of all symmetric matrices from ${\mathcal M}^{n \times n}$. Let us recall the following basic definitions and properties of Hermitian matrices (see, for example, \cite{HOJ}, \cite{GANT}, \cite{BELL}).

\begin{lemma}(see \cite{HOJ}, p. 169-172) Let ${\mathbf A}, {\mathbf B} \in {\mathcal H}$. Then
\begin{enumerate}
\item[\rm 1.] ${\mathbf A} + {\mathbf B} \in {\mathcal H}$.
\item[\rm 2.] $\alpha{\mathbf A} \in {\mathcal H}$ for every $\alpha \in {\mathbb R}$.
\item[\rm 3.]${\mathbf A}{\mathbf B} \in {\mathcal H}$ if and only if $\mathbf A$ and $\mathbf B$ commute.
\item[\rm 4.]${\mathbf A}^{-1} \in {\mathcal H}$ whenever $\mathbf A$ is nonsingular.
\item[\rm 5.] ${\mathbf S}^*{\mathbf A}{\mathbf S}\in {\mathcal H}$ for any ${\mathbf S} \in {\mathcal M}^{n \times n}({\mathbb C})$.
 \end{enumerate}
\end{lemma}

The above properties show that the class of Hermitian matrices $\mathcal H$ form a linear subspace of ${\mathcal M}^{n \times n}({\mathbb C})$. Respectively, the class of symmetric matrices $\mathcal S$ form a subspace of ${\mathcal M}^{n \times n}$.

 Two matrices ${\mathbf A} , {\mathbf B} \in {\mathcal M}^{n \times n}({\mathbb C})$ are called {\it congruent} if they are connected by the formula
$${\mathbf B} = {\mathbf S}^*{\mathbf A}{\mathbf S} $$ for some nonsingular ${\mathbf S} \in {\mathcal M}^{n \times n}({\mathbb C})$.
In case of ${\mathbf A} , {\mathbf B} \in {\mathcal M}^{n \times n}$, the congruence is defined as ${\mathbf B} = {\mathbf S}^T{\mathbf A}{\mathbf S}$, for some nonsingular ${\mathbf S} \in {\mathcal M}^{n \times n}$.

Given a matrix $\mathbf A$ with real spectrum, the {\it inertia} of $\mathbf A$ (denoted ${\rm In}({\mathbf A})$) is the ordered triple:
$$i(\mathbf A) = (i_+(\mathbf A), \ i_-(\mathbf A), \ i_0(\mathbf A)), $$
where $i_+(\mathbf A)$, $i_-(\mathbf A)$ and $i_0(\mathbf A)$ are the numbers of positive, negative and zero eigenvalues of $\mathbf A$, respectively, all counting their multiplicities. Recall a well-known fact that {\it all the eigenvalues of a Hermitian matrix are real}. The following theorem holds (see \cite{HOJ}, p. 223).

\begin{theorem}[Sylvester's law of inertia]
Let ${\mathbf A}, {\mathbf B} \in {\mathcal H}$. Then $\mathbf A$ and $\mathbf B$ are congruent if and only if $i(\mathbf A) = i(\mathbf B)$, i.e. they have the same number of positive, negative and zero eigenvalues.
\end{theorem}

\subsection{Definite matrices and their properties}
In this subsection, we consider complex Hermitian matrices. All the mentioned facts remain valid for real symmetric matrices. Here, as usual, we denote $[n] = \{1, \ \ldots, \ n\}$. Given a set of indices $\alpha = (i_1, \ \ldots, \ i_k) \subseteq [n]$ we use the notations ${\mathbf A}[\alpha]$ for the principal submatrix and ${\mathbf A}(\alpha)$ for the principal minor of ${\mathbf A}$, formed by rows and columns with the indices from $\alpha$. We use the notation ${\mathbf A}(k)$ for the $k$th leading principal minor of ${\mathbf A}$, i.e. the minor formed by rows and columns with consecutive indices $(1, \ldots, \ k)$, $1 \leq k \leq n$.

 An $n \times n$ Hermitian matrix $\mathbf A$ is called {\it positive definite (positive semidefinite)} if $\langle x, {\mathbf A}x\rangle >0$ for all nonzero $x \in {\mathbb C}^{n}$ (respectively, $\langle x, {\mathbf A}x\rangle \geq 0$ for all $x \in {\mathbb C}^{n}$). An $n \times n$ Hermitian matrix $\mathbf A$ is called {\it negative definite (semidefinite)} if $-{\mathbf A}$ is positive definite (semidefinite). Further we say that $\mathbf A$ is {\it definite} ({\it semidefinite}) if $\mathbf A$ is either positive or negative definite (semidefinite). We denote $\mathbf A \prec 0$ $(A \succ 0)$ for a negative definite (respectively, positive definite) matrix $\mathbf A$, and $\preceq 0$ $(\succeq 0)$ means negative (positive) semidefinite. The notation ${\mathbf A} \prec {\mathbf B}$ (${\mathbf A} \preceq {\mathbf B}$) means that $\mathbf A$, $\mathbf B$ are Hermitian and that ${\mathbf A} - {\mathbf B}$ is negative definite (semidefinite). For the results on definite and semidefinite matrices, we mostly refer to \cite{BHAT2} (see also \cite{BELL}).

The class ${\mathcal H}^+$ of Hermitian positive definite matrices is closed under matrix addition and multiplication by a positive constant, the same is true for the class ${\mathcal H}^+_0$ of Hermitian nonnegative definite matrices. Thus both ${\mathcal H}^+$ and ${\mathcal H}^+_0$ are convex cones in $\mathcal H$. Consider a vector space ${\mathcal M}^{n \times n}({\mathbb C})$ equipped with an operator norm
$$\|{\mathbf A}\| = \sup_{\|x\|=1}\|{\mathbf A}x\| = \sup_{\|x\|\leq1}\|{\mathbf A}x\|.$$
{\it The convex cone of positive definite matrices ${\mathcal H}^+$ is open in the subspace of Hermitian matrices ${\mathcal H}$} (\cite{BHAT2}, p. 18). The closure of ${\mathcal H}^+$ coincides with the class ${\mathcal H}^+_0$ of Hermitian positive semidefinite matrices (see, for example, \cite{HOJ3}, p. 432, Observation 7.1.9 for the Hermitian case and \cite{BOYV}, p. 43, for the symmetric case). Thus it is easy to see that ${\mathbf A} + {\mathbf A}_0 \in {\mathcal H}^+$ whenever ${\mathbf A} \in {\mathcal H}^+$, ${\mathbf A}_0 \in {\mathcal H}^+_0$. The above facts imply the corresponding properties for the following subspaces of ${\mathcal H}$. Given a partition $\alpha = (\alpha_1, \ \ldots, \ \alpha_p)$ of the set of indices $[n]$, define the subspace ${\mathcal H}(\alpha)$ of $\alpha$-diagonal matrices as follows:
 $${\mathcal H}(\alpha) = \{{\mathbf H} \in {\mathcal M}^{n \times n}({\mathbb C}): {\mathbf H} = {\rm diag}\{{\mathbf H}[\alpha_1], \ \ldots, \ {\mathbf H}[\alpha_p]\}, $$
 where ${\mathbf H}[\alpha_1], \ \ldots, \ {\mathbf H}[\alpha_p]$ are Hermitian matrices of appropriate size.
 {\it For any partition $\alpha = (\alpha_1, \ \ldots, \ \alpha_p)$ of $[n]$, the set ${\mathcal H}^+(\alpha)$ of Hermitian positive definite $\alpha$-diagonal matrices is an open convex cone in the subspace ${\mathcal H}(\alpha)$ of Hermitian $\alpha$-diagonal matrices. The cone ${\mathcal H}_0^+(\alpha)$ of Hermitian positive semidefinite $\alpha$-diagonal matrices coincides with the closure of ${\mathcal H}^+(\alpha)$}.

The class ${\mathcal H}^+$ is not closed with respect to matrix multiplication.

 \begin{lemma}(\cite{BHAT2})\label{comm}
 for ${\mathbf A}, {\mathbf B} \in {\mathcal H}^+$, the matrix product ${\mathbf A}{\mathbf B}$ belongs to ${\mathcal H}^+$ if and only if $\mathbf A$ and $\mathbf B$ commute.
 \end{lemma}

 Let us list the following equivalent characterizations of positive definite matrices (see \cite{BHAT2}, p. 1-2).
\begin{lemma}\label{propert} A Hermitian matrix ${\mathbf A} \in {\mathcal M}^{n \times n}({\mathbb C})$ is positive definite if and only if one of the following conditions holds.
\begin{enumerate}
\item[\rm 1.] ${\mathbf A} = {\mathbf B}^*{\mathbf B}$ for some nonsingular matrix ${\mathbf B} \in {\mathcal M}^{n \times n}({\mathbb C})$.
\item[\rm 2.] The principal minors of $\mathbf A$ are all positive.
\item[\rm 3.] The leading principal minors of $\mathbf A$ are all are positive, i.e. ${\mathbf A}(k) > 0$ for $k = 1, \ \ldots, \ n$ (Sylvester's criterion).
\item[\rm 4.] The eigenvalues of $\mathbf A$ are all positive.
\item[\rm 5.] ${\mathbf A} = {\mathbf B}^2$ for some positive definite matrix $\mathbf B$.
\item[\rm 6.] ${\mathbf A} = {\mathbf T}{\mathbf T}^*$ for some nonsingular lower triangular matrix $\mathbf T$ with positive principal diagonal entries (Cholesky decomposition). The matrix $\mathbf T$ is unique.
\item[\rm 7.] ${\mathbf A}^{-1}$ (the inverse of $\mathbf A$) is positive definite.
\item[\rm 8.] ${\mathbf X}^*{\mathbf A}{\mathbf X}$ is positive definite, where $\mathbf X$ is any nonsingular matrix.
\end{enumerate}
\end{lemma}

For positive semidefinite matrices, we mention the following properties (see \cite{BHAT2}).

\begin{lemma}\label{propertsemi} A Hermitian matrix ${\mathbf A} \in {\mathcal M}^{n \times n}({\mathbb C})$ is positive semidefinite if and only if one of the following conditions holds.
\begin{enumerate}
\item[\rm 1.] The principal minors of $\mathbf A$ are all nonnegative.
\item[\rm 2.] The eigenvalues of $\mathbf A$ are all nonnegative.
\item[\rm 3.] ${\mathbf X}^*{\mathbf A}{\mathbf X}$ is positive semidefinite, where $\mathbf X$ is any nonsingular matrix.
\end{enumerate}
\end{lemma}

Given a Hermitian matrix $\mathbf A$, we denote $\lambda^{\downarrow}({\mathbf A})$ the ordered set of $n$ eigenvalues of $\mathbf A$, listed in weakly decreasing order, taking into account their multiplicities:
$$ \lambda^{\downarrow}({\mathbf A}) = (\lambda_1^{\downarrow}({\mathbf A}), \ \ldots, \ \lambda_n^{\downarrow}({\mathbf A})),$$ where $\lambda_1^{\downarrow}({\mathbf A}) \leq \ldots \leq \lambda_n^{\downarrow}({\mathbf A})$.

As we see from Lemmas \ref{propert} and \ref{propertsemi}, {\it $\mathbf A$ is positive definite (semidefinite) if and only if $\lambda_n^{\downarrow}({\mathbf A})> 0$ (respectively, $\geq 0$), and $\mathbf A$ is negative definite (semidefinite) if and only if $\lambda_1^{\downarrow}({\mathbf A})< 0$ (respectively, $\leq 0$)}.

Recall the following statement (for non-strict inequalities, see \cite{BHAT}, p. 62, Theorem III.2.1, for the conditions, when Weyl's inequalities are strict, see \cite{HORR}, p. 33, Theorem 3.1).
\begin{lemma}[Weyl's inequalities]\label{Weyl}
Let ${\mathbf A}, {\mathbf B} \in {\mathcal M}^{n \times n}({\mathbb C})$ be Hermitian matrices. Then
$$\lambda_j^{\downarrow}({\mathbf A} + {\mathbf B}) \leq \lambda_i^{\downarrow}({\mathbf A}) + \lambda_{j-i+1}^{\downarrow}({\mathbf B}) \qquad \mbox{for} \ i \leq j;$$
$$\lambda_j^{\downarrow}({\mathbf A} + {\mathbf B}) \geq \lambda_i^{\downarrow}({\mathbf A}) + \lambda_{j-i+n}^{\downarrow}({\mathbf B}) \qquad \mbox{for} \ i \geq j.$$
In particular,
\begin{equation}\label{Wey}
\lambda_1^{\downarrow}({\mathbf A} + {\mathbf B}) \geq \lambda_1^{\downarrow}({\mathbf A}) + \lambda_{n}^{\downarrow}({\mathbf B});
\end{equation}
\begin{equation}\label{Wey2}
\lambda_1^{\downarrow}({\mathbf A} + {\mathbf B}) \leq \lambda_1^{\downarrow}({\mathbf A}) + \lambda_{1}^{\downarrow}({\mathbf B}).
\end{equation}
Inequality \eqref{Wey} is strict if $${\rm Ker}({\mathbf A} - \lambda_1^{\downarrow}({\mathbf A}){\mathbf I}) \cap{\rm Ker}({\mathbf B} - \lambda_1^{\downarrow}({\mathbf B}){\mathbf I}) = \{0\}.$$
In its turn, Inequality \eqref{Wey2} is strict if $${\rm Ker}({\mathbf A} - \lambda_1^{\downarrow}({\mathbf A}){\mathbf I}) \cap{\rm Ker}({\mathbf B} - \lambda_n^{\downarrow}({\mathbf B}){\mathbf I}) = \{0\}.$$
\end{lemma}
\begin{corollary}\label{strictWeyl} Inequalities \eqref{Wey} and \eqref{Wey2} are both strict if $\det({\mathbf A}{\mathbf B} - {\mathbf B}{\mathbf A}) \neq 0$.
\end{corollary}

Later we will also use the complex version of the famous Lyapunov theorem (see \cite{HOJ}, p. 96, Theorem 2.2.1).

\begin{theorem}[Lyapunov]\label{lyap} Let ${\mathbf A} \in {\mathcal M}^{n \times n}({\mathbb C})$. Then $\mathbf A$ is stable (i.e. ${\rm Re}(\lambda) < 0$ whenever $\lambda \in \sigma({\mathbf A})$) if and only if there is a Hermitian positive definite matrix ${\mathbf P} \in {\mathcal M}^{n \times n}({\mathbb C})$ such that
\begin{equation}\label{lyapeq} {\mathbf W}:={\mathbf P}{\mathbf A}^* + {\mathbf A}^*{\mathbf P}
\end{equation}
is negative definite.
\end{theorem}

For the case of real (not necessarily symmetric) matrices, further we will use the following generalization of positive definiteness, introduced in \cite{JOHN6}.

An $n \times n$ real (not necessarily symmetric) matrix $\mathbf A$ is called {\it positive definite (semidefinite)} if its symmetric part ${\rm Sym}({\mathbf A}) = \frac{{\mathbf A} + {\mathbf A}^T}{2}$ is positive definite (respectively, semidefinite). By Lyapunov theorem, such matrices are necessarily stable.

\subsection{Normal matrices and canonical forms}

A normal matrix is known to be orthogonally similar to its quasi-diagonal form (see \cite{HOJ2}, p. 101, Theorem 2.5.4, also see \cite{GANT}, p. 284).
\begin{theorem}\label{canon}
A matrix ${\mathbf A} \in {\mathcal M}^{n \times n}$ is normal if and only if
$${\mathbf A} = {\mathbf Q}{\mathbf \Lambda}_{\mathbf A}{\mathbf Q}^T,$$
where ${\mathbf Q}$ is a real orthogonal matrix (${\mathbf Q}^T = {\mathbf Q}^{-1}$),
and ${\mathbf \Lambda}_{\mathbf A}$ is a block-diagonal matrix of the following form
\begin{equation}\label{qdiag}{\mathbf \Lambda}_{\mathbf A} = {\rm diag}\{\begin{pmatrix}\mu_1 & \nu_1 \\- \nu_1 & \mu_1 \end{pmatrix}, \ldots, \begin{pmatrix}\mu_k & \nu_k \\- \nu_k & \mu_k \end{pmatrix}, \ \lambda_{2k+1}, \ \ldots, \lambda_n\}, \end{equation}
where $\lambda_{2j-1} = \mu_j+i\nu_j$, $\lambda_{2j} = \mu_j-i\nu_j$, $j = 1, \ldots, k$ are non-real eigenvalues of $\mathbf A$, $\lambda_{2k+1}, \ \ldots, \lambda_n$ are real eigenvalues of $\mathbf A$.
 \end{theorem}

Symmetric and skew-symmetric matrices are obviously normal. The following statement holds for a symmetric matrix (see \cite{HOJ}, p. 171).

\begin{theorem}\label{SpecSym}
A matrix ${\mathbf A} \in {\mathcal M}^{n \times n}$ is symmetric if and only if ${\mathbf A} = {\mathbf Q}{\mathbf \Lambda}_{\mathbf A}{\mathbf Q}^T$, where ${\mathbf Q}$ is a real orthogonal matrix and ${\mathbf \Lambda}_{\mathbf A} \in {\mathcal M}^{n \times n}$ is a diagonal matrix such that
\begin{equation}\label{canonsym}{\mathbf \Lambda}_{\mathbf A} = {\rm diag}\{\lambda_1, \ \ldots, \ \lambda_n\},\end{equation} where $\{\lambda_i\}_{i = 1}^n$ are the eigenvalues of $\mathbf A$.
\end{theorem}

 The following result holds for a skew-symmetric matrix (see \cite{HOJ2}, p. 107, Corollary 2.5.14).
 \begin{theorem}\label{canonskew}
A matrix ${\mathbf A} \in {\mathcal M}^{n \times n}$ is skew-symmetric if and only if
$${\mathbf A} = {\mathbf Q}{\mathbf \Lambda}_{\mathbf A}{\mathbf Q}^T,$$
where ${\mathbf Q}$ is a real orthogonal matrix,
and ${\mathbf \Lambda}$ is a block-diagonal matrix of the following form
\begin{equation}\label{qdiagskew}{\mathbf \Lambda}_{\mathbf A} = {\rm diag}\{\begin{pmatrix}0 & \nu_1 \\- \nu_1 & 0 \end{pmatrix}, \ldots, \begin{pmatrix}0 & \nu_k \\- \nu_k & 0 \end{pmatrix}, \ 0, \ \ldots, 0\}, \end{equation}
where $\lambda_{2j-1} = i\nu_j$, $\lambda_{2j} = -i\nu_j$, $j = 1, \ldots, k$ are non-real eigenvalues of $\mathbf A$, $0$ is the only real eigenvalues of $\mathbf A$.
 \end{theorem}

 Consider the following equivalent characteristics of normal matrices (see \cite{HOJ}, p. 109).

\begin{lemma}\label{norm} ${\mathbf A} \in {\mathcal M}^{n \times n}$ is normal if and only if:
\begin{enumerate}
\item[\rm 1.] ${\rm Sym}({\mathbf A})$ commutes with ${\rm Skew}({\mathbf A})$;
\item[\rm 2.] ${\mathbf A}$ commutes with some normal matrix with distinct eigenvalues;
\item[\rm 3.] ${\mathbf A} + t{\mathbf I}$ is normal for any $t \in {\mathbb R}$.
\item[\rm 4.] If, in addition, all the eigenvalues of $\mathbf A$ are real, $\mathbf A$ is normal if and only if $\mathbf A$ is symmetric.
\end{enumerate}
\end{lemma}
Given a symmetric matrix $\mathbf A$, let us decompose $\sigma({\mathbf A})$ as follows:
$$\sigma({\mathbf A}) = (\sigma({\mathbf A})\cap {\mathbb R}^-)\bigcup(\sigma({\mathbf A})\cap ({\mathbb R}^+\cup\{0\})) $$
and write the corresponding decomposition of ${\mathbf \Lambda}_{\mathbf A}$:
$${\mathbf \Lambda}_{\mathbf A} = {\mathbf \Lambda}_1 - x{\mathbf I} - ({\mathbf \Lambda}_2 - x{\mathbf I}),$$
where ${\mathbf \Lambda}_1$ and $- {\mathbf \Lambda}_2$ are the block-diagonal matrices consisting of block that corresponds to the negative and nonnegative eigenvalues of $\mathbf A$, respectively, $x \in {\mathbb R}, \ x > 0$. Then we obtain the following decomposition of $\mathbf A$.
\begin{equation}\label{dec}
{\mathbf A} = {\mathbf Q}{\mathbf \Lambda}_{\mathbf A}{\mathbf Q}^T = {\mathbf Q}({\mathbf \Lambda}_1 - x{\mathbf I} - ({\mathbf \Lambda}_2 - x{\mathbf I})){\mathbf Q}^T =\end{equation}$$ {\mathbf Q}({\mathbf \Lambda}_1 - x{\mathbf I}){\mathbf Q}^T - {\mathbf Q}({\mathbf \Lambda}_2- x{\mathbf I}){\mathbf Q}^T ={\mathbf A}_1 - {\mathbf A}_2.$$

Both the matrices ${\mathbf A}_1$ and ${\mathbf A}_2$ are real, symmetric (by Theorem \ref{SpecSym}) and negative definite, moreover, taking sufficiently small $x < 0$, we obtain
$\sigma({\mathbf A}_1) \rightarrow \sigma({\mathbf A})\cap {\mathbb R}^-$ and $\sigma(-{\mathbf A}_2) \rightarrow \sigma({\mathbf A})\cap ({\mathbb R}^+\cup\{0\})$.

Given an arbitrary normal matrix $\mathbf A$ and a nonsingular matrix $\mathbf S$, a congruence transformation ${\mathbf S}{\mathbf A}{\mathbf S}^T$ does not necessarily preserve normality of $\mathbf A$. Later, we are interested in the following two cases, when it does. These are:
 \begin{enumerate}
\item[\rm 1.] $\mathbf A$ is symmetric (skew-symmetric), $\mathbf S$ is arbitrary nonsingular. In this case, ${\mathbf S}{\mathbf A}{\mathbf S}^T$ preserves normality since it obviously preserves symmetry (skew-symmetry).
\item[\rm 2.] $\mathbf A$ is arbitrary normal, $\mathbf S$ is orthogonal.
\end{enumerate}

Let us recall the following well-known statement from the theory of matrices (see, for example, \cite{HOJ2}, p. 413, Corollary 7.3.3, also \cite{GANT}, \cite{BHAT}).
\begin{theorem}\label{polar} Every matrix ${\mathbf A} \in {\mathcal M}^{n \times n}$ can be written in the form
$${\mathbf A} = {\mathbf P}{\mathbf U},$$
where ${\mathbf P} \in {\mathcal M}^{n \times n}$ is positive semidefinite and ${\mathbf U} \in {\mathcal M}^{n \times n}$ is orthogonal. The matrix $\mathbf P$ is always uniquely determined as ${\mathbf P} = ({\mathbf A}{\mathbf A}^T)^{\frac{1}{2}}$, if $\mathbf A$ is nonsingular, $\mathbf U$ is also uniquely determined as ${\mathbf U} = {\mathbf P}^{-1}{\mathbf A}$.
\end{theorem}

The following fact can be easily deduced from the canonical form of normal matrices (see \cite{HOJ2}, p. 417).

\begin{lemma}\label{un} Let ${\mathbf A} \in {\mathcal M}^{n \times n}$ be normal and have the unitary diagonal representation
$${\mathbf A} = {\mathbf W}{\mathbf \Lambda}{\mathbf W}^*, \qquad {\mathbf W}{\mathbf W}^* = {\mathbf I}, $$
where ${\mathbf \Lambda} = {\rm diag}\{\rho_1e^{i\varphi_1}, \ \ldots, \ \rho_ne^{i\varphi_n}\}$, $\rho_je^{i\varphi_j} =\lambda_j \in \sigma({\mathbf A})$. Then $\mathbf A$ has the polar decomposition ${\mathbf A} = {\mathbf P}{\mathbf U}$, where ${\mathbf P} = {\mathbf W}{\mathbf \Lambda}_{\mathbb R}{\mathbf W}^*$, ${\mathbf U} = {\mathbf W}{\mathbf \Lambda}_{\varphi}{\mathbf W}^*$, ${\mathbf \Lambda}_{\mathbb R} = {\rm diag}\{\rho_1, \ \ldots, \ \rho_n\}$ and ${\mathbf \Lambda}_{\varphi} = {\rm diag}\{e^{i\varphi_1}, \ \ldots, \ e^{i\varphi_n}\}$.
\end{lemma}

\subsection{Simultaneous reduction by congruence}
Given a family ${\mathcal A} = \{{\mathbf A}_i\}_{i=1}^k$ of normal matrices from ${\mathcal M}^{n \times n}$, the matrices ${\mathbf A}_i$, $i = 1, \ \ldots, \ k$, are called {\it simultaneously quasi-diagonalizable by congruence} if there is an invertible matrix ${\mathbf S}$ such that all the matrices ${\mathbf S}{\mathbf A}_i{\mathbf S}^T$ are quasi-diagonal. If all ${\mathbf S}{\mathbf A}_i{\mathbf S}^T$ are diagonal, the matrices ${\mathbf A}_i$, $i = 1, \ \ldots, \ k$, are called {\it simultaneously diagonalizable by congruence}. If, in addition, ${\mathbf S}$ is orthogonal (unitary), we say that ${\mathbf A}_i$, $i = 1, \ \ldots, \ k$, are {\it simultaneously quasi-diagonalizable (diagonalizable) by orthogonal (unitary) congruence}. Note, than an orthogonal (unitary) congruence does not change matrix spectra.

The main problem we face in studying LMI regions is as follows.

{\bf Problem 1.} {\it Given two matrices ${\mathbf A}, {\mathbf B} \in {\mathcal M}^{n \times n}$, where ${\mathbf A}$ is normal, $\mathbf B$ is symmetric, when ${\mathbf A}$, ${\mathbf A}^T$ and ${\mathbf B}$ are simultaneously quasi-diagonalizable by (not necessarily orthogonal) congruence? Here, we may also assume $\mathbf B$ be negative definite.}

 Note, that, {\it if a congruence transformation reduces a matrix $\mathbf A$ to a quasi-diagonal form, then it also reduces ${\mathbf A}^T$ to a quasi-diagonal form.} Indeed, if ${\mathbf S}{\mathbf A}{\mathbf S}^T = {\mathbf \Lambda}$, where ${\mathbf \Lambda}$ is quasi-diagonal, by transposition we obtain $({\mathbf S}{\mathbf A}{\mathbf S})^T = {\mathbf S}{\mathbf A}^T{\mathbf S}^T = {\mathbf \Lambda}^T$, which is also quasi-diagonal.

Applying Lemma \ref{norm}, we get the following re-statement of Problem 1.

{\bf Problem 2}. {\it Given three matrices ${\mathbf A}, {\mathbf B}, {\mathbf C} \in {\mathcal M}^{n \times n}$, where ${\mathbf A}$ and $\mathbf C$ are symmetric, $\mathbf B$ is skew-symmetric and commute with ${\mathbf A}$, when all of them are simultaneously quasi-diagonalizable by (not necessarily orthogonal) congruence? Here, we may also assume $\mathbf C$ be negative definite.}

Recall the following well-known result (see \cite{HOJ2}, p. 108, Theorem 2.5.15, also \cite{GANT}, p. 292, Theorem 12').
\begin{theorem}\label{commute}
Given a commuting family ${\mathcal A} = \{{\mathbf A}_i\}_{i=1}^k$ of normal matrices from ${\mathcal M}^{n \times n}$, they can be transformed to their quasi-diagonal forms \eqref{qdiag}, using the same orthogonal transformation $Q$.
\end{theorem}
\begin{corollary}\label{commute3} Let normal matrices ${\mathbf A}, \ {\mathbf B} \in {\mathcal M}^{n \times n}$ commute. Then the matrices $\mathbf A$, ${\mathbf A}^T$ and $\mathbf B$ are simultaneously quasi-diagonalizable by orthogonal congruence.
\end{corollary}
\begin{corollary}\label{commute4} Let $\mathbf A$ be normal, $\mathbf B$ be symmetric, ${\mathbf A}{\mathbf B} = {\mathbf B}{\mathbf A}$. Then, if Form \eqref{qdiag} of $\mathbf A$ is given by $${\mathbf \Lambda}_{\mathbf A} = {\rm diag}\{\begin{pmatrix}\mu_1 & \nu_1 \\- \nu_1 & \mu_1 \end{pmatrix}, \ldots, \begin{pmatrix}\mu_k & \nu_k \\- \nu_k & \mu_k \end{pmatrix}, \ \lambda_{2k+1}({\mathbf A}), \ \ldots, \lambda_n({\mathbf A})\}, $$
where $\lambda_{2j-1}({\mathbf A}) = \mu_j+i\nu_j$, $\lambda_{2j}({\mathbf A}) = \mu_j-i\nu_j$, $j = 1, \ldots, k$ are non-real eigenvalues of $\mathbf A$, $\lambda_{2k+1}({\mathbf A}), \ \ldots, \lambda_n({\mathbf A})$ are real eigenvalues of $\mathbf A$, then the corresponding diagonal form of $\mathbf B$ is given by
$${\mathbf \Lambda}_{\mathbf B} = {\rm diag}\{\lambda_1({\mathbf B}), \ \ldots, \ \lambda_{2k}({\mathbf B}), \  \lambda_{2k+1}({\mathbf B}), \ \ldots, \lambda_n({\mathbf B})\}, $$
with $\lambda_{2j-1}({\mathbf B}) = \lambda_{2j}({\mathbf B})$, $j = 1, \ \ldots, \ k$.
\end{corollary}
{\bf Proof.} The condition ${\mathbf A}{\mathbf B} = {\mathbf B}{\mathbf A}$ implies ${\mathbf Q}({\mathbf A}{\mathbf B}){\mathbf Q}^T = {\mathbf Q}({\mathbf B}{\mathbf A}){\mathbf Q}^T$ and $({\mathbf Q}{\mathbf A}{\mathbf Q}^T)({\mathbf Q}{\mathbf B}{\mathbf Q}^T) = ({\mathbf Q}{\mathbf B}{\mathbf Q}^T)({\mathbf Q}{\mathbf A}{\mathbf Q}^T)$ for any orthogonal matrix $\mathbf Q$. Thus ${\mathbf \Lambda}_{\mathbf A}{\mathbf \Lambda}_{\mathbf B} = {\mathbf \Lambda}_{\mathbf B}{\mathbf \Lambda}_{\mathbf A}$. Since ${\mathbf \Lambda}_{\mathbf B}$ is diagonal and ${\mathbf \Lambda}_{\mathbf A}$ has block structure \eqref{qdiag}, we have the commutativity condition for each pair of $2 \times 2$ blocks:
$$\begin{pmatrix}\lambda_{2j-1}({\mathbf B}) & 0 \\ 0 & \lambda_{2j}({\mathbf B}) \end{pmatrix}\begin{pmatrix}\mu_j & \nu_j \\- \nu_j & \mu_j \end{pmatrix} = \begin{pmatrix}\mu_j & \nu_j \\- \nu_j & \mu_j \end{pmatrix}\begin{pmatrix}\lambda_{2j-1}({\mathbf B}) & 0 \\ 0 & \lambda_{2j}({\mathbf B}) \end{pmatrix}$$
for each $j = 1, \ \ldots, \ 2k$. These conditions obviously imply $\nu_j\lambda_{2j}({\mathbf B}) = \nu_j\lambda_{2j-1}({\mathbf B})$ and since $\nu_j \neq 0$, we have $\lambda_{2j-1}({\mathbf B}) = \lambda_{2j}({\mathbf B})$ for $j = 1, \ \ldots, \ k$.
 $\square$

Now, let us introduce the following notation: given an arbitrary matrix ${\mathbf A}$ and a definite matrix ${\mathbf B} \in {\mathcal M}^{n \times n}$, denote ${\mathbf A}_{\mathbf B}:= {\mathbf T}^{-1}{\mathbf A}({\mathbf T}^{-1})^T$, where $\mathbf T$ is a lower triangular matrix from the Cholesky decomposition ${\mathbf B} = \pm{\mathbf T}{\mathbf T}^T$.

We are also interested in conditions sufficient for two normal matrices ${\mathbf A}$ and ${\mathbf B}$ to be simultaneously quasi-diagonalized by congruence. Consider the following statement on definite matrices (\cite{BHAT2}, p. 23, also \cite{BER}).
\begin{lemma}\label{redsym} Let $\mathbf A$ be an arbitrary symmetric matrix, $\mathbf B$ be a symmetric positive (negative) definite matrix. Then they are simultaneously diagonalizable by congruence. Moreover, we can find a nonsingular ${\mathbf S} \in {\mathcal M}^{n \times n}$ such that
${\mathbf S}{\mathbf A}{\mathbf S}^T = {\mathbf \Lambda}_{{\mathbf A}_{\mathbf B}}$ and ${\mathbf S}{\mathbf A}{\mathbf S}^T = {\mathbf I}$ (respectively, $- {\mathbf I}$).
\end{lemma}

For semidefinite matrices, this technique fails (see \cite{NEW}).

Now we are interested in the analogous statement for skew-symmetric matrices. Due to the rich literature on matrix pencils, the following result may be well-known.

\begin{lemma}\label{redskew} Let $\mathbf A$ be a skew-symmetric matrix, $\mathbf B$ be a symmetric positive (negative) definite matrix. Then $\mathbf A$ and $\mathbf B$ are simultaneously quasi-diagonalized by congruence. Moreover, we can find a nonsingular ${\mathbf S} \in {\mathcal M}^{n \times n}$ such that
${\mathbf S}{\mathbf A}{\mathbf S}^T = {\mathbf \Lambda}_{{\mathbf A}_{\mathbf B}}$ and ${\mathbf S}{\mathbf A}{\mathbf S}^T = {\mathbf I}$ (respectively, $- {\mathbf I}$).
\end{lemma}
{\bf Proof.} Consider the case when $\mathbf B$ is symmetric negative definite (the case of positive definiteness is considered analogically). Applying Lemma \ref{propert} to $-{\mathbf B}$, we obtain that ${\mathbf B} = -{\mathbf T}{\mathbf T}^T$ for some nonsingular lower triangular matrix $\mathbf T$. Consider the real matrix ${\mathbf A}_{\mathbf B} = {\mathbf T}^{-1}{\mathbf A}({\mathbf T}^T)^{-1}$. Since congruence transformation preserve skew-symmetry, it is also skew-symmetric, hence normal. By Theorem \ref{canon}, it can be transformed to the quasi-diagonal form by an orthogonal transformation $\mathbf Q$: ${\mathbf A}_{\mathbf B} = {\mathbf Q}{\mathbf \Lambda}_{{\mathbf A}_{\mathbf B}}{\mathbf Q}^T$, where ${\mathbf \Lambda}_{{\mathbf A}_{\mathbf B}}$ is a block-diagonal matrix of Form \ref{qdiag}, ${\mathbf Q}^T = {\mathbf Q}^{-1}$. Then consider the matrix ${\mathbf Y}:={\mathbf Q}^T{\mathbf T}^{-1}$. For $\mathbf A$ and $\mathbf B$, we have
$${\mathbf Y}{\mathbf A}{\mathbf Y}^T = {\mathbf Q}^T{\mathbf T}^{-1}{\mathbf A}({\mathbf T}^{-1})^T{\mathbf Q} =  {\mathbf Q}^T{\mathbf A}_{\mathbf B}{\mathbf Q} = {\mathbf \Lambda}_{{\mathbf A}_{\mathbf B}}$$
and
$${\mathbf Y}{\mathbf B}{\mathbf Y}^T = {\mathbf Q}^T{\mathbf T}^{-1}{\mathbf B}({\mathbf T}^{-1})^T{\mathbf Q} = - {\mathbf Q}^T{\mathbf T}^{-1}{\mathbf T}{\mathbf T}^T({\mathbf T}^{-1})^T{\mathbf Q} = -{\mathbf Q}^T{\mathbf Q} = - {\mathbf I}.$$
 $\square$

Consider more cases, when ${\mathbf A}$ and ${\mathbf B}$ are simultaneously diagonalizable by congruence.
It is well-known (see, for example, \cite{BHAT2}, p. 23) that {\it two Hermitian matrices are simultaneously diagonalizable by unitary congruence if and only if they commute}. However, for the case of arbitrary congruence, this conditions may be reduced. Recall the following criterion of simultaneous diagonalization (see \cite{HHJ}, p. 215, Theorem 2.1 and also \cite{BOC2}, p. 305, Theorem 3, where this result was stated and proved in terms of quadratic forms).
\begin{theorem}\label{ns} Let $\mathbf A$ and $\mathbf B$ be real symmetric matrices with ${\mathbf A}$ being nonsingular, and let ${\mathbf C} = {\mathbf A}^{-1}{\mathbf B}$. There exists a nonsingular matrix ${\mathbf S} \in {\mathcal M}^{n \times n}({\mathbb C})$ such that both ${\mathbf S}{\mathbf A}{\mathbf S}^*$ and ${\mathbf S}{\mathbf B}{\mathbf S}^*$ are diagonal if and only if ${\mathbf C}$ has real eigenvalues and is diagonalizable (i.e. there is a nonsingular ${\mathbf R} \in {\mathcal M}^{n \times n}({\mathbb C})$ such that ${\mathbf R}^{-1}{\mathbf C}{\mathbf R}$ is a real diagonal matrix).
\end{theorem}
The proof of this result implies that the matrix ${\mathbf S}$ can be chosen to have real entries and there is an orthogonal matrix ${\mathbf Q} \in {\mathcal M}^{n \times n}$ such that ${\mathbf S} = {\mathbf Q}{\mathbf R}^T$.

\section{The basic facts about LMI regions}
Given an LMI region $\mathfrak D$, defined by its characteristic function
$$ f_{\mathfrak D}(z) =  {\mathbf L} + {\mathbf M}z+{\mathbf M}^T\overline{z} ,$$
where $\mathbf L$ and $\mathbf M$ are real matrices, such that ${\mathbf L}^T = {\mathbf L}$. The characteristic function $f_{\mathfrak D}$ can also be written in the following form \begin{equation}\label{complex}f_{\mathfrak D}(x + iy) = {\mathbf L} + x({\mathbf M}+{\mathbf M}^T) + iy({\mathbf M}-{\mathbf M}^T) =\end{equation} $$ {\mathbf L} + 2{\rm Sym}({\mathbf M})x + 2{\rm Skew}({\mathbf M})iy,$$ where $x = {\rm Re}(z), y = {\rm Im}(z)$. Here, we call matrices ${\mathbf L}$ and ${\mathbf M}$ {\it generating matrices} of an LMI region ${\mathfrak D}$. The size $m$ of the matrices $\mathbf L$ and $\mathbf M$ we call the {\it order} of a characteristic function. Note, that the characteristic function of an LMI region $\mathfrak D$ is not unique. So it is natural to define {\it the order} of an LMI region $\mathfrak D$ as the minimal possible order of its characteristic functions.
\subsection{Basic properties of LMI regions}
Let us list the following properties of LMI regions, established in \cite{CHG}.
\begin{enumerate}
\item[\rm 1.] {\bf Symmetry.} LMI regions are symmetric with respect to the real axis.
\item[\rm 2.] {\bf Convexity.} LMI regions are convex.
\item[\rm 3.] {\bf Intersection property.} Given two LMI regions ${\mathfrak D}_1$ and ${\mathfrak D}_2$ with the characteristic functions $f_{{\mathfrak D}_1} = {\mathbf L}_1 + z{\mathbf M}_1+\overline{z}{\mathbf M}_1^T$ and $f_{{\mathfrak D}_2} = {\mathbf L}_2 + z{\mathbf M}_2+\overline{z}{\mathbf M}^T_2$, respectively. Then the intersection ${\mathfrak D} = {\mathfrak D}_1\cap{\mathfrak D}_2$ is again an LMI region with the characteristic function $f_{\mathfrak D} = \widetilde{{\mathbf L}} + z\widetilde{{\mathbf M}}+\overline{z}\widetilde{{\mathbf M}}^T$, where $\widetilde{{\mathbf L}} = \diag\{{\mathbf L}_1, {\mathbf L}_2\}$ and $\widetilde{{\mathbf M}} = \diag\{{\mathbf M}_1, {\mathbf M}_2\}$.
\item[\rm 4.] {\bf Density.} LMI regions are dense in the set of convex regions that are symmetric with respect to the real axis.
\end{enumerate}
Properties 1-3 obviously follow from the geometric properties of the class of negative definite matrices (see Section 2), Property 4 is due to the well-known fact from convex analysis that a convex set can be approximated arbitrarily closely by convex polygons.

Now mention some more properties.
\begin{enumerate}
\item[\rm 5.] {\bf Openness.} LMI regions are open. Indeed, if $z \in {\mathfrak D}$, we obtain
$${\mathbf L} + {\mathbf M}z+{\mathbf M}^T\overline{z} = {\mathbf W}(z) \prec 0$$
and the openness of the set of negative definite matrices implies
$${\mathbf L} + {\mathbf M}(z + \Delta z)+{\mathbf M}^T\overline{(z + \Delta z)} = {\mathbf W}(z) - ({\mathbf M}\Delta z + {\mathbf M}^T\overline{\Delta z}) \prec 0$$
for sufficiently small $|\Delta z|$.
\item[\rm 6.] {\bf Invariance under congruence transformations of the characteristic function.} An LMI region remains the same, if we apply to its characteristic function any congruence transformation with a nonsingular matrix $\mathbf B$. I.e. for a nonsingular matrix ${\mathbf B} \in {\mathcal M}^{n \times n}({\mathbb C})$, such that both ${\mathbf B}{\mathbf L}{\mathbf B}^*$ and ${\mathbf B}{\mathbf M}{\mathbf B}^*$ are real, the characteristic functions $f_{\mathfrak D} =  {\mathbf L} + {\mathbf M}z+{\mathbf M}^T\overline{z}$ and $\widetilde{f}_{\mathfrak D} =  {\mathbf B}{\mathbf L}{\mathbf B}^*+ {\mathbf B}{\mathbf M}{\mathbf B}^*z+{\mathbf B}{\mathbf M}^T{\mathbf B}^*\overline{z}$ defines the same LMI region $\mathfrak D$. In particular, when $\mathbf B$ is a nonsingular real matrix, the characteristic functions $f_{\mathfrak D} =  {\mathbf L} + {\mathbf M}z+{\mathbf M}^T\overline{z}$ and $\widetilde{f}_{\mathfrak D} =  {\mathbf B}{\mathbf L}{\mathbf B}^T+ {\mathbf B}{\mathbf M}{\mathbf B}^Tz + {\mathbf B}{\mathbf M}^T{\mathbf B}^T\overline{z}$ defines the same $\mathfrak D$.
\end{enumerate}

\subsection{Topological properties of LMI regions}
Now we are interested in certain topological properties of LMI regions. For this, we recall the following facts from convex analysis (see \cite{WEB}, p. 61, Corollary 2.3.2 and p. 64, Corollary 2.3.9). Here, as usual, we use the notation $\overline{{\mathfrak D}}$ for the closure of $\mathfrak D$, $\partial{\mathfrak D}$ for the boundary of $\mathfrak D$, ${\mathfrak D}^c$ for the completion of $\mathfrak D$ and ${\rm int}({\mathfrak D})$ for the interior of $\mathfrak D$.

\begin{lemma}\label{empyint} A convex set in ${\mathbb R}^2$ has an empty interior if and only if it is a subset of some line in ${\mathbb R}^2$.
\end{lemma}
\begin{lemma}\label{closint} Let ${\mathfrak D} \subseteq {\mathbb R}^2$ be a convex set. Then the following equalities hold:
\begin{equation}{\rm int}({\mathfrak D}) = {\rm int}(\overline{{\mathfrak D}})
\end{equation}
and, when ${\rm int}({\mathfrak D}) \neq \emptyset$,
\begin{equation}\overline{{\mathfrak D}} = \overline{{\rm int}({\mathfrak D})}.
\end{equation}
\end{lemma}

Given an $n \times n$ matrix $\mathbf A$, a positive integer $j$, $1 \leq j \leq n$, and a set of indices $\alpha_j = (i_1, \ \ldots, \ i_j) \subseteq [n]$, $1 \leq i_1 < \ldots < i_j \leq n$, recall, that we use the notation ${\mathbf A}[\alpha_j]$ for the principal submatrix of $\mathbf A$, spanned by the rows and columns with the indices from $\alpha_j$, and the notation ${\mathbf A}(\alpha_j)$ for the principal minor of $\mathbf A$, i.e. the determinant of the corresponding principal submatrix.

\begin{lemma}\label{Clos} Given a nonempty LMI region ${\mathfrak D}$ defined by its characteristic function $ f_{\mathfrak D} =  {\mathbf L} + {\mathbf M}z+{\mathbf M}^T\overline{z}$ of the order $n$. Then
\begin{enumerate}
\item[(i)] $${\mathfrak D} = \bigcap_{j = 1}^n\bigcap_{\alpha_j} P_{\alpha_j} = \bigcap_{j = 1}^n P_{[j]},$$
where $\alpha_j = (i_1, \ \ldots, \ i_j)$, $1 \leq i_1 < \ldots < i_j \leq n$, $[j] = (1, \ \ldots, \ j)$, $P_{\alpha_j}$ is an open polynomial region of the following form:
$$P_{\alpha_j} = \{z \in {\mathbb C}: (-1)^j\det({\mathbf L}[\alpha_j] + {\mathbf M}[\alpha_j]z+{\mathbf M}^T[\alpha_j]\overline{z}) > 0\}.$$
\item[(ii)] $$\overline{{\mathfrak D}} =\{z \in {\mathbb C}: \ {\mathbf L} + {\mathbf M}z+{\mathbf M}^T\overline{z} \preceq 0\} = \bigcap_{j = 1}^n\bigcap_{\alpha_j} \overline{P}_{\alpha_j} \subseteq \bigcap_{j = 1}^n \overline{P}_{[j]},$$
where $\overline{P}_{\alpha_j}$ is a closed polynomial region of the following form:
$$\overline{P}_{\alpha_j} = \{z \in {\mathbb C}: (-1)^j\det({\mathbf L}[\alpha_j] + {\mathbf M}[\alpha_j]z+{\mathbf M}^T[\alpha_j] \overline{z}) \geq 0\}.$$
\item[(iii)] $$\partial{\mathfrak D} \subseteq \bigcup_{i = 1}^n \partial P_{[j]},$$
where $\partial P_{[j]}$ is defined by the polynomial curve $$\partial P_{[j]} = \{z \in {\mathbb C}: \det({\mathbf L}[j] + {\mathbf M}[j]z+{\mathbf M}^T[j] \overline{z}) = 0\}.$$
\item[(iv)] $${\mathfrak D}^c = \bigcup_{i = 1}^n P_{[j]}^c,$$
where $$P^c_{[j]} = \{z \in {\mathbb C}: (-1)^j\det({\mathbf L}[j] + {\mathbf M}[j]z+{\mathbf M}^T[j] \overline{z}) \leq 0\}.$$
\end{enumerate}
\end{lemma}

{\bf Proof.} (i) First, write \eqref{LMI} in the form of the equation
$${\mathbf L} + {\mathbf M}z+{\mathbf M}^T\overline{z} = {\mathbf W}(z), $$
with a negative definite matrix $\mathbf W$. Obviously, the following equality holds for the principal submatrices of ${\mathbf W}$:
$${\mathbf W}[\alpha_j] = {\mathbf L}[\alpha_j] + {\mathbf M}[\alpha_j]z+{\mathbf M}^T[\alpha_j]\overline{z},$$
for any set of indices $\alpha_j = (i_1, \ \ldots, \ i_j)$, $1 \leq i_1 < \ldots < i_j \leq n$. Applying to $-{\mathbf W}(z)$ criterion of positive definiteness (see Lemma \ref{propert}, part 2), we obtain $-{\mathbf W}(z)$ is positive definite if and only if $(-1)^j{\mathbf W}(\alpha_j) = (-1)^j\det({\mathbf L}[\alpha_j] + {\mathbf M}[\alpha_j]z+{\mathbf M}^T[\alpha_j]\overline{z}) > 0$ for all $j = 1, \ \ldots, \ n$
and all $\alpha_j = (i_1, \ \ldots, \ i_j)$, $1 \leq i_1 < \ldots < i_j \leq n$. This obviously implies the first equality in (i). To prove the second equality, we apply Sylvester's criterion (see Lemma \ref{propert}, part 3). Thus $-{\mathbf W}(z)$ is positive definite if and only if $(-1)^j{\mathbf W}([j]) = (-1)^j\det({\mathbf L}[j] + {\mathbf M}[j]z+{\mathbf M}^T[j]\overline{z}) > 0$ for all $j = 1, \ \ldots, \ n$.

(ii) Denote ${\mathfrak D}_{\preceq}:=\{z \in {\mathbb C}: \ {\mathbf L} + {\mathbf M}z+{\mathbf M}^T\overline{z} \preceq 0\}$. First, we show that \begin{equation}\label{ecclos}{\mathfrak D}_{\preceq} = \bigcap_{j = 1}^n\bigcap_{\alpha_j} \overline{P}_{\alpha_j} \subseteq \bigcap_{j = 1}^n \overline{P}_{[j]},\end{equation}
where $\overline{P}_{\alpha_j}$ is a closed polynomial region of the following form:
$$\overline{P}_{\alpha_j} = \{z \in {\mathbb C}: (-1)^j\det({\mathbf L}[\alpha_j] + {\mathbf M}[\alpha_j]z+{\mathbf M}^T[\alpha_j] \overline{z}) \geq 0\}.$$ For this, it is enough to apply the criterion of positive semi-definiteness (see Lemma \ref{propertsemi}, part 1), to the positive semidefinite matrix $-{\mathbf W}(z)$.

 Now we need to prove the equality ${\mathfrak D}_{\preceq} = \overline{{\mathfrak D}}$.
For this, we first show that ${\mathfrak D}_{\preceq}$ is a closed convex subset of $\mathbb C$. Indeed, ${\mathfrak D}_{\preceq} = \bigcap_{j = 1}^n\bigcap_{\alpha_j} \overline{P}_{\alpha_j}$, i.e. is an intersection of closed regions, thus it is closed. Its convexity easily follows from the convexity of the set of negative semidefinite matrices. Secondly, Parts (i) together with Equality \ref{ecclos} and topological identities imply that
$$ {\rm int}({\mathfrak D}_{\preceq}) = {\rm int}\left(\bigcap_{j = 1}^n\bigcap_{\alpha_j} \overline{P}_{\alpha_j}\right) = \bigcap_{j = 1}^n\bigcap_{\alpha_j} {\rm int}(\overline{P}_{\alpha_j}) = \bigcap_{j = 1}^n\bigcap_{\alpha_j} P_{\alpha_j} = {\mathfrak D}.$$
Finally, applying Lemma \ref{closint}, we obtain
$${\mathfrak D}_{\preceq} = \overline{{\rm int}({\mathfrak D}_{\preceq})} = \overline{{\mathfrak D}}.$$

(iii) and (iv) obviously follows from the preceding parts and well-known topological identities.
$\square$

Lemma \ref{Clos} represents an LMI region in the form of an intersection of a finite number of open polynomial regions. It also shows that the set of the form $\{z \in {\mathbb C}: f_{\mathfrak D}(z) = 0\}$, which represents the intersection of some polynomial curves, does not define the boundary of an LMI region $\mathfrak D$. Similarly, the completion ${\mathfrak D}^c$ does not coincide with the set ${\mathfrak D}_{\succeq} = \{z \in {\mathbb C} : f_{\mathfrak D}(z) \succeq 0\}$.

Note that, in general case, the region ${\mathfrak D}_{\preceq} =\{z \in {\mathbb C}: \ {\mathbf L} + {\mathbf M}z+{\mathbf M}^T\overline{z} \preceq 0\}$ may have empty interior. Then, by Lemma \ref{empyint}, it coincides with a closed subset of the real or imaginary axis.

\subsection{Localizations of LMI regions}
Basing on Lemma \ref{Clos}, we obtain the following localization of an LMI region $\mathfrak D$ into an intersection of LMI regions of order 1 and 2.

\begin{lemma}\label{App1} Given a nonempty LMI region ${\mathfrak D}$, defined by its characteristic function $ f_{\mathfrak D} =  {\mathbf L} + {\mathbf M}z+{\mathbf M}^T\overline{z}$ of order $n$, where ${\mathbf L} = \{l_{ij}\}_{i,j = 1}^n$, ${\mathbf M} = \{m_{ij}\}_{i,j = 1}^n$. Then ${\mathfrak D} \subseteq {\mathfrak D}_1$, where ${\mathfrak D}_1$ is a nonempty LMI region, defined by $ f_{{\mathfrak D}_1} =  {\mathbf L}_1 + {\mathbf M}_1z+{\mathbf M}_1^T\overline{z}$, where ${\mathbf L}_1 = {\rm diag}\{l_{11}, \ \ldots, \ l_{nn}\}$ and ${\mathbf M}_1 = {\rm diag}\{m_{11}, \ \ldots, \ m_{nn}\}$ are the diagonal matrices constructed by principal diagonal entries of $\mathbf L$ and $\mathbf M$, respectively.
\end{lemma}

{\bf Proof.} By Lemma \ref{Clos}, part (i), ${\mathfrak D} = \bigcap_{j = 1}^n\bigcap_{\alpha_j} P_{\alpha_j}$. Taking in the first intersection $j = 1$, we obtain the inclusion ${\mathfrak D} \subseteq \bigcap_{i=1}^n P_{i}$, where each $P_i$ is defined by
$$P_i(z) = \{z \in {\mathbb C}: l_{ii} + m_{ii}(z + \overline{z}) < 0\}, \qquad i = 1, \ \ldots, \ n.$$
Now show that $\bigcap_{i=1}^n P_{i} = {\mathfrak D}_1$. Indeed, from the definition of the LMI region ${\mathfrak D}_1$, we get: $z \in {\mathfrak D}_1$ if and only if
$${\mathbf L}_1 + {\mathbf M}_1(z + \overline{z}) = {\mathbf W}_1(z) \prec 0, $$
where ${\mathbf W}_1(z)$ is a diagonal matrix with principal diagonal entries $w_{ii}(z) = l_{ii} + m_{ii}(z + \overline{z})$. By Lemma \ref{propert}, its negative definiteness is equivalent to the negativity of all principal diagonal entries: $w_{ii}= l_{ii}+ m_{ii}(z+\overline{z})< 0$ for $i = 1, \ \ldots, n$. Thus ${\mathfrak D}_1 = \bigcap_{i=1}^n P_{i}$ and ${\mathfrak D} \subseteq \bigcap_{i=1}^n P_{i} = {\mathfrak D}_1$.
 $\square$

The localization of an LMI region in an intersection of shifted halfplanes, given by Lemma \ref{App1}, is obviously too rough. Thus we also consider a localization in an intersection of some second-order regions.

\begin{lemma}\label{App2} Given an LMI region ${\mathfrak D}$ of order $m$ defined by its characteristic function $ f_{\mathfrak D} =  {\mathbf L} + {\mathbf M}z+{\mathbf M}^T\overline{z}$ of order $n$, where ${\mathbf L} = \{l_{ij}\}_{i,j = 1}^n$, ${\mathbf M} = \{m_{ij}\}_{i,j = 1}^n$. Then $${\mathfrak D} \subseteq {\mathfrak D}_1\bigcap\left(\bigcap_{(i,j)} P_{(i,j)}\right),$$ where $1 \leq i < j \leq n$, ${\mathfrak D}_1$ is defined in Lemma \ref{App1}, $P_{(i,j)}$ is a region, bounded by a second-order curve:
\begin{equation}\label{poly} P_{(i,j)} = \{z = x+iy \in {\mathbb C}: a_{11}^{(i,j)}x^2 + a_{22}^{(i,j)}y^2 + 2a_{13}^{(i,j)}x + a_{33}^{(i,j)} > 0\},\end{equation}
with $a_{11}^{(i,j)} = ({\mathbf M} + {\mathbf M}^T)(i,j)$, $a_{22}^{(i,j)} = -({\mathbf M} - {\mathbf M}^T)(i,j)$, $a_{33}^{(i,j)} = {\mathbf L}(i,j)$, $a_{13}^{(i,j)} =
({\mathbf L}\wedge{\mathbf M})(i,j)$, where $({\mathbf L}\wedge{\mathbf M})(i,j)$ denotes so-called mixed minor of matrices $\mathbf L$ and $\mathbf M$ defined as follows:
$$({\mathbf L}\wedge{\mathbf M})(i,j) = \begin{vmatrix}m_{ii} & l_{ij} \\ m_{ji} & l_{jj} \\ \end{vmatrix} + \begin{vmatrix}l_{ii} & m_{ij} \\ l_{ji} & m_{jj} \\ \end{vmatrix}.$$
\end{lemma}

{\bf Proof.} By Lemma \ref{Clos}, part (i), ${\mathfrak D} = \bigcap_{j = 1}^n\bigcap_{\alpha_j} P_{\alpha_j}$. Taking in the first intersection $j = 2$, we obtain the inclusion ${\mathfrak D} \subseteq \bigcap_{(i,j)}^n P_{(i,j)}$, where each $P_i$ is defined by
$$P_{(i,j)}(z) = \{z \in {\mathbb C}: \det({\mathbf L}[i,j] + {\mathbf M}[i,j]z + {\mathbf M}^T[i,j]\overline{z}) > 0\}.$$
Transform the inequality $\det({\mathbf L}[i,j] + {\mathbf M}[i,j]z + {\mathbf M}^T[i,j]\overline{z}) > 0$ into the following form:
$$\det({\mathbf L}[i,j] + ({\mathbf M}+{\mathbf M}^T)[i,j]x + ({\mathbf M}-{\mathbf M}^T)[i,j]iy)>0, $$
i.e. $$\det\left(\begin{pmatrix}l_{ii} & l_{ij} \\
l_{ji} & l_{jj} \end{pmatrix} + \begin{pmatrix}2m_{ii} & m_{ij} + m_{ji} \\  m_{ij} + m_{ji} & 2m_{jj} \end{pmatrix}x + \begin{pmatrix} 0 & m_{ij} - m_{ji} \\ m_{ji} - m_{ij} & 0 \end{pmatrix}iy\right)>0. $$
By expanding the above determinant, we get:
$$\begin{vmatrix}l_{ii} + 2m_{ii}x & l_{ij} + (m_{ij}+m_{ji})x +iy(m_{ij} - m_{ji})\\
l_{ij} + (m_{ij}+m_{ji})x +iy(m_{ji} - m_{ij}) &  l_{jj} + 2m_{jj}x \end{vmatrix}=$$
$$(l_{ii} + 2m_{ii}x)(l_{jj} + 2m_{jj}x) - $$ $$ (l_{ij} + (m_{ij}+m_{ji})x +iy(m_{ij} - m_{ji}))(l_{ij} + (m_{ij}+m_{ji})x +iy(m_{ji} - m_{ij})) =$$
$$x^2(4m_{ii}m_{jj}-(m_{ij}+m_{ji})^2) - y^2(m_{ij} - m_{ji})^2 $$ $$ + 2x(l_{ii}m_{jj} + l_{jj}m_{ii} - l_{ij}(m_{ij}+m_{ji})) + (l_{ii}l_{jj} - l_{ij}^2)=$$
$$x^2\begin{vmatrix} 2m_{ii} & (m_{ij}+m_{ji}) \\ (m_{ij}+m_{ji}) & 2m_{jj} \\ \end{vmatrix} + y^2\begin{vmatrix} 0 & (m_{ij}-m_{ji}) \\ (m_{ij}-m_{ji}) & 0 \\ \end{vmatrix}+ $$ $$ 2x\left(\begin{vmatrix}m_{ii} & l_{ij} \\ m_{ji} & l_{jj} \\ \end{vmatrix} + \begin{vmatrix}l_{ii} & m_{ij} \\ l_{ji} & m_{jj} \\ \end{vmatrix}\right) + \begin{vmatrix}l_{ii} & l_{ij} \\ l_{ji} & l_{jj} \\ \end{vmatrix}=$$
$$x^2({\mathbf M} + {\mathbf M}^T)(i,j) - y^2({\mathbf M} - {\mathbf M}^T)(i,j)  + 2x({\mathbf L}\wedge{\mathbf M})(i,j) + {\mathbf L}(i,j).$$
Applying Lemma \ref{App1}, we complete the proof. $\square$

\begin{lemma}\label{App3} Given an LMI region ${\mathfrak D}$, defined by its characteristic function $ f_{\mathfrak D} =  {\mathbf L} + {\mathbf M}z+{\mathbf M}^T\overline{z}$ of order $n$, where ${\mathbf L} = \{l_{ij}\}_{i,j = 1}^n$, ${\mathbf M} = \{m_{ij}\}_{i,j = 1}^n$. Then ${\mathfrak D} \subseteq {\mathfrak D}_{\alpha_j}$, for any $\alpha_j = (i_1, \ \ldots, \ i_j)$, $1 \leq i_1 < \ldots < i_j \leq n$, and any $j$, $1 \leq j \leq n$, where ${\mathfrak D}_{\alpha_j}$ is an LMI region, defined by its characteristic function $ f_{{\mathfrak D}_{\alpha_j}} =  {\mathbf L}[\alpha_j] + {\mathbf M}[\alpha_j]z+{\mathbf M}^T[\alpha_j]\overline{z}$ of order $j$.
\end{lemma}

{\bf Proof.} From the definition of the LMI region ${\mathfrak D}_{\alpha_j}$, we get: $z \in {\mathfrak D}_{\alpha_j}$ implies
${\mathbf L}[\alpha_j] + {\mathbf M}[\alpha_j]z+{\mathbf M}^T[\alpha_j]\overline{z} = {\mathbf W}[\alpha_j](z) \prec 0,$
where ${\mathbf W}[\alpha_j](z)$ is a principal submatrix of ${\mathbf W}(z)$, spanned by the rows and columns with the indices from $\alpha_j$. Since every principal minor of ${\mathbf W}[\alpha_j]$ is a principal minor of ${\mathbf W}(z)$, the inclusion ${\mathfrak D} \subseteq {\mathfrak D}_{\alpha_j}$ obviously follows from Lemma \ref{Clos}, part (i). $\square$

In fact, for the intersection property we have even a stronger statement.

\begin{lemma}\label{intersect}
An LMI region $\mathfrak D$ can be defined by the characteristic function
\begin{equation}\label{charf} f_{\mathfrak D}(z) =  {\mathbf L} + {\mathbf M}z+{\mathbf M}^T\overline{z},\end{equation} where the matrices $\mathbf L$ and $\mathbf M$ share the same block-diagonal structure: $${\mathbf L}  = {\rm diag}\{{\mathbf L}_{11}, \ \ldots, \ {\mathbf L}_{pp}\},$$
$${\mathbf M}  = {\rm diag}\{{\mathbf M}_{11}, \ \ldots, \ {\mathbf M}_{pp}\},$$
where $1 < p \leq n$, ${\rm dim}({\mathbf L}_{ii}) = {\rm dim}({\mathbf M}_{ii}) = n_i$, $\sum_{i=1}^pn_i = n$, if and only if $${\mathfrak D} = \bigcap_{i=1}^p{\mathfrak D}_i,$$ where each ${\mathfrak D}_i$ is an LMI region, defined by the characteristic function
$$ f_{{\mathfrak D}_i}(z) =  {\mathbf L}_{ii} + {\mathbf M}_{ii}z+{\mathbf M}_{ii}^T\overline{z},$$
of order $n_i$.
\end{lemma}
{\bf Proof.} $\Leftarrow$ This implication is given by the intersection property.

$\Rightarrow$ Let $\mathfrak D$ be defined by characteristic function of Form \eqref{charf}. The inclusion ${\mathfrak D} \subseteq \bigcap_{i=1}^p{\mathfrak D}_i$ follows from Lemma \ref{Clos}. Now let us show the reverse inclusion $\bigcap_{i=1}^p{\mathfrak D}_i \subseteq {\mathfrak D}$. Indeed, consider $ f_{\mathfrak D}(z)$ for $z \in \bigcap_{i=1}^p{\mathfrak D}_i$. We have ${\mathbf L} + {\mathbf M}z+{\mathbf M}^T\overline{z} = {\mathbf W}(z)$, where ${\mathbf W}(z)  = {\rm diag}\{{\mathbf W}_{11}(z), \ \ldots, \ {\mathbf W}_{pp}(z)\}.$ Since all the diagonal blocks ${\mathbf W}_{ii}(z)$ $i = 1, \ \ldots, \ p$, are negative definite, by Sylvester's criterion (Lemma \ref{propert}) so is the matrix ${\mathbf W}(z)$. $\square$

\section{Convex geometry of LMI regions}
The study of robust stability problems requires a deep analysis of the geometric properties of LMI regions. Here, we study LMI regions from the point of view of convex geometry. We consider the questions, when an LMI region $\mathfrak D$ has a conic structure, is it bounded or unbounded, and study such characteristics of unboundedness as the recession cone and the lineality space, through the properties of the generating matrices $\mathbf L$ and $\mathbf M$.
\subsection{Basic definitions and facts} Here, we recall the following definitions and facts from convex analysis (see, for example, \cite{WEB}, \cite{BOYV}).

A nonempty set ${\mathfrak D} \subseteq {\mathbb C}$ is called a {\it cone} if $t z \in {\mathfrak D}$ whenever $z \in {\mathfrak D}$ and $t \geq 0$. A cone $\mathfrak D$ is called {\it solid} if ${\rm int}({\mathfrak D}) \neq \emptyset$. A cone $\mathfrak D$ is called {\it proper} if it is closed, convex, solid and pointed (i.e. ${\mathfrak D}\cap(-{\mathfrak D}) = \{0\}$).

 Given $z \in {\mathbb C}$, a ray $l_0^+$, defined by $l_0^+:=\{t z\}_{t \geq 0}$ is called a {\it direction}. A non-empty convex set ${\mathfrak D} \subseteq {\mathbb C}$ is said {\it to recede in a direction} $l_0^+$ or to have a {\it direction of recession} $l_0^+$ if every half-line of the form $z_0+l_0^+$, where $z_0 \in {\mathfrak D}$, lies in $\mathfrak D$, i.e. ${\mathfrak D} + l_0^+ \subseteq {\mathfrak D}.$ The union of all directions of recession of $\mathfrak D$ together with zero vector is called the {\it recession cone} of $\mathfrak D$ and denoted ${\mathfrak D}_{rc}$.

Consider the following properties of a recession cone (see \cite{WEB}, Theorem 2.5.6).
\begin{lemma}\label{rec}
Let ${\mathfrak D} \subseteq {\mathbb C}$ be a nonempty convex set. Then $${\mathfrak D}_{rc} = \{z \in {\mathbb C}: {\mathfrak D} + z \subseteq {\mathfrak D}\}.$$ Moreover, the recession cone ${\mathfrak D}_{rc}$ is a convex cone, which is closed when $\mathfrak D$ is closed.
\end{lemma}

Later we will use the following criterion of boundedness of a convex set (see \cite{WEB}, p. 74, Theorem 2.5.1).
\begin{lemma}\label{bound}
A non-empty closed convex set ${\mathfrak D}$ is bounded if and only if its recession cone consists of zero vector alone, i.e. ${\mathfrak D}_{rc} = \{0\}$.
\end{lemma}

 Given $z \in {\mathbb C}$, and a line $l_0$, defined by $l_0:=\{\lambda z\}_{\lambda \in {\mathbb R}}$.  A non-empty convex set ${\mathfrak D} \subseteq {\mathbb C}$ is said to be {\it linear in the direction $l_0$} or to have a {\it direction of linearity} $l_0$ if every line, meeting $\mathfrak D$, which has a direction $l_0$, entirely lies in $\mathfrak D$. The union of all directions of linearity together with the zero vector is called the {\it lineality space} of $\mathfrak D$ and denoted $L_{\mathfrak D}$.

The following equality holds (see \cite{WEB}, Theorem 2.5.7.):
$$L_{\mathfrak D} = \{z \in {\mathbb C}: {\mathfrak D} + z = {\mathfrak D}\}.$$

Let us consider the intersection property of recession cones and lineality spaces, which is of importance for studying LMI regions.
\begin{lemma}\label{inter1} Given two convex sets ${\mathfrak D}_1, \ {\mathfrak D}_2 \in {\mathbb C}$. Let ${\mathfrak D} ={\mathfrak D}_1 \bigcap {\mathfrak D}_2 \neq \emptyset$. Then ${\mathfrak D}_{rc} = {\mathfrak D}^1_{rc} \bigcap {\mathfrak D}^2_{rc}$, where ${\mathfrak D}^1_{rc}, \ {\mathfrak D}^2_{rc}$ are the recession cones of ${\mathfrak D}_1$ and ${\mathfrak D}_2$, respectively, and $L_{\mathfrak D} = L_{{\mathfrak D}_1}\bigcap L_{{\mathfrak D}_2}$, where $L_{{\mathfrak D}_1}$ and $L_{{\mathfrak D}_2}$ are the lineality spaces of ${\mathfrak D}_1$ and ${\mathfrak D}_2$, respectively.
\end{lemma}
{\bf Proof}. Since ${\mathfrak D}_1 \bigcap {\mathfrak D}_2$ is convex whenever ${\mathfrak D}_1, \ {\mathfrak D}_2$ are convex, the proof immediately follows from the definitions and Lemma \ref{rec}.

\subsection{Conic LMI regions} Given an LMI region ${\mathfrak D}$, defined by its characteristic function $ f_{\mathfrak D} =  {\mathbf L} + {\mathbf M}z+{\mathbf M}^T\overline{z}$, we call $\mathfrak D$ a {\it uniform region}, if ${\mathbf L} = 0$, i.e $ f_{\mathfrak D} = {\mathbf M}z+{\mathbf M}^T\overline{z}$. Let us prove the following statement, describing which LMI regions are cones in ${\mathbb C}$. Recall, that, as already mentioned in Section 2, a (not necessarily symmetric) matrix $\mathbf M$ is called definite if ${\rm Sym}({\mathbf M})$ is definite.
\begin{theorem}\label{Cone} A nonempty LMI region $\mathfrak D \subset {\mathbb C}$ is a cone in $\mathbb C$ if and only if $\mathfrak D$ is uniform. In this case, $\mathfrak D$ is an open convex cone, symmetric around the negative (positive) direction of the real axis.
\end{theorem}
{\bf Proof.} $\Rightarrow$ First check, that $\mathfrak D$ is a cone, i.e. that $t z \in {\mathfrak D}$ for any $z \in {\mathfrak D}$ and any $t > 0$. Indeed, by definition, $z \in {\mathfrak D}$ if and only if $f_{\mathfrak D}(z) \prec 0$. Thus $f_{\mathfrak D}(t z) = {\mathbf M}t z + {\mathbf M}^T\overline{t z} = t f_{\mathfrak D}(z) \prec 0$ for any $z \in {\mathfrak D}$ and any $t > 0$.
 Now let us show that $\mathfrak D$ is an open convex cone, symmetric around the negative (positive) direction of the real axis.
 Any nonempty LMI region ${\mathfrak D}$ is open and convex (see Properties 2 and 5 of LMI regions). By symmetry (Property 1 of LMI regions), if $z \in {\mathfrak D}$ then $\overline{z} \in {\mathfrak D}$. By convexity (Property 2 of LMI regions), $2{\rm Re}(z) = z + \overline{z} \in {\mathfrak D}$ for any $z \in {\mathfrak D}$. Taking small values of $t$, we obtain $t \in {\mathfrak D}$ can be arbitrarily close to $0$. Thus all the negative (or positive, depends of the sign of ${\rm Re}(z)$) direction of the real axis belongs to $\mathfrak D$, but not all the real line, otherwise it is easy to show that ${\mathfrak D} = {\mathbb C}$. Hence the cone $\mathfrak D$ is symmetric with respect to the negative (or positive) direction of the real axis.

$\Leftarrow$ Let an LMI region $\mathfrak D$ be a cone in ${\mathbb C}$. In this case, as it is shown above, $\mathfrak D$ is an open convex cone in $\mathbb C$ symmetric around the negative (positive) direction of the real axis. Denote its inner angle $2\phi$, $0 < \phi \leq \frac{\pi}{2}$. It is easy to check that $\mathfrak D$ can be defined by either $$f_{\mathfrak D} = \begin{pmatrix} \sin \phi & \cos \phi \\ - \cos\phi & \sin\phi \\ \end{pmatrix}z + \begin{pmatrix} \sin\phi & -\cos\phi \\  \cos\phi & \sin\phi \\ \end{pmatrix}\overline{z},$$
when it is symmetric around the negative direction of the real axis, or
$$f_{\mathfrak D} = \begin{pmatrix} -\sin \phi & \cos \phi \\ - \cos\phi & -\sin\phi \\ \end{pmatrix}z + \begin{pmatrix} -\sin\phi & -\cos\phi \\  \cos\phi & -\sin\phi \\ \end{pmatrix}\overline{z},$$
when it is symmetric around the positive direction.
 $\square$

\begin{corollary}\label{copos} Let ${\mathfrak D}$ be a uniform LMI region, defined by its characteristic function $ f_{{\mathfrak D}} =  {\mathbf M}z+{\mathbf M}^T\overline{z}$. Then ${\mathfrak D} \neq \emptyset$  if and only if $\mathbf M$ is (negative or positive) definite.
\end{corollary}
$\Rightarrow$ Given a uniform LMI region ${\mathfrak D} \neq \emptyset$. Let $z_0 \in {\mathfrak D}$. Then both $\overline{z}_0 \in {\mathfrak D}$ and $z_0+ \overline{z}_0 \in {\mathfrak D}$. Thus
$${\rm Sym}({\mathbf M})(2{\rm Re}(z_0)) = {\rm Sym}({\mathbf M})(z_0+ \overline{z}_0) = {\mathbf W}(z) + {\mathbf W}(\overline{z}) \prec 0$$
and we get that ${\rm Sym}({\mathbf M})$ is either positive or negative definite (according to the sign of ${\rm Re}(z_0)$).

$\Leftarrow$ Given an LMI region $\mathfrak D$ with $\mathbf M$ be negative definite (the case of positive definite $\mathbf M$ is considered analogically). Let us show that the positive direction of the real axis belong to $\mathfrak D$. Indeed, by the substitution $z = x$ to the LMI $$2{\rm Sym}({\mathbf M})x + 2{\rm Skew}({\mathbf M})iy \prec 0,$$ we obtain $2{\rm Sym}({\mathbf M})x \prec 0$, which holds for all $x > 0$.
\begin{corollary} An open cone in $\mathbb C$ with the inner angle $2\phi$ around the positive (negative) direction of the real axis is an LMI region of order $2$ if $0 < \phi < \frac{\pi}{2}$ and of order $1$, if $\phi = \frac{\pi}{2}$.
\end{corollary}

\subsection{Recession cones of LMI regions}
Given an LMI region ${\mathfrak D}$, defined by its characteristic function $ f_{\mathfrak D} =  {\mathbf L} + {\mathbf M}z+{\mathbf M}^T\overline{z}$, consider a uniform LMI region $\overline{{\mathfrak D}}_U$, defined by
  \begin{equation}\label{DU}\overline{{\mathfrak D}}_U = \{z \in {\mathbb C}: {\mathbf M}z+{\mathbf M}^T\overline{z} \preceq 0 \}, \end{equation}
with the characteristic functions $ f_{\overline{{\mathfrak D}}_{U}} =  {\mathbf M}z+{\mathbf M}^T\overline{z}$. By definition, $\overline{{\mathfrak D}}_{U} \neq \emptyset$ (it always contain at least one point $0$). By Lemma \ref{Clos}, $\overline{{\mathfrak D}}_{U}$ is closed, and, by Theorem \ref{Cone}, ${\rm int}(\overline{{\mathfrak D}}_{U})$, if non-empty, is an open convex cone in $\mathbb C$, symmetric around positive or negative direction of the real axis. The properties of the cone of negative semidefinite matrices easily imply, that $\overline{{\mathfrak D}}_{U}$ is a closed convex cone in ${\mathbb C}$, which may have an empty interior or, moreover, consist of only one point $0$.

The following theorem describes the recession cone of an LMI region $\mathfrak D$.
\begin{theorem}\label{resco} Let an LMI region ${\mathfrak D}$ be defined by its characteristic function $ f_{\mathfrak D} =  {\mathbf L} + {\mathbf M}z+{\mathbf M}^T\overline{z}$. Then ${\mathfrak D}_{rc} = \{z \in {\mathbb C}: {\mathbf M}z+{\mathbf M}^T\overline{z} \preceq 0 \}$.
\end{theorem}
{\bf Proof}. $\Rightarrow$ First, let us prove the inclusion ${\mathfrak D}_{rc} \subseteq \overline{{\mathfrak D}}_{U}$. For this, let us take any direction of recession $t z \in {\mathfrak D}_{rc}$, $t > 0$. By definition, we have $z_0 + t z \in {\mathfrak D}$ for any $z_0 \in {\mathfrak D}$ and any $t > 0$. Re-writing the above inclusion in terms of characteristic functions, we obtain $$0 \succ f_{\mathfrak D}(z_0 + tz) =  {\mathbf L} + {\mathbf M}(z_0 + t z)+{\mathbf M}^T(\overline{z_0 + t z}) = $$ $${\mathbf L} + {\mathbf M}z_0 +{\mathbf M}^T\overline{z_0} + t({\mathbf M}z+{\mathbf M}^T\overline{z}) = f_{\mathfrak D}(z_0) + t f_{\overline{{\mathfrak D}}_{U}}(z).$$
Thus we have the following equality for the Hermitian matrices $f_{\mathfrak D}(z_0 + z)$, $f_{\mathfrak D}(z_0)$ and $f_{\overline{{\mathfrak D}}_{U}}(z)$: $$f_{\mathfrak D}(z_0 + z) = t f_{\overline{{\mathfrak D}}_{U}}(z) + f_{\mathfrak D}(z_0).$$ Applying Weyl's inequality \eqref{Wey} (see Lemma \ref{Weyl}), we obtain the following inequality for the eigenvalues:
$$\lambda_1^{\downarrow}(f_{\mathfrak D}(z_0 + z)) \geq t \lambda_1^{\downarrow}(f_{\overline{{\mathfrak D}}_{U}}(z)) + \lambda_n^{\downarrow}(f_{\mathfrak D}(z_0)).$$ Taking into account negative definiteness of $f_{\mathfrak D}(z_0 + z)$, we obtain the inequality $\lambda_1^{\downarrow}(f_{\mathfrak D}(z_0 + z)) < 0$ which implies $$0 > t \lambda_1^{\downarrow}(f_{\overline{{\mathfrak D}}_{U}}(z)) + \lambda_n^{\downarrow}(f_{\mathfrak D}(z_0))$$ for any $t > 0$.
Thus $\lambda_1^{\downarrow}(f_{\overline{{\mathfrak D}}_{U}}(z)) < -\frac{\lambda_n^{\downarrow}(f_{\mathfrak D}(z_0))}{t} \rightarrow 0$ as $t \rightarrow \infty$. So we obtain $$\lambda_1^{\downarrow}(f_{\overline{{\mathfrak D}}_{U}}(z)) \leq 0,$$ which obviously implies $f_{\overline{{\mathfrak D}}_{U}}(z)$ be negative semidefinite and $z \in \overline{{\mathfrak D}}_{U}$.

 $\Leftarrow$ Now let us prove the inclusion $\overline{{\mathfrak D}}_{U} \subseteq {\mathfrak D}_{rc}$. By Lemma \ref{rec}, it is enough to show that ${\mathfrak D} + z \subseteq {\mathfrak D}$ for any $z \in \overline{{\mathfrak D}}_{U}$. Indeed, taking $z_0 + z$, where $z_0 \in {\mathfrak D}$, $z \in \overline{{\mathfrak D}}_{U}$ and considering $f_{{\mathfrak D}}(z_0+z)$, we obtain
$$f_{{\mathfrak D}}(z_0 + z) = {\mathbf L} + {\mathbf M}z_0+{\mathbf M}^T\overline{z}_0 + {\mathbf M}z+{\mathbf M}^T\overline{z} = $$ $$ f_{{\mathfrak D}}(z_0) + f_{\overline{{\mathfrak D}}_{U}}(z).$$
Since $ f_{{\mathfrak D}}(z_0)$ is negative definite and $f_{\overline{{\mathfrak D}}_{U}}(z)$ is negative semidefinite, it follows from the properties of the cone of negative semidefinite matrices (see Subsection 2.2) that their sum $f_{{\mathfrak D}}(z_0 + z)$ is negative definite. $\square$
\begin{corollary} The recession cone ${\mathfrak D}_{rc}$ of an LMI region $\mathfrak D$ is closed and coincides with the recession cone of $\overline{{\mathfrak D}}$.
\end{corollary}
{\bf Proof.} The proof follows from the equality $$\overline{{\mathfrak D}} =\{z \in {\mathbb C}: \ {\mathbf L} + {\mathbf M}z+{\mathbf M}^T\overline{z} \preceq 0\},$$ established in Lemma \ref{Clos}, and the proof of Theorem \ref{resco}. $\square$

Now let us study the cases, when the recession cone ${\mathfrak D}_{rc} \neq \{0\}$.
For the case of positive (negative) semidefinite matrix $\mathbf M$, the following statement holds.

\begin{theorem}\label{semdef} Let an LMI region $\mathfrak D$ be defined by its characteristic function $ f_{{\mathfrak D}} = {\mathbf L} + {\mathbf M}z+{\mathbf M}^T\overline{z}$ with ${\rm Sym}({\mathbf M}) \neq 0$. Then the following statements are equivalent.
\begin{enumerate}
\item[\rm (i)] The matrix ${\mathbf M}$ is negative (respectively, positive) semidefinite.
\item[\rm (ii)] The recession cone ${\mathfrak D}_{rc}$ contains the positive (respectively, negative) direction of the real axis (including 0).
\end{enumerate}
\end{theorem}
{\bf Proof.} $(i)\Rightarrow (ii)$. Let ${\mathbf M}$ be negative semidefinite (the case of positive semidefinite $\mathbf M$ is considered analogically). By Theorem \ref{resco}, we have the equality
\begin{equation}\label{resco1}{\mathfrak D}_{rc} = \{z \in {\mathbb C}: {\mathbf M}z+{\mathbf M}^T\overline{z} \preceq 0 \},\end{equation}
which implies
\begin{equation}\label{resco2}{\mathfrak D}_{rc}\cap{\mathbb R} = \{x \in {\mathbb R}: x({\mathbf M}+{\mathbf M}^T) \preceq 0\}.\end{equation}
The above equality shows that $x \geq 0$ implies $x \in {\mathfrak D}_{rc}$.

$(ii)\Rightarrow (i)$. Let $x \in {\mathfrak D}_{rc}$ for some nonzero $x \in {\mathbb R}$. By Theorem \ref{resco}, we have Equality \eqref{resco2}, which implies $2{\rm Sym}({\mathbf M})x \preceq 0 $. This obviously implies $\mathbf M$ be positive or negative semidefinite, according to the sign of $x$.
$\square$

\subsection{Lineality spaces of LMI regions}
 The structure of the possible lineality spaces of LMI regions can be described by the following statement.
\begin{theorem}\label{lin} Let a nonempty LMI region ${\mathfrak D} \neq {\mathbb C}$ be defined by its characteristic function $ f_{{\mathfrak D}} = {\mathbf L} + {\mathbf M}z+{\mathbf M}^T\overline{z}$. Then its lineality space $L_{\mathfrak D} \neq \{0\}$ if and only if one of the following two cases holds.
\begin{enumerate}
\item[\rm 1.] $\mathbf M$ is symmetric. In this case, $L_{\mathfrak D} = {\mathbb I}$ and $\mathfrak D$ coincides with an open vertical stripe or half-plane in $\mathbb C$, defined by
    $${\mathfrak D} = \{z = x + iy \in {\mathbb C}: x_{min} < x < x_{max}\}, $$ for some values $x_{min}, x_{max} \in \overline{{\mathbb R}}$.
\item[\rm 2.] $\mathbf M$ is skew-symmetric. In this case, $L_{\mathfrak D} = {\mathbb R}$ and $\mathfrak D$ coincides with an open horizontal stripe  in $\mathbb C$, defined by
    $${\mathfrak D} = \{z = x + iy \in {\mathbb C}: |y| < \delta\},$$
    for some $\delta > 0$.
\end{enumerate}
\end{theorem}
{\bf Proof.} For $n = 1$ the statement is obvious. Suppose $n \geq 2$. $\Leftarrow$ Consider Case 1. Since ${\mathbf M} = {\mathbf M}^T$, we obtain that ${\mathfrak D}$ is defined by the following inequality:
\begin{equation}\label{D}{\mathfrak D} = \{z=x + iy \in {\mathbb C}: {\mathbf L} + 2{\mathbf M}x \prec 0\}. \end{equation}
Taking into account that ${\mathfrak D}$ is convex, we observe that $l = \{t iy\}_{t \in {\mathbb R}}$ is obviously a direction of lineality for $\mathfrak D$, and if ${\mathfrak D} \neq {\mathbb C}$ there are no other directions of lineality. Now let us show that $\mathfrak D$ is either a stripe or a half-plane. Indeed, by Lemma \ref{Clos}, $${\mathfrak D} = \bigcap_{j = 1}^n P_{[j]},$$
where $[j] = (1, \ \ldots, \ j)$, $P_{[j]}$ is a region, defined by a polynomial inequality on $x$:
 $$P_{[j]} = \{z \in {\mathbb C}: (-1)^j\det({\mathbf L}[j] + 2x{\rm Sym}({\mathbf M})[j]) > 0\}.$$

The solution of each polynomial inequality $P_{[j]}$ in $x$ is either an empty set or a union of open (finite or infinite) intervals on the real axis. Being viewed in ${\mathbb C}$, it gives a union of open vertical stripes and halfplanes. Thus their intersection $\bigcap_{j = 1}^n P_{[j]},$ if non-empty, also gives an open vertical stripe or halfplane (taking into account convexity). Putting $x_{min} = \inf\{x \in{\mathbb R}: {\mathbf L} + 2{\mathbf M}x \prec 0\}$ and $x_{max} = \sup\{x \in{\mathbb R}: {\mathbf L} + 2{\mathbf M}x \prec 0\}$ (this values may be infinite), we complete the proof.

Consider Case 2. Since ${\mathbf M} + {\mathbf M}^T = 0$, we get $${\mathfrak D} = \{z = x + iy \in {\mathbb C}: {\mathbf L} + iy({\mathbf M} - {\mathbf M}^T)\prec 0\}. $$ In this case, applying Lemma \ref{Clos}, we get the intersection of the polynomial regions $${\mathfrak D} = \bigcap_{j = 1}^n P_{[j]},$$
where $[j] = (1, \ \ldots, \ j)$, $P_{[j]}$ is a region, defined by a polynomial inequality on $x$:
 $$P_{[j]} = \{z \in {\mathbb C}: (-1)^j\det({\mathbf L}[j] + 2iy{\rm Skew}({\mathbf M})[j]) > 0\}.$$
Expanding the determinants, we obtain the solutions of each polynomial inequality with respect to $y$: $$\{iy \in {\mathbb I}: (-1)^j\det({\mathbf L}[j] + 2iy{\rm Skew}({\mathbf M})[j]) > 0\},$$

 Each solution, if non-empty, is a union of (horyzontal) stripes, symmetric with respect to the real axis.  Taking into account convexity, we obtain a ssymmetric with respect to the real axis horyzontal stripe as their intersection. Putting $\delta: = \inf\{y \in {\mathbb R}: (-1)^j\det({\mathbf L}[j] + 2iy{\rm Skew}({\mathbf M})[j]) > 0\}$, we complete the proof. In this case, $L_{\mathfrak D} = {\mathbb R}$ is obviously the direction of lineality.

$\Rightarrow$ Let $L_{\mathfrak D} \neq \{0\}$ for some nonempty LMI region ${\mathfrak D}$. Then, due to the convexity and symmetry of $\mathfrak D$ with respect to the real axis, we have the following two options: $L_{\mathfrak D} = {\mathbb R}$ or $L_{\mathfrak D} = {\mathbb I}$. Assume that the matrix $\mathbf M$ is neither symmetric no skew-symmetric. Then, if $L_{\mathfrak D} = {\mathbb R}$, we have two directions of recession: positive and negative directions of the real axis. Taking into account that
$${\mathfrak D}\cap{\mathbb R} = \{x \in {\mathbb R}: {\mathbf L} + 2{\rm Sym}{\mathbf M}x \prec 0\},$$
and repeating the reasoning of the proof of Theorem \ref{semdef}, we get that ${\rm Sym}({\mathbf M})$ is positive and negative semidefinite at the same time, consequently, ${\rm Sym}({\mathbf M}) = 0$, and we get the contradiction. Now consider the second option $L_{\mathfrak D} = {\mathbb I}$. Applying Lemma \ref{App2}, we obtain the localization of $\mathfrak D$ in the intersection of regions $P_{(i,j)}$, $1 \leq i < j \leq n$, of the following form:
\begin{equation}\label{poly} P_{(i,j)} = \{z = x+iy \in {\mathbb C}: a_{11}^{(i,j)}x^2 + a_{22}^{(i,j)}y^2 + 2a_{13}^{(i,j)}x + a_{33}^{(i,j)} > 0\},\end{equation}
with $a_{11}^{(i,j)} = ({\mathbf M} + {\mathbf M}^T)(i,j)$, $a_{22}^{(i,j)} = -({\mathbf M} - {\mathbf M}^T)(i,j)$.
Let us reduce the corresponding second-order curves to their canonical forms (for the techniques, see, for example, \cite{BOC}). First consider the case $a_{11} \neq 0$. Then the inequality
$$a_{11}^{(i,j)}x^2 + a_{22}^{(i,j)}y^2 + 2a_{13}^{(i,j)}x + a_{33}^{(i,j)} > 0 $$
implies
$$a_{11}^{(i,j)}\left(x^2 + 2\frac{a_{13}^{(i,j)}}{a_{11}^{(i,j)}}x + \frac{(a_{13}^{(i,j)})^2}{(a_{11}^{(i,j)})^2}\right) - \frac{(a_{13}^{(i,j)})^2}{a_{11}^{(i,j)}} + a_{22}^{(i,j)}y^2 + a_{33}^{(i,j)} > 0 $$
and
\begin{equation}\label{region} a_{11}^{(i,j)}\left(x + \frac{a_{13}^{(i,j)}}{a_{11}^{(i,j)}}\right)^2  + a_{22}^{(i,j)}y^2 > \frac{(a_{13}^{(i,j)})^2 - a_{33}^{(i,j)}a_{11}^{(i,j)}}{a_{11}^{(i,j)}} \end{equation}

Then, if $(a_{13}^{(i,j)})^2 - a_{33}^{(i,j)}a_{11}^{(i,j)} = 0$, Inequality \eqref{region} can be transformed to
$$a_{11}^{(i,j)}(x')^2 + a_{22}^{(i,j)}y^2 > 0 $$
which gives an empty region, or an open region, bounded by a pair of intersecting lines with nonzero slope. Since $a_{22}^{(i,j)} < 0$, it does not contain any lines, parallel to $\mathbb I$.

Now denote $\gamma = \frac{(a_{13}^{(i,j)})^2 - a_{33}^{(i,j)}a_{11}^{(i,j)}}{a_{11}^{(i,j)}}$ and assume $\gamma > 0$.
Then Inequality \eqref{region} can be transformed to
$$\frac{a_{11}^{(i,j)}}{\gamma}(x')^2 + \frac{a_{22}^{(i,j)}}{\gamma}y^2 > 1, $$
which gives an empty region or interior part of hyperbola, defined by inequality $(x^{''})^2 - (y^{''})^2 > 1$. This region does not contain any lines, parallel $\mathbb I$.

If $\gamma < 0$, we get $$\frac{a_{11}^{(i,j)}}{\gamma}(x')^2 + \frac{a_{22}^{(i,j)}}{\gamma}y^2 < 1, $$
which gives an interior of an ellipse or the exterior part of hyperbola, defined by inequality $(y^{''})^2 - (x^{''})^2 < 1$. Both of this regions does not contain any lines, parallel $\mathbb I$.

The last case corresponds to $a_{11} = 0$. Then we have
$$ a_{22}^{(i,j)}y^2 + 2a_{13}^{(i,j)}x + a_{33}^{(i,j)} > 0 $$
Since $a_{22}^{(i,j)} < 0$, this region, which is either an interior part of a parabola, or a horizontal stripe, also do not contain any lines parallel to $\mathbb I$.

Since for $a_{22}^{(i,j)} < 0$, all of the regions are bounded by second-order curves on $\mathbb C$, do not contain a line parallel to $\mathbb I$, we have $L_{\mathfrak D} = {\mathbb I}$ implies $-({\mathbf M} - {\mathbf M}^T)(i,j) = a_{22}^{(i,j)} = 0$ for any $1 \leq i < j \leq n$. Due to skew-symmetry of ${\mathbf M} - {\mathbf M}^T$ this is possible if and only if when ${\mathbf M} - {\mathbf M}^T = 0$ and consequently, $\mathbf M$ is symmetric.
 $\square$

\begin{corollary}\label{linsp} Given a non-empty LMI region ${\mathfrak D} \neg {\mathbb C}$, defined by its characteristic function $ f_{{\mathfrak D}} = {\mathbf L} + {\mathbf M}z+{\mathbf M}^T\overline{z}$. Then $L_{\mathfrak D} = \{0\}$ if and only if the matrix $\mathbf M$ is neither symmetric, no skew-symmetric.
\end{corollary}

Now for the case of definite matrix $\mathbf M$, we can show that the recession cone ${\mathfrak D}_{rc}$ is a proper cone in $\mathbb C$.

\begin{theorem}\label{Mdef} Let an LMI region $\mathfrak D$ be defined by its characteristic function $ f_{{\mathfrak D}} = {\mathbf L} + {\mathbf M}z+{\mathbf M}^T\overline{z}$. Then the following statements are equivalent.
\begin{enumerate}
\item[\rm (i)] The matrix ${\mathbf M}$ is non-symmetric negative (respectively, positive) definite.
\item[\rm (ii)] The recession cone ${\mathfrak D}_{rc}$ is a proper cone in $\mathbb C$.
\end{enumerate}
\end{theorem}
{\bf Proof.} $(i)\Rightarrow(ii)$. Let ${\mathbf M}$ be positive (negative) definite. Applying Theorem \ref{resco}, we get ${\mathfrak D}_{rc} = \overline{{\mathfrak D}}_U$, which is, as mentioned above, a closed convex cone in $\mathbb C$. By Corollary \ref{linsp}, $\overline{{\mathfrak D}}_U$ is pointed (i.e. $\overline{{\mathfrak D}}_U \cap(-\overline{{\mathfrak D}}_U ) = \{0\}$) if and only if ${\rm Sym}({\mathbf M}) \neq 0$ and ${\rm Skew}({\mathbf M}) \neq 0$. Now show that $\overline{{\mathfrak D}}_U$ is solid, i.e. ${\rm int}(\overline{{\mathfrak D}}_U)$ is non-empty. By Corollary \ref{copos}, we have, that the set $\{z \in {\mathbb C}: {\mathbf M}z+{\mathbf M}^T\overline{z} \prec 0\}$ is non-empty if and only if $\mathbf M$ is positive (negative) definite.
 Then applying Lemma \ref{Clos}, we obtain, that $$\overline{\{z \in {\mathbb C}: {\mathbf M}z+{\mathbf M}^T\overline{z} \prec 0\}} = \{z \in {\mathbb C}: {\mathbf M}z+{\mathbf M}^T\overline{z} \preceq 0\} = \overline{{\mathfrak D}}_U,$$ which implies
\begin{equation}\label{Intu}{\rm int}(\overline{{\mathfrak D}}_U) = \{z \in {\mathbb C}: {\mathbf M}z+{\mathbf M}^T\overline{z} \prec 0\} \neq \emptyset.\end{equation} Thus $\overline{{\mathfrak D}}_U$ is a solid cone. Together with the properties mentioned above it means, that $\overline{{\mathfrak D}}_U$ is a proper cone.

$(ii)\Rightarrow (i)$. Let for some $z = x+iy \in {\mathbb C}$ we have $z \in {\rm int}({\mathfrak D}_{rc})$. By Theorem \ref{resco}, we have ${\mathfrak D}_{rc} = \overline{{\mathfrak D}}_{U}$. By symmetry and convexity of $\overline{{\mathfrak D}}_{U}$, we conclude that $2x = z+\overline{z} \in {\rm int}(\overline{{\mathfrak D}}_{U})\cap{\mathbb R}$ and Equality \eqref{Intu} implies $x({\mathbf M} + {\mathbf M}^T) \prec 0 $. This obviously implies $\mathbf M$ be positive or negative definite, according to the sign of $x$.
 $\square$
\begin{corollary}\label{realcone} Let an LMI region $\mathfrak D$ be defined by its characteristic function $ f_{{\mathfrak D}} = {\mathbf L} + {\mathbf M}z+{\mathbf M}^T\overline{z}$, with non-symmetric matrix $\mathbf M$. Then the following statements are equivalent.
\begin{enumerate}
\item[\rm (i)] The matrix ${\mathbf M}$ is singular negative (respectively, positive) semidefinite.
\item[\rm (ii)] The recession cone ${\mathfrak D}_{rc} = {\mathbb R}^+$ (respectively, ${\mathfrak D}_{rc} = {\mathbb R}^-$).
\end{enumerate}
\end{corollary}
{\bf Proof.} $\Rightarrow$ Let ${\mathbf M}$ be non-symmetric singular negative semidefinite (the case of positive semidefiniteness can be considered analogically). Then, by Theorem \ref{semdef}, we get ${\mathbb R}^+ \subseteq {\mathfrak D}_{rc}$. Assume that ${\mathbb R}^+ \subset {\mathfrak D}_{rc}$. Then it is easy to see that ${\mathfrak D}_{rc}$ is a proper cone in $\mathbb C$. Applying Theorem \ref{Mdef}, we obtain that $\mathbf M$ is negative definite, hence nonsingular. Contradiction.

$\Leftarrow$ Let ${\mathfrak D}_{rc} = {\mathbb R}^+$ (the case ${\mathfrak D}_{rc} = {\mathbb R}^-$ is considered analogically). By Theorem  \ref{semdef}, $\mathbf M$ is negative semidefinite, hence by Lemma \ref{propertsemi} all its eigenvalues are nonpositive. Assume that $\mathbf M$ is nonsingular. Then by Lemma \ref{propert} it is negative definite. Applying Theorem \ref{Mdef}, we get that ${\mathfrak D}_{rc}$ is proper. Contradiction. $\square$

\subsection{Boundedness of LMI regions} Summarizing the results of previous subsection, we provide the following criterion of the boundness of an LMI region.

\begin{theorem}\label{bounded} A nonempty LMI region $\mathfrak D$, defined by its characteristic function $ f_{{\mathfrak D}} = {\mathbf L} + {\mathbf M}z+{\mathbf M}^T\overline{z}$, is bounded if and only if $\mathbf M$ is indefinite and ${\rm Skew}({\mathbb M}) \neq 0$.
\end{theorem}
{\bf Proof.} $\Rightarrow$ Let $\mathfrak D$ be bounded. Then, by Lemma \ref{bound}, its recession cone ${\mathfrak D}_{rc} = \{0\}$ and, as it immediately follows, its lineality space $L_{\mathfrak D} = \{0\}$. Since $L_{\mathfrak D} = \{0\}$, we get by Theorem \ref{lin}, that $\mathbf M$ is neither symmetric (i.e ${\rm Skew}({\mathbf M}) \neq 0$) no skew-symmetric (i.e. ${\rm Sym}({\mathbf M}) \neq 0$).  Applying Theorem \ref{semdef} to the matrix ${\mathbf M}$ with ${\rm Sym}({\mathbf M}) \neq 0$, we get ${\mathfrak D}_{rc} = \{0\}$ implies $\mathbf M$ be indefinite.

$\Leftarrow$ Let $\mathbf M$ be indefinite and ${\rm Skew}({\mathbf M}) \neq 0$. Assume $\mathfrak D$ is unbounded. Again by Theorem \ref{semdef}, the matrix $\mathbf M$ be indefinite implies ${\mathfrak D}\cap{\mathbb R}$ be bounded. Convexity and symmetry with respect to real axis imply, that the only directions of recession ${\mathfrak D}$ may have are along the imaginary axis, and $L_{\mathfrak D} = OY$. Then by Theorem \ref{lin} ${\mathfrak D}$ is an open stripe or halfplane and ${\rm Skew}({\mathbf M}) = 0$. We came to the contradiction. $\square$

\section{Canonical forms of matrices and complete description of LMI regions}
Using the results of Subsection 2.4, here we consider the cases, when an LMI region $\mathfrak D$ coincides with an intersection of certain regions, bounded by first- and second-order curves. The results of this section allows us to find the lowest order characteristic functions for certain LMI regions.
\subsection{Case of commuting matrices}
 First, let us consider the following simple special case.

\begin{theorem}\label{comp} Let a nonempty LMI region $\mathfrak D$ be defined by its characteristic function $f_{\mathfrak D} = {\mathbf L} + {\mathbf M}z + {\mathbf M}^T\overline{z}$ with $\mathbf M$ be normal and ${\mathbf L}{\mathbf M} = {\mathbf M}{\mathbf L}$. Let $\sigma({\mathbf L}) = \{\lambda_i({\mathbf L})\}_{i = 1}^n$ and, respectively, $\sigma({\mathbf M}) = \{\lambda_i({\mathbf M})\}_{i = 1}^n$.
Then the LMI region $\mathfrak D$ coincides with the intersection of the following types of regions, defined by the eigenvalues of matrices $\mathbf L$ and $\mathbf M$:
\begin{enumerate}
\item[\rm 1.] Shifted halfplanes ${\mathfrak D}_1^i$ of the form $${\mathfrak D}_1^i = \{z = x + iy \in {\mathbb C}: \lambda_i({\mathbf L}) + 2 {\rm Re}(\lambda_i({\mathbf M})) x < 0\},$$ where the indices $i \in [n]$ are such that ${\rm Re}(\lambda_i({\mathbf M})) \neq 0$.
\item[\rm 2.] Shifted cones ${\mathfrak D}_2^i$ with the vertex at the point $(-\frac{\lambda_i({\mathbf L})}{2{\rm Re}(\lambda_i({\mathbf M}))},0)$ and the inner angle $2\theta$, $\theta = |\frac{\pi}{2} - {\rm arg}(\lambda_i(\mathbf M))|$ around the negative direction of the real axis if $0< {\rm arg}(\lambda_i(\mathbf M)) < \frac{\pi}{2}$ and around the positive direction of the real axis if $\frac{\pi}{2}< {\rm arg}(\lambda_i(\mathbf M)) < \pi$. Here $i \in [n]$ are such that ${\rm Re}(\lambda) \neq 0$ and ${\rm Im}(\lambda) \neq 0$, one cone ${\mathfrak D}_2^i$ corresponds to a pair of the complex conjugate eigenvalues of $\mathbf M$.
\item[\rm 3.] Horizontal stripes ${\mathfrak D}^i_{3}$ of the form $${\mathfrak D}^i_{3} = \{z = x+iy \in {\mathbb C}:  -\frac{\lambda_i({\mathbf L})}{2{\rm Im}(\lambda_i({\mathbf M}))} < y < \frac{\lambda_i({\mathbf L})}{2{\rm Im}(\lambda_i({\mathbf M}))}\},$$ where $i \in [n]$ are such that ${\rm Re}(\lambda_i({\mathbf M})) = 0$.
\end{enumerate}
\end{theorem}
{\bf Proof.} By Corollary \ref{commute3}, $\mathbf L$, $\mathbf M$ and ${\mathbf M}^T$ can be simultaneously quasi-diagonalized by some orthogonal congruence $\mathbf Q$:
 $${\mathbf L} = {\mathbf Q}{\mathbf \Lambda}_{\mathbf L}{\mathbf Q}^T,$$
 where ${\mathbf \Lambda}_{\mathbf L} = {\rm diag}\{\lambda_1({\mathbf L}), \ \ldots, \ \lambda_n({\mathbf L})\}$;
 and
   $${\mathbf M} = {\mathbf Q}{\mathbf \Lambda}_{\mathbf M}{\mathbf Q}^T,$$
where $$ {\mathbf \Lambda}_{\mathbf M} ={\rm diag}\{\begin{pmatrix}\mu_1 & \nu_1 \\- \nu_1 & \mu_1 \end{pmatrix}, \ldots, \begin{pmatrix}\mu_k & \nu_k \\- \nu_k & \mu_k \end{pmatrix}, \ \lambda_{2k+1}({\mathbf M}), \ \ldots, \lambda_n({\mathbf M})\},$$
where $\lambda_{2j-1}({\mathbf M}) = \mu_j+i\nu_j$, $\lambda_{2j}({\mathbf M}) = \mu_j-i\nu_j$, $j = 1, \ldots, k$ are non-real eigenvalues of $\mathbf M$, $\lambda_{2k+1}({\mathbf M}), \ \ldots, \lambda_n({\mathbf M})$ are real eigenvalues of $\mathbf M$.

Moreover, the matrices ${\mathbf \Lambda}_{\mathbf L}$ and ${\mathbf \Lambda}_{\mathbf M}$ are connected by Corollary \ref{commute4}.
Then, by Property 6 of LMI regions, $\mathfrak D$ can be defined by the characteristic function of the following form:
$$\widetilde{f}_{\mathfrak D} = {\mathbf \Lambda}_{\mathbf L} + {\mathbf \Lambda}_{\mathbf M}z + {\mathbf \Lambda}_{\mathbf M}^T\overline{z}. $$
 Applying Lemma \ref{App1}, we obtain $\mathfrak D$ is contained in the LMI region ${\mathfrak D}_1$, defined by $ f_{{\mathfrak D}_1} =  {\mathbf \Lambda}_{\mathbf L} + 2\widetilde{{\mathbf \Lambda}}_{\mathbf M}x$, where $$\widetilde{{\mathbf \Lambda}}_{\mathbf M} = {\rm diag}\{\mu_1, \ \mu_1, \ldots, \ \mu_k, \mu_k, \  \ \lambda_{2k+1}({\mathbf M}), \ \ldots, \lambda_n({\mathbf M})\}$$ is the diagonal matrices constructed by principal diagonal entries of ${\mathbf \Lambda}_{\mathbf M}$, where each $\mu_i = {\rm Re}(\lambda_{2i-1}({\mathbf M})) =  {\rm Re}(\lambda_{2i}({\mathbf M}))$, $i = 1, \ \ldots, \ k$. Thus ${\mathfrak D}_1$ coinsides the intersection of half-planes of the form $${\mathfrak D}_1^i = \{z = x + iy \in {\mathbb C}: \lambda_i({\mathbf L}) + 2 {\rm Re}(\lambda_i({\mathbf M})) x < 0\},$$
 $1 \leq i \leq n$ and ${\rm Re}(\lambda_i({\mathbf M})) \neq 0$.

 Both matrices ${\mathbf \Lambda}_{\mathbf L}$ and ${\mathbf \Lambda}_{\mathbf M}$ have a block-diagonal structure with the size of the blocks $\leq 2$. Thus, by Lemma \ref{App3}, ${\mathfrak D}$ coincides with the intersection of regions ${\mathfrak D}^j_i$ of order 1 and 2, defined by $f_{{\mathfrak D}_{i}} = ({\mathbf \Lambda}_{\mathbf L})_{ii} + ({\mathbf \Lambda}_{\mathbf M})_{ii}z + ({\mathbf \Lambda}^T_{\mathbf M})_{ii}\overline{z}$, where $({\mathbf \Lambda}_{\mathbf L})_{ii}$ and $({\mathbf \Lambda}_{\mathbf M})_{ii}$ are diagonal blocks of ${\mathbf \Lambda}_{\mathbf L}$ and ${\mathbf \Lambda}_{\mathbf M}$, respectively. Consider the following two cases.

{\bf Case I.}$({\mathbf \Lambda}_{\mathbf M})_{ii}$ is a real $1 \times 1$ matrix, which corresponds to a real eigenvalue $\lambda_i({\mathbf M})$. Then the corresponding LMI region is of the form
 $${\mathfrak D}^i_{1} = \{z = x + iy \in {\mathbb C}: \lambda_i({\mathbf L}) + 2 \lambda_i({\mathbf M}) x < 0\}.$$
If $\lambda_i({\mathbf M}) \neq 0$, it obviously represents a shifted half-plane $${\mathfrak D}^{i}_1 = \{z = x + iy \in {\mathbb C}: x < - \frac{\lambda_i({\mathbf L})}{2 \lambda_i({\mathbf M})}\}.$$ The number of such half-planes obviously coincides with the number of real nonzero eigenvalues $\lambda_i({\mathbf M})$. In case of $\lambda_i({\mathbf M}) = 0$, ${\mathfrak D}_1^{i} = {\mathbb C}$ if $\lambda_i({\mathbf L})< 0$ and ${\mathfrak D}_1^{i} = \emptyset$ if $\lambda_i({\mathbf L})> 0$.

 {\bf Case II.} $({\mathbf \Lambda}_{\mathbf M})_{ii}$ is a real $2 \times 2$ matrix of the form:
 $$({\mathbf \Lambda}_{\mathbf M})_{ii} = \begin{pmatrix}{\rm Re}(\lambda_i({\mathbf M})) & {\rm Im}(\lambda_i({\mathbf M})) \\ - {\rm Im}(\lambda_i({\mathbf M})) & {\rm Re}(\lambda_i({\mathbf M})) \\ \end{pmatrix} = \rho\begin{pmatrix}\cos\varphi & \sin\varphi \\ -\sin\varphi & \cos\varphi \\ \end{pmatrix},$$
 which corresponds to the pair $\lambda_{i, i+1}({\mathbf M})) ={\rm Re}(\lambda_i({\mathbf M})) \pm {\rm Im}(\lambda_i({\mathbf M}))= \rho(\cos\varphi \pm i\sin\varphi)$ of the complex conjugate eigenvalues of ${\mathbf M}$. Note, that in this case, ${\rm Im}(\lambda_i({\mathbf M})) \neq 0$. Taking into account Corollary \ref{commute4}, we obtain that the corresponding LMI region is of the form
 $$\tiny {\mathfrak D}^{i}_2 = \left\{z \in {\mathbb C}: \begin{pmatrix}\lambda_i({\mathbf L}) & 0 \\ 0 & \lambda_i({\mathbf L}) \\ \end{pmatrix} + \begin{pmatrix}{\rm Re}(\lambda_i({\mathbf M})) & {\rm Im}(\lambda_i({\mathbf M})) \\ - {\rm Im}(\lambda_i({\mathbf M})) & {\rm Re}(\lambda_i({\mathbf M})) \\ \end{pmatrix}z + \begin{pmatrix}{\rm Re}(\lambda_i({\mathbf M})) & -{\rm Im}(\lambda_i({\mathbf M})) \\  {\rm Im}(\lambda_i({\mathbf M})) & {\rm Re}(\lambda_i({\mathbf M})) \\ \end{pmatrix}\overline{z}\prec 0\right\}.$$
 Using Lemma \ref{App2} and Lemma \ref{intersect}, we get that
 $${\mathfrak D}^{i}_2 = P_1^i(z) \bigcap P_{(1,2)}^i(z),$$
Since all the regions of the first order are already considered above, we consider $P_{(1,2)}^i(z)$. By Formula \eqref{poly},
\begin{equation}\label{poly1} P^i_{(1,2)}(z) = \{z = x+iy \in {\mathbb C}: a_{11}x^2 + a_{22}y^2 + 2a_{13}x + a_{33} > 0\},\end{equation}
where $a_{11} = \det(({\mathbf \Lambda}_{\mathbf M})_{ii} + ({\mathbf \Lambda}_{\mathbf M})_{ii}^T)$, $a_{22} = -\det(({\mathbf \Lambda}_{\mathbf M})_{ii} - ({\mathbf \Lambda}_{\mathbf M})_{ii}^T)$, $a_{13} = ({\mathbf \Lambda}_{\mathbf M})_{ii}\wedge({\mathbf \Lambda}_{\mathbf L})_{ii}$, $a_{33} = \det(({\mathbf \Lambda}_{\mathbf L})_{ii})$. Calculating the coefficients, we obtain
$$a_{11} = \begin{vmatrix}2{\rm Re}(\lambda_i({\mathbf M})) & 0 \\ 0 & 2{\rm Re}(\lambda_i({\mathbf M})) \\ \end{vmatrix} = 4 ({\rm Re}(\lambda_i({\mathbf M})))^2;$$
$$a_{22} = -\begin{vmatrix}0 & 2{\rm Im}(\lambda_i({\mathbf M})) \\ - 2{\rm Im}(\lambda_i({\mathbf M})) & 0 \\ \end{vmatrix} = - 4({\rm Im}(\lambda_i({\mathbf M})))^2;$$
$$a_{13} = \begin{vmatrix}{\rm Re}(\lambda_i({\mathbf M})) & 0 \\ - {\rm Im}(\lambda_i({\mathbf M})) & \lambda_i({\mathbf L}) \\ \end{vmatrix} + \begin{vmatrix}\lambda_i({\mathbf L}) & {\rm Im}(\lambda_i({\mathbf M})) \\ 0 & {\rm Re}(\lambda_i({\mathbf M})) \\ \end{vmatrix} =  2{\rm Re}(\lambda_i)\lambda_i({\mathbf L}).$$
After substitution, Inequality \eqref{poly1} get the form
\begin{equation}\label{poly2} P^i_{(1,2)}(z) =\end{equation}
$$ \{ z = x+iy \in {\mathbb C}:  4 ({\rm Re}(\lambda_i({\mathbf M})))^2 x^2 - 4({\rm Im}(\lambda_i({\mathbf M})))^2y^2 + 4{\rm Re}(\lambda_i({\mathbf M}))\lambda_i({\mathbf L})x + \lambda_i^2({\mathbf L}) > 0\},$$
Calculating the second-order and third-order determinants
$$\delta = \begin{vmatrix}a_{11} & a_{12} \\ a_{21} & a_{22} \\ \end{vmatrix} =$$ $$ \begin{vmatrix}4 ({\rm Re}(\lambda_i({\mathbf M})))^2 & 0 \\ 0 & - 4({\rm Im}(\lambda_i({\mathbf M})))^2 \\ \end{vmatrix} = -16({\rm Re}(\lambda_i({\mathbf M})))^2({\rm Im}(\lambda_i({\mathbf M})))^2;$$
$$\Delta = \begin{vmatrix}a_{11} & a_{12} & a_{13}\\ a_{21} & a_{22} & a_{23} \\ a_{31} & a_{32} & a_{33} \\\end{vmatrix} =$$ $$ \begin{vmatrix}4 ({\rm Re}(\lambda_i({\mathbf M})))^2 & 0 &  2{\rm Re}(\lambda_i({\mathbf M}))\lambda_i({\mathbf L}) \\ 0 & - 4({\rm Im}(\lambda_i({\mathbf M})))^2 & 0 \\  2{\rm Re}(\lambda_i({\mathbf M}))\lambda_i({\mathbf L}) & 0 &  \lambda^2_i({\mathbf L}) \\ \end{vmatrix} = $$ $$ - 16({\rm Re}(\lambda_i))^2({\rm Im}(\lambda_i))^2\lambda^2_i({\mathbf L}) + 16({\rm Re}(\lambda_i))^2({\rm Im}(\lambda_i))^2\lambda^2_i({\mathbf L}) =0;$$
and applying well-known results on the classification of second-order curves (see, for example, \cite{BOC}, p. 182), we obtain the following cases, that correspond to Parts 2 and 3 of the statement of the theorem, respectively.

{\bf Case IIa}. ${\rm Re}(\lambda_i({\mathbf M})) \neq 0$. In this case, as we have already assumed ${\rm Im}(\lambda_i({\mathbf M})) \neq 0$, we get $\delta < 0$ and the corresponding curve is of hyperbolic type. Since $\Delta = 0$, the corresponding curve is a pair of intersecting lines.

Transforming Inequality \eqref{poly2}, we obtain the following boundary conditions:
\begin{equation}\label{poly3}(2 {\rm Re}(\lambda_i({\mathbf M}))x + \lambda_i({\mathbf L}))^2 - 4({\rm Im}(\lambda_i({\mathbf M})))^2y^2 = 0. \end{equation}
After a shift along the real axis $x_1 = x + \frac{\lambda_i({\mathbf L})}{2 {\rm Re}(\lambda_i({\mathbf M}))}$, which put the intersection point of the lines to zero, we transform Equation \eqref{poly3} into its canonical form:
$$\frac{x_1^2}{4({\rm Im}(\lambda_i({\mathbf M})))^2 } - \frac{y^2}{4({\rm Re}(\lambda_i({\mathbf M})))^2} = 0 $$

Well-known formulae give the equality for the angle between the line and the positive direction of the real axis: $\tan(\theta) = - \frac{{\rm Re}(\lambda_i({\mathbf M}))}{{\rm Im}(\lambda_i({\mathbf M}))}= - \frac{\cos \varphi}{\sin \varphi},$ where $\lambda_i = \rho e^{i\varphi}$. If $0 < \varphi < \frac{\pi}{2}$, by using trigonometric formulae, we obtain that $\theta$, $0 < \theta < \frac{\pi}{2}$ is an angle around the negative direction of the real axis, satisfying $\theta = \frac{\pi}{2} - \varphi$. In case $\frac{\pi}{2} < \varphi < \pi$, we consider $\theta$ to be the angle around the positive direction and $\theta = \varphi - \frac{\pi}{2}$.

 {\bf Case IIb}. ${\rm Re}(\lambda_i) = 0$. In this case, $\delta = 0$ and $\Delta = 0$. By \cite{BOC}, we get that the corresponding curve is a pair of parallel lines, and the region $P^i_{(1,2)}(z)$ is a horizontal stripe. The boundary conditions may be easily derived directly from the inequality for the characteristic function:
 $$f_{{\mathfrak D}^3_{i}}(x+iy) = \begin{pmatrix} \lambda_i({\mathbf L}) & 0 \\ 0 & \lambda_i({\mathbf L}) \\ \end{pmatrix} + \begin{pmatrix} 0 & 2{\rm Im}(\lambda_i({\mathbf M})) \\ -2{\rm Im}(\lambda_i({\mathbf M})) & 0 \\ \end{pmatrix}iy \prec 0,$$
 which is equivalent to
 $$ \lambda_i^2({\mathbf L}) -  4({\rm Im}(\lambda_i({\mathbf M})))^2 = (\lambda_i({\mathbf L}) - 2{\rm Im}(\lambda_i({\mathbf M}))y)(\lambda_i({\mathbf L}) + 2{\rm Im}(\lambda_i({\mathbf M}))y)>0 .$$
 Taking into account that ${\rm Im}(\lambda_i)$ is assumed to be positive (and the corresponding pair of complex conjugate eigenvalues is defined by $\lambda_{i,i+1}({\mathbf M}) = \pm {\rm Im}(\lambda_i({\mathbf M}))$) we obtain the conditions $${\mathfrak D}^{i}_3 = \{z = x+iy \in {\mathbb C}:  -\frac{\lambda_i({\mathbf L})}{2{\rm Im}(\lambda_i({\mathbf M}))} < y < \frac{\lambda_i({\mathbf L})}{2{\rm Im}(\lambda_i({\mathbf M}))}\}.$$
 This gives exactly Case 3 of the theorem statement.
 $\square$

\subsection{Case of simultaneously quasi-diagonalizable matrices}

Since the condition of commutativity of $\mathbf L$ and $\mathbf M$ is sufficient, but not necessarily for simultaneous reduction by congruence to a diagonal and quasi-diagonal forms, respectively, in the following statement, we assume simultaneous reduction as is. Note, that under this assumption, the eigenvalues of the corresponding forms may not coincide with the eigenvalues of the initial matrices $\mathbf L$ and $\mathbf M$.

\begin{theorem}\label{comp1} Let a nonempty LMI region $\mathfrak D$ be defined by its characteristic function $f_{\mathfrak D} = {\mathbf L} + {\mathbf M}z + {\mathbf M}^T\overline{z}$ with $\mathbf M$ be normal. Let $\mathbf M$ and $\mathbf L$ be simultaneously reduced by congruence to a diagonal and quasi-diagonal forms, respectively.
Then the LMI region $\mathfrak D$ coincides with the intersection of the following four types of regions:
\begin{enumerate}
\item[\rm 1.] shifted halfplanes;
\item[\rm 2.] shifted cones around the positive or negative direction of the real axis;
\item[\rm 3.] horizontal stripes symmetric with respect to the real axis;
\item[\rm 4.] hyperbolas.
\end{enumerate}
\end{theorem}
{\bf Proof.} Let $\mathbf S$ be a nonsingular matrix such that
$${\mathbf S}{\mathbf L}{\mathbf S}^T = {\mathbf \Lambda}_{\mathbf L}^{\mathbf S} = {\rm diag}\{\widetilde{\lambda}_1, \ \ldots, \ \widetilde{\lambda}_n\}, $$
$${\mathbf S}{\mathbf M}{\mathbf S}^T = {\mathbf \Lambda}_{\mathbf M}^{\mathbf S} ={\rm diag}\{\begin{pmatrix}\widetilde{\mu}_1 & \widetilde{\nu}_1 \\- \widetilde{\nu}_1 & \widetilde{\mu}_1 \end{pmatrix}, \ldots, \begin{pmatrix}\widetilde{\mu}_k & \widetilde{\nu}_k \\- \widetilde{\nu}_k & \widetilde{\mu}_k \end{pmatrix}, \ \kappa_{2k+1}, \ \ldots, \kappa_n\}.$$

Consider the diagonal blocks ${\mathbf \Lambda}_{\mathbf M}^{\mathbf S}[i,i+1]$ of ${\mathbf \Lambda}_{\mathbf M}^{\mathbf S}$, where $i = 1, \ 3, \ldots, \ 2k-1$, and the corresponding diagonal blocks ${\mathbf \Lambda}_{\mathbf L}^{\mathbf S}[i,i+1]$ of ${\mathbf \Lambda}_{\mathbf L}^{\mathbf S}$.
  Then the corresponding LMI region is of the form
 $${\mathfrak D}_{i} = \left\{z \in {\mathbb C}: \pm\begin{pmatrix}\widetilde{\lambda}_i& 0 \\ 0 & \widetilde{\lambda}_{i+1}  \\ \end{pmatrix} + \begin{pmatrix}\widetilde{\mu}_i & \widetilde{\nu}_i \\ - \widetilde{\nu}_i & \widetilde{\mu}_i \\ \end{pmatrix}z + \begin{pmatrix}\widetilde{\mu}_i & -\widetilde{\nu}_i \\  \widetilde{\nu}_i & \widetilde{\mu}_i \\ \end{pmatrix}\overline{z}\right\}.$$
 If $\widetilde{\lambda}_i = \widetilde{\lambda}_{i+1}$ we just repeat the reasoning of the previous proof and obtain one of the cases 1, 2 or 3.

Now consider the case $\widetilde{\lambda}_i \neq \widetilde{\lambda}_{i+1}$. Applying Lemma \ref{App2} and Lemma \ref{intersect} to ${\mathfrak D}_{i}$, we get that
 $${\mathfrak D}_{i} = P_1^i(z) \bigcap P_1^{i+1}(z) \bigcap P^i_{(1,2)}(z),$$
 where $P_1^i(z)$ is the first-order region, defined by $$P_1^i(z) = \{z = x + iy \in {\mathbb C}: \widetilde{\lambda}_i+ 2 \widetilde{\mu}_i x < 0\}.$$
If $\widetilde{\mu}_i \neq 0$, it obviously represents a shifted half-plane $$P_1^i(z) = \{z = x + iy \in {\mathbb C}: x < -\frac{\widetilde{\lambda}_i}{\widetilde{\mu}_i}\}.$$ The case $\widetilde{\mu}_i = 0$ gives either the whole $\mathbb C$ or an empty region. By analogy, $$P_1^{i+1}(z) = \{z = x + iy \in {\mathbb C}: x < -\frac{\widetilde{\lambda}_{i+1}}{\widetilde{\mu}_i}\}.$$

In its turn, by Formula \eqref{poly},
\begin{equation}\label{poly1} P^i_{(1,2)}(z) = \{z = x+iy \in {\mathbb C}: a_{11}x^2 + a_{22}y^2 + 2a_{13}x + a_{33} > 0\},\end{equation}
where $a_{11} = \det({\mathbf \Lambda}_{\mathbf M}^{\mathbf S}[i,i+1] + ({\mathbf \Lambda}_{\mathbf M}^{\mathbf S}[i,i+1])^T)$, $a_{22} = -\det({\mathbf \Lambda}_{\mathbf M}^{\mathbf S}[i,i+1] - ({\mathbf \Lambda}_{\mathbf M}^{\mathbf S}[i,i+1])^T)$, $a_{13} = ({\mathbf \Lambda}_{\mathbf M}^{\mathbf S}[i,i+1])\wedge ({\mathbf \Lambda}_{\mathbf L}^{\mathbf S}[i,i+1])$, $a_{33} = {\mathbf \Lambda}_{\mathbf L}^{\mathbf S}(i,i+1) = \widetilde{\lambda}_{i}\widetilde{\lambda}_{i+1}$. Calculating and substituting the coefficients, we get
\begin{equation}\label{poly2} P^i_{(1,2)}(z) = \end{equation}
$$\{z = x+iy \in {\mathbb C}:  4 (\widetilde{\mu}_i)^2 x^2 - 4(\widetilde{\nu}_i)^2y^2 + 2\widetilde{\mu}_i(\widetilde{\lambda}_{i} + \widetilde{\lambda}_{i+1}) x + \widetilde{\lambda}_{i}\widetilde{\lambda}_{i+1} > 0\},$$

Assuming $\widetilde{\mu}_i \neq 0$, $\widetilde{\nu}_i \neq 0$ and $\widetilde{\lambda}_{i} \neq \widetilde{\lambda}_{i+1}$, we estimate the second-order and third-order determinants
$$\delta = \begin{vmatrix}a_{11} & a_{12} \\ a_{21} & a_{22} \\ \end{vmatrix} = \begin{vmatrix}4 (\widetilde{\mu}_i)^2 & 0 \\ 0 & - 4(\widetilde{\nu}_i)^2 \\ \end{vmatrix} = -16(\widetilde{\mu}_i)^2(\widetilde{\nu}_i)^2 < 0;$$
$$\Delta = \begin{vmatrix}a_{11} & a_{12} & a_{13}\\ a_{21} & a_{22} & a_{23} \\ a_{31} & a_{32} & a_{33} \\\end{vmatrix} = \begin{vmatrix}4 (\widetilde{\mu}_i)^2 & 0 & \widetilde{\mu}_i(\widetilde{\lambda}_{i}+\widetilde{\lambda}_{i+1}) \\ 0 & - 4(\widetilde{\nu}_i)^2 & 0 \\ \widetilde{\mu}_i(\widetilde{\lambda}_{i}+\widetilde{\lambda}_{i+1}) & 0 &  \widetilde{\lambda}_{i}\widetilde{\lambda}_{i+1} \\ \end{vmatrix} = $$ $$ -16\widetilde{\lambda}_{i}\widetilde{\lambda}_{i+1}(\widetilde{\mu}_i)^2(\widetilde{\nu}_i)^2 + 4(\widetilde{\mu}_i)^2(\widetilde{\lambda}_{i}+\widetilde{\lambda}_{i+1})^2(\widetilde{\nu}_i)^2 =$$
$$ 4(\widetilde{\mu}_i)^2(\widetilde{\lambda}_{i}-\widetilde{\lambda}_{i+1})^2(\widetilde{\nu}_i)^2 > 0,$$

Applying well-known results on the classification of second-order curves (see, for example, \cite{BOC}, p. 182), we get, that the corresponding line is a hyperbola. $\square$

\subsection{Arbitrary case}
In the most general case, we have the following statement.

\begin{theorem}\label{arbit} Let a nonempty LMI region $\mathfrak D$ be defined by its characteristic function $f_{\mathfrak D} = {\mathbf L} + {\mathbf M}z + {\mathbf M}^T\overline{z}$ of order $n\geq 2$ with ${\rm In}({\rm Sym}({\mathbf M})) = (i_+, i_-, i_0)$.
Then \begin{equation}\label{2ndinter}{\mathfrak D} \subseteq \bigcap_{(i,j)}P_{(i,j)},\end{equation}
where $1 \leq i < j \leq n$, $P_{(i,j)}$ is a polynomial region defined by \eqref{poly}.
The intersection $\bigcap_{(i,j)}P_{(i,j)}$ contains at most $n \choose 2$ different regions, including
\begin{enumerate}
\item[\rm 1.] At most $i_+i_-$ regions of elliptic type;
\item[\rm 2.] At most ${{i_+}\choose 2} + {{i_-}\choose 2}$ regions of hyperbolic type.
\item[\rm 3.] More than ${{i_0}\choose 2} + i_0(i_+ + i_-)$ regions of parabolic type.
\end{enumerate}
\end{theorem}
{\bf Proof.} Consider the characteristic function of $\mathfrak D$, written in the form: $$f_{\mathfrak D}(x + iy) = {\mathbf L} + 2{\rm Sym}({\mathbf M})x + 2{\rm Skew}({\mathbf M})iy.$$
By Theorem \ref{SpecSym}, the symmetric matrix ${\rm Sym}({\mathbf M})$ can be reduced to its diagonal form by an orthogonal transformation ${\mathbf Q}$. Applying this transformation to $f_{\mathfrak D}$, we obtain:
$$f_{\mathfrak D}(x + iy) = {\mathbf Q}{\mathbf L}{\mathbf Q}^T + 2{\mathbf Q}{\rm Sym}({\mathbf M}){\mathbf Q}^Tx + 2{\mathbf Q}{\rm Skew}({\mathbf M}){\mathbf Q}^Tiy =$$
$$\widetilde{{\mathbf L}} + 2{\mathbf \Lambda}x + 2{\mathbf N}iy,$$
where ${\mathbf \Lambda}_{{\rm Sym}({\mathbf M})} = {\rm diag}\{\lambda_1, \ \ldots, \ \lambda_n\}$ is the canonical form of ${\rm Sym}({\mathbf M})$ and ${\mathbf N} = {\mathbf Q}{\rm Skew}({\mathbf M}){\mathbf Q}^T$ is a skew-symmetric matrix.

Applying Lemma \ref{App2}, we obtain the inclusion: $${\mathfrak D} \subseteq \bigcap_{(i,j)} P_{(i,j)},$$ where $1 \leq i < j \leq n$, $P_{(i,j)}$ is a region, bounded by a second-order curve:
\begin{equation} P_{(i,j)} = \{z = x+iy \in {\mathbb C}: a_{11}^{(i,j)}x^2 + a_{22}^{(i,j)}y^2 + 2a_{13}^{(i,j)}x + a_{33}^{(i,j)} > 0\},\end{equation}
where $a_{11}^{(i,j)} = {\mathbf \Lambda}_{{\rm Sym}({\mathbf M})}(i,j)$, $a_{22}^{(i,j)} = -{\mathbf N}(i,j)$.
Taking into account that ${\mathbf \Lambda}_{{\rm Sym}({\mathbf M})}(i,j) = \lambda_i \lambda_j$ and ${\mathbf N}(i,j) = \nu \geq 0$ (by skew-symmetry of $\mathbf N$), we calculate the second-order determinant $$\delta(i,j) = \begin{vmatrix}a_{11} & a_{12} \\ a_{21} & a_{22} \\ \end{vmatrix} =
\begin{vmatrix}\lambda_i\lambda_j & 0 \\ 0 & -\nu \\ \end{vmatrix} = -\nu\lambda_i\lambda_j. $$
Now let us consider the following three cases.
\begin{enumerate}
\item[\rm 1.] $\delta(i,j) > 0$. By \cite{BOC}, this case corresponds to $P_{(i,j)}$, bounded by a curve of ellyptic type. This happens when $\lambda_i\lambda_j < 0$, i.e. there is at least one pair of eigenvalues of different signs.
\item[\rm 2.] $\delta(i,j) < 0$. This case corresponds to $P_{(i,j)}$, bounded by a curve of hyperbolic type. This happens when $\lambda_i\lambda_j > 0$. Obviously, if ${\mathbf M}$ is semidefinite with ${\rm rank}({\rm Sym}({\mathbf M})) \geq 3$, at least one such pair do exists.
\item[\rm 3.] $\delta(i,j) = 0$. This case corresponds to $P_{(i,j)}$, bounded by a curve of parabolic type. This happens when there is at least one $\lambda_i =  0$, i.e. ${\rm Sym}({\mathbf M})$ is singular.
\end{enumerate}
 $\square$
 \begin{corollary} A necessary condition that the intersection $\bigcap_{(i,j)}P_{(i,j)}$ contains:
  \begin{enumerate}
  \item[\rm 1.] At least one region of parabolic type is that ${\rm Sym}({\mathbf M})$ is singular;
  \item[\rm 2.] At least one region of elliptic type is that ${\rm Sym}({\mathbf M})$ is indefinite;
  \item[\rm 3.] At least one region of hyperbolic type is that ${\rm Sym}({\mathbf M})$ is definite or ${\rm Rank}({\rm Sym}({\mathbf M})) \geq 3$.
  \end{enumerate}
 \end{corollary}

 Let us consider the following special case of Theorem \label{arbit}.

\begin{theorem} Let a nonempty LMI region $\mathfrak D$ be defined by its characteristic function $f_{\mathfrak D} = {\mathbf L} + {\mathbf M}z + {\mathbf M}^T\overline{z}$ of order $n\geq 2$ with ${\mathbf M}$ being normal.
Then $${\mathfrak D} \subseteq \bigcap_{(i,j)}P_{(i,j)},$$
where $1 \leq i < j \leq n$, each $P_{(i,j)}$ is a polynomial region of hyperbolic or parabolic type defined by \eqref{poly}.
\end{theorem}
{\bf Proof.} Since $\mathbf M$ is normal, by Lemma \ref{norm}, ${\rm Sym}({\mathbf M})$ and ${\rm Skew({\mathbf M})}$ commute. Thus, by Theorem \ref{commute}, they are simultaneously reduced to a diagonal and quasi-diagonal form, respectively, by an orthogonal transformation $\mathbf Q$. Applying this transformation to $f_{\mathfrak D}$, we obtain
$$f_{\mathfrak D}(x + iy) = \widetilde{{\mathbf L}} + 2 {\mathbf \Lambda}_{{\rm Sym}({\mathbf M})}x + 2 {\mathbf \Lambda}_{{\rm Skew}({\mathbf M})}iy,$$
where ${\mathbf \Lambda}_{{\rm Sym}({\mathbf M})}$ and ${\mathbf \Lambda}_{{\rm Skew}({\mathbf M})}$ are the canonical forms \eqref{canonsym} and \eqref{canonskew}, respectively. Consider the determinant $\delta(i,j) = -\nu\lambda_i\lambda_j$. By Corollary \ref{commute4}, $\nu > 0$ implies $\lambda_i\lambda_j$ and then $\delta(i,j) = -\nu\lambda_i\lambda_j < 0$ whenever $\lambda_i \neq 0$. Then the case of $\lambda_i \neq 0$ corresponds to a hyperbolic region and the case of $\lambda_i = 0$ --- to a parabolic region.

\section{Characteristics and maps of LMI regions}
\subsection{Embedding relations between LMI regions}

Let us start with the following theorem which shows the embedding relations between LMI regions.
\begin{theorem}\label{incl1} Given two LMI regions ${\mathfrak D}_1$ and ${\mathfrak D}_2$, defined by their characteristic functions $f_{{\mathfrak D}_1} = {\mathbf L}_1 + 2{\rm Sym}({\mathbf M}_1)x + 2{\rm Skew}({\mathbf M}_1)iy$ and $f_{{\mathfrak D}_2} = {\mathbf L}_2 + 2{\rm Sym}({\mathbf M}_2)x + 2{\rm Skew}({\mathbf M}_2)iy$. Let ${\mathbf L}_1 \preceq {\mathbf L}_2$ and one of the following cases holds:
\begin{enumerate}
\item[\rm 1.] ${\mathfrak D}_2 \in {\mathbb C}^+$, ${\rm Sym}({\mathbf M}_1) \preceq {\rm Sym}({\mathbf M}_2)$, ${\rm Skew}({\mathbf M}_1) = {\rm Skew}({\mathbf M}_2)$;
\item[\rm 2.] ${\mathfrak D}_2 \in {\mathbb C}^-$ ${\rm Sym}({\mathbf M}_2) \preceq {\rm Sym}({\mathbf M}_1)$, ${\rm Skew}({\mathbf M}_1) = {\rm Skew}({\mathbf M}_2)$;
\item[\rm 3.] no conditions on ${\mathfrak D}_2$, ${\mathbf M}_1 = {\mathbf M}_2$.
\end{enumerate}
 Then ${\mathfrak D}_2 \subseteq {\mathfrak D}_1$.
\end{theorem}
{\bf Proof.} Let Case 1 holds. Take an arbitrary $z \in {\mathfrak D}_2$. Then $z = x+iy \in {\mathfrak D}_2$ implies $x > 0$. Consider ${\mathbf W}_2(z): = f_{{\mathfrak D}_2}(z) = {\mathbf L}_2 + 2 {\rm Sym}({\mathbf M}_2)x + 2{\rm Skew}({\mathbf M}_2)iy \prec 0$. Then for ${\mathbf W}_1(z): = f_{{\mathfrak D}_1}(z)$ we obtain $${\mathbf W}_1(z)= {\mathbf L}_1 + 2 {\rm Sym}({\mathbf M}_1)x + 2{\rm Skew}({\mathbf M}_1)iy = $$ $${\mathbf L}_2 + ({\mathbf L}_1 - {\mathbf L}_2) + 2 {\rm Sym}({\mathbf M}_2)x + 2({\rm Sym}({\mathbf M}_1) - {\rm Sym}({\mathbf M}_2))x +2{\rm Skew}({\mathbf M}_2)iy = $$
$$({\mathbf L}_1 - {\mathbf L}_2 ) +2({\rm Sym}({\mathbf M}_1) - {\rm Sym}({\mathbf M}_2))x + {\mathbf W}_2(z).$$
Since ${\mathbf L}_1 - {\mathbf L}_2 \preceq 0$ and ${\rm Sym}({\mathbf M}_1) - {\rm Sym}({\mathbf M}_2) \preceq 0$, we have $\widetilde{{\mathbf W}}(z) := ({\mathbf L}_1 - {\mathbf L}_2) + 2 ({\rm Sym}({\mathbf M}_2) - {\rm Sym}({\mathbf M}_1))x \preceq 0$ whenever $x > 0$. Thus ${\mathbf W}_1(z) = \widetilde{{\mathbf W}}(z) + {\mathbf W}_2(z) \prec 0$ as the sum of $ \widetilde{{\mathbf W}}(z) \preceq 0$ and ${\mathbf W}_2(z) \prec 0$, and we conclude that $z \in {\mathfrak D}_1$.

The Cases 2 and 3 are considered analogically. $\square$

Now let us prove the following important lemma, which shows the role of the generating matrix $\mathbf L$ in the location of an LMI region on $\mathbb C$.

\begin{lemma}\label{zero} An LMI region $\mathfrak D$ contains $0$ if and only if $\mathbf L$ is negative definite. Its closure $\overline{{\mathfrak D}}$ contains $0$ if and only if $\mathbf L$ is negative semidefinite.
\end{lemma}
The proof follows immediately from the substitution $z = 0$.
\begin{corollary} Given an LMI region $\mathfrak D$, let $z, - z \in {\mathfrak D}$ for some $z \in {\mathbb C}$. Then $\mathbf L$ is negative definite.
\end{corollary}

Now we can answer the question, when an LMI region $\mathfrak D$ or its closure $\overline{{\mathfrak D}}$ contains the recession cone ${\mathfrak D}_{rc}$. For this, we first prove the following lemma.

\begin{theorem}\label{incl} Let $\mathfrak D$ be an LMI region, defined by its characteristic function $ f_{{\mathfrak D}} = {\mathbf L} + {\mathbf M}z+{\mathbf M}^T\overline{z}$. Then the following statements hold.
\begin{enumerate}
\item[\rm 1.]  $\mathbf L$ is negative definite if and only if ${\mathfrak D}_{rc} \subset {\mathfrak D}$.
\item[\rm 2.]  $\mathbf L$ is negative semidefinite if and only if ${\mathfrak D}_{rc} \subset \overline{{\mathfrak D}}$.
\end{enumerate}
\end{theorem}
{\bf Proof.}
{\bf Case 1.} $\Rightarrow$ By Lemma \ref{zero}, $\mathbf L$ is negative definite implies $0 \in {\mathfrak D}$. Thus by definition of the recession cone, $0 + {\mathfrak D}_{rc} = {\mathfrak D}_{rc} \subseteq {\mathfrak D}$.

$\Leftarrow$ Let ${\mathfrak D}_{rc} \subseteq {\mathfrak D}$. By definition, $0 \in {\mathfrak D}_{rc}$. Thus $0 \in {\mathfrak D}$ and by Lemma \ref{zero}, $\mathbf L$ is negative definite.

{\bf Case 2.} The case of semidefinite $\mathbf L$ is considered analogically.
 $\square$

\subsection{Maps of LMI regions}
Here, let us ask a general question: which classes of maps $\varphi: {\mathbb C} \rightarrow {\mathbb C}$ possess the following property: $\varphi(\mathfrak D)$ is an LMI region whenever $\mathfrak D$ is an LMI region? We consider the two simplest classes of maps with the above property.

{\bf Class 1.} $\varphi: z \rightarrow \alpha z$. This class includes stretching ($\alpha > 1$), contraction ($\alpha < 1$) and reflection with respect to the origin ($\alpha = -1$). For this class, the following statement holds.

\begin{theorem} Let $\mathfrak D$ be an LMI region, defined by its characteristic function $f_{\mathfrak D} = {\mathbf L} + {\mathbf M}z + {\mathbf M}^T\overline{z}$. Then, for any $\alpha \in {\mathbb R}, \ \alpha \neq 0$, $\alpha{\mathfrak D}$ is also an LMI region, which is nonempty if and only if $\mathfrak D$ is nonempty, with the characteristic function $f_{\alpha\mathfrak D} = |\alpha| {\mathbf L} + {\rm Sign} (\alpha){\mathbf M}z + {\rm Sign} (\alpha){\mathbf M}^T\overline{z}$. Moreover, if one of the following cases hold:
\begin{enumerate}
\item[\rm 1.] ${\mathbf L} \succeq 0$, $\alpha \in [-1,0)\cup [1, + \infty)$;
\item[\rm 2.] ${\mathbf L} \preceq 0$, $\alpha \in (-\infty,-1]\cup (0, 1]$,
\end{enumerate}
then $\alpha{\mathfrak D} \subseteq {\rm Sign} (\alpha){\mathfrak D}$.
\end{theorem}
{\bf Proof.} Since $z \in \alpha{\mathfrak D}$ if and only if $\frac{1}{\alpha}z \in {\mathfrak D}$, we obtain $\alpha{\mathfrak D} \neq \emptyset$ if and only if $\mathfrak D \neq \emptyset$. Then, for $z \in \alpha{\mathfrak D}$, we have ${\mathbf L} + {\mathbf M}\frac{1}{\alpha}z + {\mathbf M}^T\frac{1}{\alpha}\overline{z} \prec 0$, which is equivalent to $\alpha{\mathbf L} + {\mathbf M}z + {\mathbf M}^T\overline{z} \prec 0$ if $\alpha > 0$ and $-\alpha{\mathbf L} - {\mathbf M}z - {\mathbf M}^T\overline{z} \prec 0$ if $\alpha < 0$. Next,
${\mathbf L} \succeq 0$, $\alpha \geq 1$ implies $(1 - \alpha){\mathbf L} \preceq 0$ and ${\mathbf L} \preceq \alpha{\mathbf L}$. In its turn, ${\mathbf L} \preceq 0$, $0 < \alpha \leq 1$ also implies $(1 - \alpha){\mathbf L} \preceq 0$ and ${\mathbf L} \preceq \alpha{\mathbf L}$. Applying Theorem \ref{incl1} to both of the cases, we obtain the inclusion $\alpha{\mathfrak D} \subseteq {\mathfrak D}$.

The case 2 is considered analogically. $\square$

 {\bf Class 2.} $\varphi: z \rightarrow z+ \alpha$, where $\alpha \in {\mathbb R}$. This class consists of shifts along the real axis.
 \begin{theorem}\label{Shift} Let $\mathfrak D$ be an LMI region, defined by its characteristic function $f_{\mathfrak D} = {\mathbf L} + {\mathbf M}z + {\mathbf M}^T\overline{z}$. Then $\alpha + {\mathfrak D}$ is an LMI region with the characteristic function \begin{equation}\label{shiftfun}f_{\alpha + \mathfrak D} = \widetilde{{\mathbf L}} + {\mathbf M}z + {\mathbf M}^T\overline{z}, \end{equation} where $\widetilde{{\mathbf L}} = {\mathbf L} - 2\alpha {\rm Sym}{\mathbf M}$, for any $\alpha \in {\mathbb R}$. Moreover, if one of the following cases hold:
\begin{enumerate}
\item[\rm 1.] ${\mathbf M} \succeq 0$, $\alpha \leq 0$;
\item[\rm 2.] ${\mathbf M} \preceq 0$, $\alpha \geq 0$,
\end{enumerate}
then the following inclusion holds: $\alpha + {\mathfrak D} \subseteq {\mathfrak D}$.
\end{theorem}
 {\bf Proof.} Since $z \in \alpha + {\mathfrak D}$ if and only if $z - \alpha \in {\mathfrak D}$, we have ${\mathbf L} + {\mathbf M}(z -\alpha) + {\mathbf M}^T(\overline{z} - \alpha) \prec 0$, which is equivalent to ${\mathbf L} - \alpha({\mathbf M}+{\mathbf M}^T) + {\mathbf M}z + {\mathbf M}^T\overline{z} \prec 0$. Then, both of the cases ${\mathbf M} \succeq 0$, $\alpha \leq 0$ and ${\mathbf M} \preceq 0$, $\alpha \geq 0$ imply $\alpha({\mathbf M}+{\mathbf M}^T)\preceq 0$ and further, ${\mathbf L} \preceq {\mathbf L} - \alpha({\mathbf M}+{\mathbf M}^T)$. Applying Lemma \ref{incl1} to both of the cases, we complete the proof. $\square$

 \begin{corollary}\label{shiftdef} Let $f_{\mathfrak D} = {\mathbf L} + 2x{\rm Sym}{\mathbf M} + 2iy{\rm Skew}{\mathbf M}$ defines a nonempty LMI region $\mathfrak D$. Then there is $\alpha \in {\mathbb R}$ such that ${\mathbf L} + \alpha{\rm Sym}{\mathbf M} \prec 0$.
 \end{corollary}
 {\bf Proof}. Let ${\mathfrak D} \neq \emptyset$. Then there is $z_0 = x_0+iy_0 \in {\mathfrak D}$. By symmetry (Property 1 of LMI regions), $\overline{z}_0 \in {\mathfrak D}$ and by convexity (Property 2 of LMI regions), $x_0 = \frac{z_0 + \overline{z}_0}{2} \in {\mathfrak D}$. Consider the shift $\varphi: z \rightarrow z - x_0$, which maps $x_0$ to $0$. By Theorem \ref{Shift}, ${\mathfrak D} - x_0$ is an LMI region, defined by the characteristic function \ref{shiftfun}. Since  ${\mathfrak D} - x_0$ contains $0$, by Lemma \ref{zero}, $\widetilde{{\mathbf L}} = {\mathbf L} - 2x_0 {\rm Sym}{\mathbf M} \prec 0$. Putting $\alpha = x_0$, we complete the proof. $\square$

\subsection{When an LMI region is empty?}

Given an arbitrary matrix-valued function of the form $f(z) = {\mathbf L} + {\mathbf M}z + {\mathbf M}^T\overline{z}$, here we clarify the question, when the inequality $f(z) \prec 0$ defines an empty set. For this, we use the following construction. By substitution $z = x \in {\mathbb R}$, we obtain the function $f(x) = {\mathbf L} + 2x{\rm Sym}({\mathbf M})$. The inequality $f(x) \prec 0$ defines an (empty or non-empty) interval ${\mathfrak D}_{\mathbb R}$ of the real line as follows:
\begin{equation}\label{dr}{\mathfrak D}_{\mathbb R} = \{x \in {\mathbb R}: {\mathbf L} + 2x{\rm Sym}({\mathbf M})\prec0\}. \end{equation}
In general, ${\mathfrak D}_{\mathbb R}$ is not an LMI region.

An LMI region $\mathfrak D$ and the interval ${\mathfrak D}_{\mathbb R}$ of the real line are connected by the following lemma.

 \begin{lemma}\label{critempty}
Given an LMI region ${\mathfrak D} \subseteq {\mathbb C}$, defined by its characteristic function $f_{\mathfrak D} = {\mathbf L} + 2{\rm Sym}({\mathbf M})x + 2{\rm Skew}({\mathbf M})iy$. Then $\mathfrak D$ is non-empty if and only if the interval ${\mathfrak D}_{\mathbb R} \subseteq {\mathbb R}$, defined by \eqref{dr}, is nonempty, and ${\mathfrak D}_{\mathbb R} = {\mathfrak D} \cap {\mathbb R}$.
\end{lemma}
{\bf Proof.} $\Rightarrow$ Let an LMI region $\mathfrak D$ be non-empty. Then there is $z = x+iy \in {\mathfrak D}$. By Property 1 (Symmetry) of LMI regions, $z \in {\mathfrak D}$ implies $\overline{z} \in {\mathfrak D}$ and by Property 2 (Convexity), $z + \overline{z} = 2x \in {\mathfrak D}$. By substitution $z = x$ into the inequality $f_{\mathfrak D} \prec 0$, we get ${\mathbf L} + 2x{\rm Sym}({\mathbf M})\prec 0$ and conclude that ${\mathfrak D}_{\mathbb R}$ is non-empty and ${\mathfrak D}\cap {\mathbb R} \subseteq {\mathfrak D}_{\mathbb R}$.

$\Leftarrow$ Let ${\mathfrak D}_{\mathbb R} \neq \emptyset$. Then, there is $x_0 \in {\mathbb R}$ such that ${\mathbf W}(x_0) = {\mathbf L} + 2x_0{\rm Sym}({\mathbf M})\prec0$. Consider $z = x_0 + iy$, where $y$ will be chosen later to be sufficiently small. Then
$${\mathbf W}(z) := {\mathbf L}= + 2{\rm Sym}({\mathbf M})x_0 + 2{\rm Skew}({\mathbf M})iy =$$ $${\mathbf W}_0 + 2{\rm Skew}({\mathbf M})iy .$$
Since ${\mathbf W}(x_0)$ is negative definite and the set of Hermitian negative definite matrices is open in ${\mathcal M}^{n \times n}({\mathbb C})$, we can choose $\delta > 0$ such that for all skew-symmetric ${\mathbf \Delta} \in {\mathcal M}^{n \times n}{\mathbb C}$ satisfying $\|\Delta\| < \delta$, the matrix ${\mathbf W}(x_0) + i{\mathbf \Delta}$ is Hermitian negative definite. Then, taking $y<\frac{\delta}{2\|{\rm Skew}({\mathbf M})\|}$, we obtain $${\mathbf W}(z) = {\mathbf W}(x_0)+ \frac{2\delta{\rm Skew}({\mathbf M})}{2\|{\rm Skew}({\mathbf M})\|}i,$$
and $\|{\mathbf W}(z) - {\mathbf W}(x_0)\| < \delta$. Thus ${\mathbf W}(z)$ is also negative definite.

Hence ${\mathfrak D}_{\mathbb R} \subseteq {\mathfrak D}\cap {\mathbb R}$.
 $\square$

Lemma \ref{critempty} shows that checking if an LMI region $\mathfrak D$ is nonempty, is equivalent to checking the feasibility of the following LMI:
\begin{equation}\label{fesp}F(x):={\mathbf L} + 2x{\rm Sym}({\mathbf M}) \prec 0, \end{equation}
where both the matrices $\mathbf L$ and ${\rm Sym}({\mathbf M})$ are symmetric. Such LMI are widely studied (see, for example, \cite{BOYE}, \cite{VB}), and LMI feasibility problem (i.e. finding $x \in {\mathbf R}$ such that $F(x) \prec 0$) is considered by numerical methods. Thus checking feasibility of \label{fesp} is now practically possible for any given symmetric matrix ${\mathbf L}$ and arbitrary matrix $\mathbf M$. The ellipsoid algorithm (see, for example, \cite{BOYE}, \cite{IB}), guarantees to solve this problem. The LMI solver in MathLab is based on the interior-point methods (see \cite{NN}, \cite{NG}, \cite{GNLC}) which allows us to check feasibility in polynomial time.
Thus we can find a point $\hat{x} \in {\mathfrak D}\cap {\mathbb R}$ with the help of feasibility solver in MathLab (feasp).

In this paper, we suggest another method of checking if an LMI region is nonempty and finding a point $\hat{x} \in {\mathfrak D}\cap {\mathbb R}$ based on the following necessary condition.

\begin{lemma}\label{condempt}  Given an LMI region ${\mathfrak D}$, defined by its characteristic function $ f_{\mathfrak D} = {\mathbf L}+{\mathbf M}z+{\mathbf M}^T\overline{z}$ with ${\rm Sym}({\mathbf M})$ being nonsingular. Then it is necessary for ${\mathfrak D} \neq \emptyset$ that ${\mathbf C}:=({\rm Sym}({\mathbf M}))^{-1}{\mathbf L}$ has real eigenvalues and is diagonalizable (i.e. there is a nonsingular ${\mathbf R} \in {\mathcal M}^{n \times n}$ such that ${\mathbf R}^{-1}{\mathbf C}{\mathbf R}$ is a real diagonal matrix).
\end{lemma}
{\bf Proof.} Let ${\mathfrak D}$ be non-empty. Then, by Corollary \ref{shiftdef}, we obtain that there is $\alpha \in {\mathbb R}$ such that ${\mathbf L} + \alpha{\rm Sym}{\mathbf M} \prec 0$. Since ${\mathbf L} + \alpha{\rm Sym}{\mathbf M}$ is negative definite, by Lemma \ref{redsym}, $\widetilde{{\mathbf L}}:= {\mathbf L} + \alpha{\rm Sym}{\mathbf M}$ and ${\rm Sym}{\mathbf M}$ are simultaneously diagonalizable by congruence. By Theorem \ref{ns}, this implies that the matrix $\widetilde{{\mathbf C}} := ({\rm Sym}({\mathbf M}))^{-1}\widetilde{{\mathbf L}} = ({\rm Sym}({\mathbf M}))^{-1}{\mathbf L} + \alpha{\mathbf I}$ has real eigenvalues and is diagonalizable. It means, that for some nonsingular ${\mathbf R} \in {\mathcal M}^{n \times n}$, we have $${\mathbf R}^{-1}\widetilde{{\mathbf C}}{\mathbf R} = {\mathbf \Lambda},$$
where ${\mathbf \Lambda}$ is a real diagonal matrix. However,
$${\mathbf R}^{-1}{\mathbf C}{\mathbf R} = {\mathbf R}^{-1}(\widetilde{{\mathbf C}}-\alpha{\mathbf I}){\mathbf R} = {\mathbf \Lambda}-\alpha{\mathbf I},$$
thus ${\mathbf C}$ is diagonalizable if and only if $\widetilde{{\mathbf C}}$ is diagonalizable.
$\square$

Simple examples with diagonal matrices show that this condition is not sufficient.

\subsection{Characteristics of inner angles of LMI regions} Consider a uniform LMI region $\mathfrak D$. By Theorem \ref{Cone} it coincides with a cone in $\mathbb C$. Now we can find its inner angle using the polar decomposition of the generating matrix $\mathbf M$.
\begin{theorem}\label{angle} Given a non-empty uniform LMI region ${\mathfrak D}$, defined by its characteristic function $ f_{\mathfrak D} = {\mathbf M}z+{\mathbf M}^T\overline{z}$. Then $\mathbf M$ is positive (or negative) definite, $\mathfrak D$ is a cone around the negative (respectively, positive) direction of the real axis with the inner angle
 $$\theta = \min_j|\frac{\pi}{2} - \arg (\lambda_j({\mathbf U}))|,$$
where ${\mathbf U} = {\mathbf P}^{-1}{\mathbf M}$, ${\mathbf P} = ({\mathbf M}{\mathbf M}^T)^{\frac{1}{2}}$, and
 $\{\lambda_j({\mathbf U})\}_{j = 1}^n =\sigma({\mathbf U})$.
\end{theorem}
{\bf Proof.} Let $\mathbf D$ be a non-empty uniform LMI region. Then, by Corollary \ref{copos}, $\mathbf M$ is (positive or negative) definite. Consider the case, when $\mathbf M$ is positive definite (the case of negative definite $\mathbf M$ can be easily considered by analogy). By Theorem \ref{Cone}, $\mathfrak D$ coincides with an open cone around the negative direction of the real axis. Let $z \in {\mathfrak D}$. Without losing the generality of the reasoning, we may assume $|z| = 1$ and $z = e^{i\varphi}$ for some $\varphi \in [0, \pi]$. Thus
\begin{equation}\label{con} {\mathbf M}e^{i\varphi} + {\mathbf M}^Te^{-i\varphi} \prec 0.\end{equation}
Consider the polar decomposition of $\mathbf M$ (see Theorem \ref{polar}):
${\mathbf M} = {\mathbf P}{\mathbf U}$, where ${\mathbf P} = ({\mathbf M}{\mathbf M}^T)^{\frac{1}{2}}$ is symmetric positive definite, ${\mathbf U}$ is a unitary matrix. Note, that the Lyapunov theorem (Theorem \ref{lyap}) implies that ${\mathbf U}$ is positive stable, i.e. $-\frac{\pi}{2} < {\rm arg}(\lambda_j({\mathbf U})) < \frac{\pi}{2}$, $j = 1, \ \ldots, \ n$.

Re-writing \eqref{con}, we obtain the equation
$${\mathbf P}(e^{i\varphi}{\mathbf U}) + (e^{i\varphi}{\mathbf U})^*{\mathbf P} = {\mathbf W}(z) \prec 0. $$
Then, applying the Lyapunov theorem (Theorem \eqref{lyap}), we obtain that
$e^{i\varphi}{\mathbf U}$ is stable if and only if $z = e^{i\varphi} \in {\mathfrak D}$. This gives the following condition:
$$ -\frac{\pi}{2} < {\rm arg}(\lambda_j({\mathbf U})) + \varphi < \frac{\pi}{2}, \qquad j = 1, \ \ldots, \ n.$$
Thus
$$ -\frac{\pi}{2}- {\rm arg}(\lambda_j({\mathbf U})) < \varphi < \frac{\pi}{2}-{\rm arg}(\lambda_j({\mathbf U})), \qquad j = 1, \ \ldots, \ n.$$
Taking into account that the complex eigenvalues of a real matrix appears in conjugate pairs, we obtain
$$ \pi \leq \varphi < \frac{\pi}{2}+\max_j{\rm arg}(\lambda_j({\mathbf U})). $$
From the equality $\theta = \pi - \inf\{\varphi: \pi \leq \varphi < {\pi}{2}; e^{i\varphi} \in {\mathfrak D}\}$, we obtain $\theta = \frac{\pi}{2} - \max_j({\rm arg}(\lambda_j({\mathbf U}))) = \min_j|\frac{\pi}{2} - \arg (\lambda_j({\mathbf U}))|$.
 $\square$
\begin{corollary}\label{normcone} Given a non-empty uniform LMI region ${\mathfrak D}$, defined by its characteristic function $ f_{\mathfrak D} = {\mathbf M}z+{\mathbf M}^T\overline{z}$, with $\mathbf M$ being normal. Then $\mathbf M$ is positive (or negative) definite, $\mathfrak D$ is a cone around the negative (respectively, positive) direction of the real axis with the inner angle
 $$\theta = \min_j|\frac{\pi}{2} - \arg (\lambda_j({\mathbf M}))|,$$
where $\{\lambda_j({\mathbf M})\}_{j = 1}^n =\sigma({\mathbf M})$.
\end{corollary}
{\bf Proof.} If $\mathbf M$ is normal, the eigenvalues of $\mathbf U$ and $\mathbf M$ are connected by Lemma \ref{un}. Thus we may replace the eigenvalues of $\mathbf U$ by eigenvalues of $\mathbf M$. $\square$

Corollary \ref{normcone} also can be easily deduced from Theorem \ref{comp1}.

Note, that the condition of definiteness of ${\mathbf M}$ implies, that $\mathbf M$ and ${\mathbf M}^T$ can be simultaneously reduced to their canonical forms. Indeed, applying Lemma \ref{redskew} to ${\rm Sym}({\mathbf M})$, which is definite and ${\rm Skew}({\mathbf M})$, which is skew-symmetric, and taking into account that ${\mathbf M} = \frac{1}{2}({\rm Sym}({\mathbf M}) + {\rm Skew}({\mathbf M}))$ and ${\mathbf M}^T = \frac{1}{2}({\rm Sym}({\mathbf M}) - {\rm Skew}({\mathbf M}))$, we get that $\mathbf M$ and ${\mathbf M}^T$ are also simultaneously reducible. Thus definiteness of ${\mathbf M}$ is closed to normality of ${\mathbf M}$, however, is not equivalent. Indeed, consider a shift of the form ${\mathbf M} + \alpha{\mathbf I}$, $\alpha \in {\mathbb R}$. For sufficiently big $|\alpha|$, we obtain ${\rm Sym}({\mathbf M} + \alpha{\mathbf I}) = {\rm Sym}({\mathbf M}) + 2\alpha{\mathbf I}$ is definite, but by Lemma \ref{norm}, is normal if and only if the initial matrix ${\mathbf M}$ is normal.

Also note that we can provide another way of calculating the inner angle of a non-empty uniform region $\mathfrak D$. For this, we write the characteristic function $f_{\mathfrak D}$ in the form: $$ f_{\mathfrak D} = 2{\rm Sym}({\mathbf M})x + 2{\rm Skew}({\mathbf M})y.$$
By Corollary \ref{copos}, $\mathbf M$ is (positive or negative) definite. Consider the case, when $\mathbf M$ is positive definite (the case of negative definite $\mathbf M$ can be easily considered by analogy). By Lemma \ref{redskew}, we simultaneously reduce by congruence ${\rm Sym}({\mathbf M})$ to ${\mathbf I}$ and ${\rm Skew}({\mathbf M})$ to ${\rm Skew}({\mathbf M})_{{\rm Sym}({\mathbf M})}$, which is a quasi-diagonal matrix of the form \eqref{qdiagskew} (note, that its eigenvalues are different from those of ${\rm Skew}({\mathbf M})$). Then by Lemma \ref{App3}, ${\mathfrak D} = \bigcap{\mathfrak D}_i$, where each $\bigcap{\mathfrak D}_i$ is either the left hand-side halfplane ${\mathbb C}^{-}$ or a cone around the negative direction of the real axis, defined by the following LMI:
$$2\begin{pmatrix} 1 & 0 \\ 0 & 1 \\ \end{pmatrix}x + 2\begin{pmatrix}0 & \nu_i \\ -\nu_i & 0 \\ \end{pmatrix}iy \prec 0,$$
where $\pm \nu_i$ are the nonzero pure imaginary eigenvalues of ${\rm Skew}({\mathbf M})_{{\rm Sym}({\mathbf M})}$. Thus we get
$\theta = \min_i\theta_i$, where $\tan(\theta_i) = \frac{1}{\nu_i}$.

The above reasoning also provides the link between the eigenvalues of the unitary matrix $\mathbf U$ in the polar decomposition of ${\mathbf M}$, and the eigenvalues of the symmetric and skew-symmetric parts ${\rm Sym}({\mathbf M})$ and ${\rm Skew}({\mathbf M})$, respectively.

Using the above results, we can calculate the angle $\theta$ of the recession cone ${\mathfrak D}_{rc}$ of an arbitrary unbounded LMI region $\mathfrak D$.

\begin{theorem}Given a non-empty unbounded LMI region ${\mathfrak D}$, defined by its characteristic function $ f_{\mathfrak D} = {\mathbf L}+ {\mathbf M}z+{\mathbf M}^T\overline{z}$. Then one of the following two cases holds:
\begin{enumerate}
 \item[\rm 1.] ${\rm Sym}({\mathbf M})$ is singular and the angle $\theta$ of the cone ${\mathfrak D}_{rc}$ is $0$.
  \item[\rm 2.] ${\rm Sym}({\mathbf M})$ is nonsingular and the angle $\theta$ of the cone ${\mathfrak D}_{rc}$ satisfies $$\theta = \min_j|\frac{\pi}{2} - \arg (\lambda_j({\mathbf U}))|,$$
where ${\mathbf U} = {\mathbf P}^{-1}{\mathbf M}$, ${\mathbf P} = ({\mathbf M}{\mathbf M}^T)^{\frac{1}{2}}$, and
 $\{\lambda_j({\mathbf U})\}_{j = 1}^n =\sigma({\mathbf U})$.
\end{enumerate}
In case ${\mathbf M}$ being normal, eigenvalues of ${\mathbf U}$ can be replaced with the eigenvalues of ${\mathbf M}$.
\end{theorem}
{\bf Proof.} By Theorem \ref{resco}, the recession cone ${\mathfrak D}_{rc}$ is defined by LMI \eqref{DU}, In its turn, by Lemma \ref{Clos},
$${\rm int}({\mathfrak D}_{rc}) = \{z \in {\mathbb C}: {\mathbf M}z+{\mathbf M}^T\overline{z} \prec 0\}.$$
Applying to ${\rm int}({\mathfrak D}_{rc})$ Theorem \ref{angle} and Corollary \ref{normcone}, we complete the proof. $\square$

\subsection{Characteristics of inner radii of LMI regions} Let us consider the problem of disk placement, which often arises in robust control: given an LMI region $\mathfrak D$ defined by its characteristic function $ f_{\mathfrak D} = {\mathbf L} + {\mathbf M}z+{\mathbf M}^T\overline{z}$, find a disk $D(x,r)$ of radius $r>0$, centered at the point $x \in {\mathbb R}$, such that $D(x,r) \subseteq {\mathfrak D}$? Here, we consider the solution of the problem in terms of certain characteristics of matrices ${\mathbf M}$ and ${\mathbf L}$, using the following geometric constructions.

\begin{enumerate}
\item[{\rm Step 1.}] We check the feasibility of $f_{\mathfrak D}$ and find $\hat{x} \in {\mathfrak D}\cap{\mathbb R}$. We also calculate the boundaries of ${\mathfrak D}_{\mathbb R}$, i.e. we find a (finite or infinite) interval $(x_{min}, x_{max}) \subseteq {\mathbb R}$ such that $(x_{min}, x_{max}) = {\mathfrak D}_{\mathbb R}$.
\item[{\rm Step 2.}] We check, if $0 \in {\mathfrak D}$ by checking the negative definiteness of ${\mathbf L}$. If $0$ does not belong to $\mathfrak D$, then, fixing $\hat{x} \in {\mathfrak D}_{\mathbb R}$, we apply the shift along the real axes $x':=x - \hat{x}$, moving $\hat{x}$ to $0$ and denote ${\mathfrak D}':={\mathfrak D} - \hat{x}$. We calculate the new boundaries $x'_{min}, \ x'_{max}$ for
 ${\mathfrak D}'_{\mathbb R} = {\mathfrak D}'\cap{\mathbb R}$.
\item[{\rm Step 3.}] We fix $x = x_0 \in (x'_{min}, x'_{max})$, to be the center of the inscribed disk. Note that if the point $x_0$ is given, we can check if $x_0 \in {\mathfrak D}$ by checking the negative definiteness of the matrix ${\mathbf L}(x_0) = {\mathbf L} + 2{\rm Sym}({\mathbf M})x_0$. For the case ${\rm Skew}({\mathbf M}) \neq 0$, we find the lower bound $y$ for the intersection of the vertical line $x = x_0$ with $\partial{\mathfrak D}'$.
\item[{\rm Step 4.}] For $x=x_0$, we find $\min((x - x'_{min}), (x'_{max} - x))$. Without loss the generality of the reasoning, we assume that it will be $x'_{max} - x$. From the right-angled triangle $x_{max}xy$ (see Figure 3), we find the altitude $r$ to the hypothenuse $yx_{max}$, using well-known formulae: \begin{equation}\label{radius}r({\mathfrak D}', x) = \frac{(x_{max} - x)y}{\sqrt{(x_{max} - x)^2 +y^2}}.\end{equation} Note, that the radius of an inscribed circle is invariant under the shifts, hence $r({\mathfrak D}',x) = r({\mathfrak D},x)$. Straightforward algorithms allows us to find the optimal placement of $x$ to maximize the radius of an inscribed circle, if such a problem will arise.
\end{enumerate}

\begin{figure}[h]
\center{\includegraphics[scale=0.5]{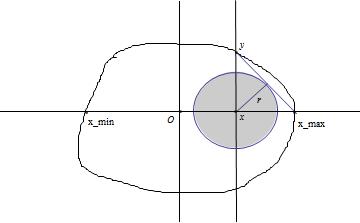}}
\caption{Circle, inscribed in an LMI region}
\end{figure}

Now consider each step in details.

{\bf Step 1}. Given an LMI region ${\mathfrak D}$, we consider its intersection with real line ${\mathfrak D}_{\mathbb R} := {\mathfrak D}\cap{\mathbb R}$.

The results of Section 4 allows us to find out, when ${\mathfrak D}_{\mathbb R} = {\mathbb R}$ and when it is an unbounded interval of the form $(-\infty, x_{max})$ or $(x_{min}, +\infty)$.

\begin{theorem} Let an LMI region $\mathfrak D \neq \emptyset$ be defined by its characteristic function $ f_{{\mathfrak D}} = {\mathbf L} + {\mathbf M}z+{\mathbf M}^T\overline{z}$ Then
\begin{enumerate}
\item[\rm 1.] ${\mathfrak D}_{\mathbb R} = {\mathbb R}$ if and only if ${\mathbf M}$ is skew-symmetric.
\item[\rm 2.] ${\mathfrak D}_{\mathbb R} = (-\infty, x_{max})$ for some value $ x_{max} \in {\mathbb  R}$ if and only if ${\mathbf M}$ is positive semidefinite and ${\rm Sym}(\mathbf M) \neq 0$.
\item[\rm 3.] ${\mathfrak D}_{\mathbb R} = (x_{min}, +\infty)$ for some value $ x_{min} \in {\mathbb  R}$ if and only if ${\mathbf M}$ is negative semidefinite and ${\rm Sym}(\mathbf M) \neq 0$.
\end{enumerate}
\end{theorem}
{\bf Proof.} Case {\bf 1} immediately follows from Theorem \ref{lin}, in this case ${\mathfrak D}_{\mathbb R} = L_{\mathfrak D}$.

Cases {\bf 2} and {\bf 3} immediately follows from Theorem \ref{semdef} and Theorem \ref{lin}.
$\square$

The following theorem based on Lemma \ref{App1}, provides the outer estimates for ${\mathfrak D}_{\mathbb R}$.

\begin{theorem} Let an LMI region ${\mathfrak D} \neq \emptyset$ be defined by its characteristic function $ f_{{\mathfrak D}} = {\mathbf L} + {\mathbf M}z+{\mathbf M}^T\overline{z}$ with ${\mathbf M} = \{m_{ij}\}_{i,j = 1}^n$ and ${\mathbf L} = \{l_{ij}\}_{i,j = 1}^n$. Define the subsets $I_1, I_2 \subseteq [n]$, where $i \in I_1$ if $m_{ii} > 0$ and $i \in I_2$ if $m_{ii} < 0$. Then
\begin{equation}\label{increal}{\mathfrak D}_{\mathbb R} \subseteq (x_{min}, x_{max}),\end{equation}
where
\begin{equation}\label{est1} x_{min} = \left\{\begin{array}{ccc} \max_{i \in I_2}\frac{-l_{ii}}{2m_{ii}} & \mbox{if} \ I_2 \neq \emptyset; \\
- \infty & \mbox{otherwise} \end{array}\right.\end{equation}
\begin{equation}\label{est2} x_{max} = \left\{\begin{array}{ccc} \min_{i \in I_1}\frac{-l_{ii}}{2m_{ii}} & \mbox{if} \ I_1 \neq \emptyset; \\
+ \infty & \mbox{otherwise} \end{array}\right.\end{equation}
If $\mathbf L$ and ${\mathbf M}$ are diagonal, then Inclusion \ref{increal} turns to the equality.
\end{theorem}
{\bf Proof.} By Lemma \ref{App1}, ${\mathfrak D} \subseteq {\mathfrak D}_1$, where ${\mathfrak D}_1$ is a nonempty LMI region, defined by $ f_{{\mathfrak D}_1} =  {\mathbf L}_1 + {\mathbf M}_1z+{\mathbf M}_1^T\overline{z}$, where ${\mathbf L}_1 = \{l_{11}, \ \ldots, \ l_{nn}\}$ and ${\mathbf M}_1 = \{m_{11}, \ \ldots, \ m_{nn}\}$. Thus ${\mathfrak D}_{\mathbb R} \subseteq {\mathfrak D}^1_{\mathbb R} = {\mathfrak D}_1\bigcap{\mathbb R}$. By Lemma \ref{Clos}, ${\mathfrak D}_1 = \bigcap_{i=1}^nP_i(z)$, where
$$P_i(z) = \{z = x+iy \in {\mathbb C}: l_{ii} + 2m_{ii}x <0 \}. $$
Thus ${\mathfrak D}^1_{\mathbb R} = \bigcap_{i=1}^n(P_i(z)\cap{\mathbb R})$.  From here we get
$x < -\frac{l_{ii}}{2m_{ii}}$ when $m_{ii} >0 $
that implies $x < \min_{i \in I_1}\frac{-l_{ii}}{m_{ii}},$
and $x > -\frac{l_{ii}}{2m_{ii}}$ when $m_{ii} < 0$
that implies $x > \max_{i \in I_2}\frac{-l_{ii}}{2m_{ii}}$.
 $\square$
\begin{corollary}\label{coest} Let an LMI region ${\mathfrak D} \neq \emptyset$ be defined by its characteristic function $ f_{{\mathfrak D}} = {\mathbf L} + {\mathbf M}z+{\mathbf M}^T\overline{z}$ and ${\mathbf L}$ commute with ${\rm Sym}({\mathbf M})$. Define the subsets $I_1, I_2 \subseteq [n]$, where $i \in I_1$ if $\lambda_i({\rm Sym}({\mathbf M})) > 0$ and $i \in I_2$ if $\lambda_i({\rm Sym}({\mathbf M})) < 0$. Then \begin{equation}{\mathfrak D}_{\mathbb R} = (x_{min}, x_{max}),\end{equation}
where
\begin{equation} x_{min} = \left\{\begin{array}{ccc} \max_{i \in I_2}\frac{-\lambda_i({\mathbf L})}{2\lambda_i({\rm Sym}({\mathbf M}))} & \mbox{if} \ I_2 \neq \emptyset; \\
- \infty & \mbox{otherwise} \end{array}\right.\end{equation}
\begin{equation} x_{max} = \left\{\begin{array}{ccc} \min_{i \in I_1}\frac{-\lambda_i({\mathbf L})}{2\lambda_{i}({\rm Sym}({\mathbf M}))} & \mbox{if} \ I_1 \neq \emptyset; \\
+ \infty & \mbox{otherwise} \end{array}\right.\end{equation}
\end{corollary}

Let us consider an LMI region ${\mathfrak D}$, defined by its characteristic function $ f_{{\mathfrak D}} = {\mathbf L} + {\mathbf M}z+{\mathbf M}^T\overline{z}$, with ${\mathbf M}$ being definite. In this case, Theorem \ref{Mdef} implies ${\mathfrak D} \neq \emptyset$, and we have the following statement.

\begin{theorem} Let an LMI region ${\mathfrak D}$ be defined by its characteristic function $ f_{{\mathfrak D}} = {\mathbf L} + {\mathbf M}z+{\mathbf M}^T\overline{z}$, with ${\mathbf M}$ being definite. Then one of the following cases holds.
\begin{enumerate}
\item[\rm Case 1.] $\mathbf M$ is positive definite. Then ${\mathfrak D}_{\mathbb R} = (-\infty, x_{max})$, where  $$x_{max} =  \min_{i \in [n]}\frac{-\lambda_i({\mathbf L}_{\mathbf M})}2.$$
\item[\rm Case 2.] $\mathbf M$ is negative definite. Then ${\mathfrak D}_{\mathbb R} = (x_{min}, +\infty)$, where $$x_{min} = \max_{i \in [n]}\frac{-\lambda_i({\mathbf L}_{\mathbf M})}2. $$
Here $\{\lambda_i({\mathbf L}_{\mathbf M})\}_{i=1}^n$ are the eigenvalues of the matrix ${\mathbf L}_{\mathbf M}$.
\end{enumerate}
\end{theorem}

{\bf Proof.} The proof obviously follows from Lemma \ref{redsym} and the previous reasoning. $\square$

Now let us consider the case of an arbitrary region $\mathfrak D$, defined by its characteristic function $ f_{{\mathfrak D}} = {\mathbf L} + {\mathbf M}z+{\mathbf M}^T\overline{z}$. In this case both the matrices $\mathbf L$ and $\mathbf M$ may be indefinite. Consider the case when ${\rm Sym}({\mathbf M})$ is nonsingular. Lemma \ref{condempt} shows, that if ${\mathfrak D} \neq \emptyset$, we still have that matrices ${\mathbf L}$ and ${\rm Sym}({\mathbf M})$ are simultaneously diagonalizable by congruence.

\begin{theorem} Given an LMI region ${\mathfrak D}$, defined by its characteristic function $ f_{{\mathfrak D}} = {\mathbf L} + {\mathbf M}z+{\mathbf M}^T\overline{z}$ with ${\rm Sym}({\mathbf M})$ being nonsingular. Then $\mathfrak D$ is non-empty if and only if the following two conditions hold:
\begin{enumerate}
\item[\rm 1.] ${\mathbf C} = ({\rm Sym}(\mathbf M))^{-1}{\mathbf L}$ is diagonalizable and has real eigenvalues, i.e. there is a nonsingular ${\mathbf R} \in {\mathcal M}^{n \times n}$ such that ${\mathbf R}^{-1}{\mathbf C}{\mathbf R}$ is a real diagonal matrix.
\item[\rm 2.] $x_{min} < x_{max},$
where
\begin{equation}\label{genbound1} x_{min} = \left\{\begin{array}{ccc} \max_{i \in I_2}\frac{-\lambda_i(\widetilde{{\mathbf L}})}{2\lambda_i(\widetilde{{\rm Sym}({\mathbf M})})} & \mbox{if} \ I_2 \neq \emptyset; \\
- \infty & \mbox{otherwise} \end{array}\right.\end{equation}
\begin{equation}\label{genbound2} x_{max} = \left\{\begin{array}{ccc} \min_{i \in I_1}\frac{-\lambda_i(\widetilde{{\mathbf L}})}{2\lambda_{i}(\widetilde{{\rm Sym}({\mathbf M})})} & \mbox{if} \ I_1 \neq \emptyset; \\
+ \infty & \mbox{otherwise} \end{array}\right.\end{equation}
\end{enumerate}
where $\widetilde{{\mathbf L}} = {\mathbf R}{\mathbf L}{\mathbf R}^T$, $\widetilde{{\rm Sym}({\mathbf M})} = {\mathbf R}{\rm Sym}({\mathbf M}){\mathbf R}^T$.
In this case, ${\mathfrak D}_{\mathbb R} = (x_{min}, x_{max})$.
\end{theorem}
{\bf Proof.} $\Rightarrow$ Let ${\mathfrak D} \neq \emptyset$. Then, by Lemma \ref{condempt}, Condition 1 holds. By Theorem \ref{ns}, we obtain that $\mathbf L$ and ${\rm Sym}({\mathbf M})$ are simultaneously diagonalizable by congruence:
$${\Lambda}_{\mathbf L} = {\mathbf S}{\mathbf L}{\mathbf S}^T,$$
$${\Lambda}_{{\rm Sym}({\mathbf M})} = {\mathbf S}{\rm Sym}({\mathbf M}){\mathbf S}^T,$$
where ${\mathbf S} = {\mathbf Q}{\mathbf R}$, for some orthogonal matrix $\mathbf Q$. Thus
$${\Lambda}_{\mathbf L} = {\mathbf Q}{\mathbf R}{\mathbf L}{\mathbf R}^T{\mathbf Q}^T,$$
$${\Lambda}_{{\rm Sym}({\mathbf M})} = {\mathbf Q}{\mathbf R}{\rm Sym}({\mathbf M}){\mathbf R}^T{\mathbf Q}^T,$$
are similar to the matrices ${\mathbf R}{\mathbf L}{\mathbf R}^T$ and ${\mathbf R}{\rm Sym}({\mathbf M}){\mathbf R}^T$, respectively.

Since the diagonal matrices ${\Lambda}_{\mathbf L}$ and ${\Lambda}_{{\rm Sym}({\mathbf M})}$ obviously commute, we apply Corollary \ref{coest} and obtain the required estimates.

$\Leftarrow$ The inverse direction obviously follows from Theorem \ref{ns} and the invariance of ${\mathfrak D}_{\mathbb R}$ under congruence transformations. $\square$

Note that if the LMI region $\mathfrak D$ is composite, i.e. ${\mathfrak D} = {\mathfrak D}^1\cap {\mathfrak D}^2$ then ${\mathfrak D}_{\mathbb R} = {\mathfrak D}^1_{\mathbb R}\cap {\mathfrak D}^2_{\mathbb R}$.

{\bf Step 2.} We fix $\hat{x} \in {\mathfrak D}_{\mathbb R}$, and apply the shift along the real axes $x':=x - \hat{x}$. We calculate  $x'_{min} = x_{min} - \hat{x}$ and $x'_{max} = x_{max} - \hat{x}$. We calculate ${\mathbf L}(\hat{x}) = {\mathbf L} + 2\hat{x}{\rm Sym}({\mathbf M})$. By Lemma \ref{zero}, it is negative definite.

 {\bf Step 3.} Now, having fixed $x = x_0 \in {\mathfrak D}$, we consider the intersection $${\mathfrak D}_{x} = {\mathfrak D}\cap\{z = x+iy \in {\mathbb C}: x = x^0\}.$$ It can be described by substitution $x = x_0$ into \eqref{LMI}:
 $${\mathfrak D}_{x_0} = \{y \in {\mathbb R}: {\mathbf L} + 2x_0{\rm Sym}({\mathbf M}) + 2iy({\rm Skew}({\mathbf M})) \prec 0 \}.$$

If ${\rm Skew}({\mathbf M}) = 0$, we easily obtain that ${\mathfrak D}_{x_0} = \{z = x+iy \in {\mathbb C}: x = x^0\}$ whenever $x_0 \in {\mathfrak D}$. Now we consider the case when ${\rm Skew}({\mathbf M}) \neq 0$.

\begin{theorem}\label{imag} Let an LMI region ${\mathfrak D}$ be defined by its characteristic function $ f_{{\mathfrak D}} = {\mathbf L} + {\mathbf M}z+{\mathbf M}^T\overline{z}$ with ${\rm Skew}({\mathbf M}) \neq 0$. Let $x_0 \in {\mathfrak D}_{\mathbb R}$. Then
$${\mathfrak D}_{x_0} = \{z=(x_0,y) \in {\mathbb C}: |y| < \frac{1}{2 \max_j|\lambda_j({\rm Skew}({\mathbf M})_{{\mathbf L}(x_0)})|}\},$$
where $\lambda_j({\rm Skew}({\mathbf M})_{{\mathbf L}(x_0)})$  are the eigenvalues of the matrix ${\rm Skew}({\mathbf M})_{{\mathbf L}(x_0)}$.
\end{theorem}
{\bf Proof.} Since $x_0 \in {\mathfrak D}_{\mathbb R} \subset {\mathfrak D}$, the matrix ${\mathbf L}(x_0) = {\mathbf L} + 2x_0{\rm Sym}({\mathbf M})$ is obviously negative definite. Since ${\mathbf L}(x_0)$ is negative definite and ${\rm Skew}({\mathbf M})$ is skew-symmetric, by Lemma \ref{redskew} they can be simultaneously reduced by congruence to $- {\mathbf I}$ and some matrix of the form \eqref{qdiagskew}, respectively. Here the Form \eqref{qdiagskew} corresponds to the skew-symmetric matrix ${\rm Skew}({\mathbf M})_{{\mathbf L}(x_0)} = {\mathbf T}^{-1}{\rm Skew}({\mathbf M})({\mathbf T}^{-1})^T$, where ${\mathbf T}{\mathbf T}^T = - {\mathbf L}(x_0)$. Consider
$${\mathbf \Lambda}_{{\rm Skew}({\mathbf M})_{{\mathbf L}(x_0)}} = {\rm diag}\left\{\begin{vmatrix}0 & \nu_1 \\ -\nu_1 & 0 \end{vmatrix}, \ \ldots, \ \begin{vmatrix}0 & \nu_k \\ -\nu_k & 0 \end{vmatrix}, \ 0, \ldots, \ 0\right\},$$
where $\pm i\nu_j$ are the pairs of pure imaginary conjugate eigenvalues of ${\rm Skew}({\mathbf M})_{{\mathbf L}(x_0)}$. Without loss the generality, we assume $\nu_j > 0$.

By Lemmas \ref{App3} and \ref{App2}, we have the set of conditions of the form
$$\begin{vmatrix}-1 & 2iy\nu_j \\ -2iy\nu_j & -1 \end{vmatrix}> 0, \qquad j = 1, \ldots, \ k,$$
which imply $1 - 4\nu_j^2y^2 = (1 - 2\nu_jy)(1+ 2\nu_jy)> 0$. From here we derive $|y| < \frac{1}{2\nu_j}$ for any $j = 1, \ \ldots, \ k$ and thus $|y| < \frac{1}{2\max_j \nu_j}$.
$\square$

{\bf Step 4}. Summarizing the results, we get the following statement.
\begin{theorem}\label{placement} Let a  LMI region ${\mathfrak D}$ be defined by its characteristic function $ f_{{\mathfrak D}} = {\mathbf L} + {\mathbf M}z+{\mathbf M}^T\overline{z}$ with ${\rm Sym}({\mathbf M})$ being nonsingular and ${\rm Skew}({\mathbf M}) \neq 0$. Then the following inclusion holds:
 $$D(x,r) \subseteq {\mathfrak D},$$
where $D(x_0,r)$ is a closed disk with the center $x_0 \in (x_{min}, x_{max})$, where
$x_{min}$ and $x_{max}$ are defined by Formulae \eqref{genbound1} and \eqref{genbound2}, respectively
and the radius $r$ is defined by Formula \eqref{radius}.
\end{theorem}

Let us introduce the following characteristic of an LMI region $\mathfrak D$:
\begin{equation}\omega_{\mathfrak D}:= \inf_{x \in {\mathfrak D}_{\mathbb R}}\frac{r(x)}{|x|}. \end{equation}

\section{Examples of LMI regions with a view to applications}
Here, we focus on the following seven most studied regions.

\subsection{Conic sector with apex at the origin and inner angle $2\theta$} Recall that the simplest characteristic function, which defines the conic region (see Figure 4)
\begin{equation}\label{conereg1}
{\mathfrak D} = \{z = x+iy \in {\mathbb C}: x < 0; -x\tan\theta < y < x\tan\theta\},
\end{equation}
with $0 < \theta < \frac{\pi}{2}$ is as follows (see, for example, \cite{CGA}, \cite{LIU}):
$$f_{\mathfrak D} = \begin{pmatrix} \sin(\theta) & \cos(\theta) \\ - \cos(\theta) & \sin(\theta) \\ \end{pmatrix}z + \begin{pmatrix} \sin(\theta) & -\cos(\theta) \\  \cos(\theta) & \sin(\theta) \\ \end{pmatrix}\overline{z}.$$

\begin{figure}[h]
\center{\includegraphics[scale=0.5]{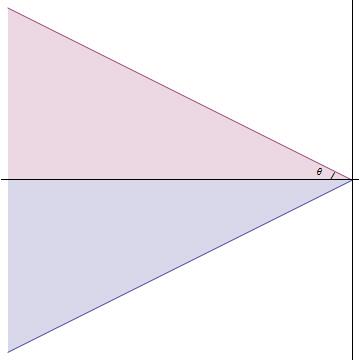}}
\caption{Conic sector}
\end{figure}

In this case, the main characteristics of $\mathfrak D$ are as follows:
\begin{enumerate}
\item[\rm 1.] ${\mathfrak D}_{rc} = {\mathfrak D};$
\item[\rm 2.] $L_{{\mathfrak D}} = \{0\}$;
\item[\rm 3.] $\theta({\mathfrak D}) = \theta$;
\item[\rm 4.] ${\mathfrak D}_{R} = (-\infty, \ 0)$;
\item[\rm 5.] $r({\mathfrak D},x) = |x|\cos(\theta)$.
\item[\rm 6.] $\omega_{\mathfrak D} = \cos(\theta)$.
\end{enumerate}

Consider the examples of problems which lead to the localization of matrix eigenvalues inside Region \ref{conereg1}.

{\bf Example 1. Transient properties of a first-order dynamical system.} Given a continuous-time system of the form
\begin{equation}\label{systprop} \dot{x}(t) = {\mathbf A}x(t), \end{equation}
where ${\mathbf A} \in {\mathcal M}^{n \times n}$, $x(t) \in {\mathbb R}^n$. Then the condition $\sigma({\mathbf A}) \subset {\mathfrak D}$ is referred as {\it relative (sector) stability} of System \ref{systprop} and $\tan(\theta) > 0$ measures the {\it minimal damping ratio} of System \ref{systprop} (see \cite{GUJU}, \cite{DBE} and many others).

{\bf Example 2. Asymptotic stability of a fractional-order system.} Given a fractional-order system of the form
\begin{equation}\label{systpropfr} x^{(\alpha)}(t) = {\mathbf A}x(t), \end{equation}
where $1 < \alpha < 2$, ${\mathbf A} \in {\mathcal M}^{n \times n}$, $x(t) \in {\mathbb R}^n$. It is known to be asymptotically stable if and only if $|\arg(\lambda)| > \frac{\alpha\pi}{2}$ (see \cite{SZQ}).

\subsection{Sliced conic sector}  Consider a region, defined by the following inequalities (see, for example, \cite{STA}, \cite{RAS}, \cite{MAOC}):
\begin{equation}\label{conereg}
{\mathfrak D} = \{z = x+iy \in {\mathbb C}: -x\tan\theta < y < x\tan\theta; \ x < \delta\},
\end{equation}
with $0 < \theta < \frac{\pi}{2}$, $\delta < 0$. This is a part of a conic sector \eqref{conereg}, bounded by a line $x = \delta$ (see Figure 5). It is easy to see that the simplest characteristic function, which defines this LMI region is as follows (see, for example, \cite{CGA}):
$$f_{\mathfrak D} = \begin{pmatrix} 0 & 0 & 0 \\
0 & 0 & 0 \\
0 & 0 & -2\delta \\
 \end{pmatrix} + \begin{pmatrix} \sin(\theta) & \cos(\theta) & 0 \\ - \cos(\theta) & \sin(\theta) & 0 \\
  0 & 0 & 1\\ \end{pmatrix}z + \begin{pmatrix} \sin(\theta) & -\cos(\theta) & 0\\  \cos(\theta) & \sin(\theta)& 0 \\ 0 & 0 & 1\\ \end{pmatrix}\overline{z}.$$

\begin{figure}[h]\label{slice}
\center{\includegraphics[scale=0.5]{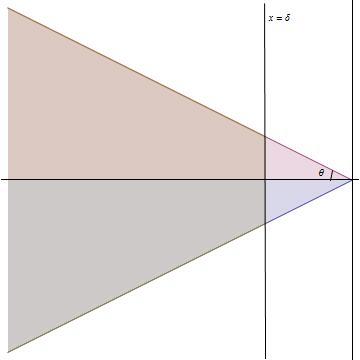}}
\caption{Sliced conic sector}
\end{figure}

This LMI region is a desired stability region for preserving specified settling time and damping ratio (see \cite{GUJU}, \cite{GUT2})
Here, we ensure minimum decay rate $\delta$ and minimum damping ratio $\tan(\theta)$.

In this case, the main characteristics of $\mathfrak D$ are as follows:
\begin{enumerate}
\item[\rm 1.] ${\mathfrak D}_{rc} = {\mathfrak D};$
\item[\rm 2.] $L_{{\mathfrak D}} = \{0\}$;
\item[\rm 3.] $\theta({\mathfrak D}) = \theta$;
\item[\rm 4.] ${\mathfrak D}_{R} = (-\infty, \ \delta)$;
\item[\rm 5.] $r({\mathfrak D}, x) = |x - \sigma|\cos(\theta)$.
\end{enumerate}

\subsection{Shifted disk} The following LMI region received particular attention in literature (see, for example, \cite{FK}, \cite{HON}, \cite{GU}, \cite{SKK}, \cite{XSX}, \cite{ZHAS} and many others). Given an (open) disk $D(a,r)$, centered at $a \in {\mathbb R}$ of the radius $r$ (see Figure 6), it can defined by the following characteristic function (see \cite{CHG}, \cite{CGA}, \cite{LIU}):
$$f_{D(a,r)} = \begin{pmatrix} -r & -a \\  -a & -r \\ \end{pmatrix} + \begin{pmatrix} 0 & 1 \\  0 & 0 \\ \end{pmatrix}z + \begin{pmatrix} 0 & 0 \\  1 & 0 \\ \end{pmatrix}\overline{z}.$$
A special case $a=0$, $r =1$ gives the well-studied unit disk $D(0,1)$:
$$f_{D(0,1)} = \begin{pmatrix} -1 & 0 \\  0 & -1 \\ \end{pmatrix} + \begin{pmatrix} 0 & 1 \\  0 & 0 \\ \end{pmatrix}z + \begin{pmatrix} 0 & 0 \\  1 & 0 \\ \end{pmatrix}\overline{z}.$$
Here ${\rm Sym}({\mathbf M}) = \begin{pmatrix} 0 & 1 \\  1 & 0 \\ \end{pmatrix}$ is obviously indefinite, thus by Theorem \ref{bounded}, the LMI region is bounded. By Lemma \ref{zero}, if $\det(\mathbf L) = (r-a)(r+a) > 0$, then $0 \in D(a,r)$.
\begin{figure}[h]
\center{\includegraphics[scale=0.5]{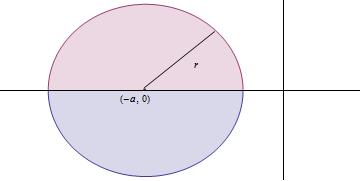}}
\caption{Shifted disk $D(a,r)$}
\end{figure}
Due to its boundedness, the main characteristics of $\mathfrak D$ are as follows:
\begin{enumerate}
\item[\rm 1.] ${\mathfrak D}_{rc} = \{0\};$
\item[\rm 2.] $L_{{\mathfrak D}} = \{0\}$;
\item[\rm 3.] $\theta({\mathfrak D}) = \{0\}$;
\end{enumerate}

Using the shift $x' = x - a$, which maps $a$ to $0$, we obtain $D' = D(0,r)$ with ${\mathbf L}'= \begin{pmatrix} -r & 0 \\  0 & -r \\ \end{pmatrix}$, which is obviously negative definite. Calculating the eigenvalues of the matrix ${\rm Sym}({\mathbf M})_ = \begin{pmatrix} 0 & \frac{1}{2} \\  \frac{1}{2} & 0 \\ \end{pmatrix}$, that are $\pm \frac{1}{2}$, by Corollary \ref{coest} we obtain ${\mathfrak D}'_{\mathbb R} = (-r, \ r)$ and ${\mathfrak D}_{\mathbb R} = (a-r, \ a+r)$.

 Now, using the results of Subsection 6.5, we find $r({\mathfrak D},x)$. Fixing $x \in (-r, r)$, we get ${\mathbf L}(x) = \begin{pmatrix} -r & x \\  x & -r \\ \end{pmatrix}$. By Cholesky decomposition, we get ${\mathbf L}(x) = {\mathbf B}^T{\mathbf B}$, where ${\mathbf B} = \begin{pmatrix}\sqrt{r} & -\frac{x}{\sqrt{r}} \\ 0 & \sqrt{r - \frac{x^2}{r}} \end{pmatrix}$.

  Calculating ${\rm Skew}({\mathbf M})_{{\mathbf L}'(x)} :=({\mathbf B}^T)^{-1}{\rm Skew}({\mathbf M}){\mathbf B}^{-1}$ and its eigenvalues, we get $${\rm Skew}({\mathbf M})_{{\mathbf L}'(x)} = \begin{pmatrix} 0 & \frac{1}{2\sqrt{r^2 - x^2}} \\ -\frac{1}{2\sqrt{r^2 - x^2}} & 0 \end{pmatrix}$$ with the eigenvalues $\pm \frac{i}{2\sqrt{r^2 - x^2}}$ and the corresponding bounds for $y$ are $$-\sqrt{r^2 - x^2} < y < \sqrt{r^2 - x^2}$$ (note, that the exact substitution to the formula $x^2 + y^2 = r^2$ gives us the same result). Then, using Formula \eqref{radius}, we get $$r({\mathfrak D}, x) = \frac{|x-a| \sqrt{r^2 - x^2}}{\sqrt{(x-a)^2 +(r^2-x^2)}}.$$

{\bf Example.} It is well-known (see, for example, \cite{CHEN2}) that stability of discrete-time system
 \begin{equation}\label{systdis} x(k + 1) = {\mathbf A}x(k) \end{equation}
where ${\mathbf A} \in {\mathcal M}^{n \times n}$, $x(k) \in {\mathbb R}^n$ denotes the state vector,
 is equivalent to the localization of the eigenvalues of a system matrix inside the unit disk $D(0,1)$. Now consider spectra localization inside a shifted disk $D(a,r)$. In the case when $|a| + r < 1$, this is a desired stability region for shaping dynamic responses of System \ref{systdis} (see \cite{CHL}).

{\bf Example.} The same concept is considered for time-delay systems. Given a linear discrete time-delay system:
\begin{equation}\label{systdel} x(k + 1) = {\mathbf A}x(k) + {\mathbf A}_dx(k-d), \end{equation}
where  $x(k) \in {\mathbb R}^n$ denotes the state vector, ${\mathbf A}, {\mathbf A}_d \in {\mathcal M}^{n \times n}$, $d > 0$ is a known positive integer. The system \eqref{systdel} is said to be $D(a,r)$-stable if all the (finite) solutions of its characteristic equation satisfy
$$|(z - a)/r|<1 $$
for $r > 0$ and $|a|+r < 1$ (see \cite{HHP}, \cite{LLK}, also see \cite{MAO}, \cite{CHCH} for the case of singular time-delay systems).

\subsection{Vertical strip (real bounding)} Consider the region ${\mathfrak D} = {\mathbb C}^{-\alpha}_{-\beta}$ defined as follows $${\mathbb C}^{-\alpha}_{-\beta}= \{z \in {\mathbb C}; \ -\beta< {\rm Re}(z) < -\alpha; \ 0 < \alpha < \beta \}.$$ This LMI region (see Figure 7) can be represented as an intersection of two first-order LMI regions (see \cite{ZHA}): $${\mathbb C}^{-\alpha}_{-\beta}= {\mathbb C}^{-\alpha} \cap {\mathbb C}^{\beta},$$ where
$${\mathbb C}^{-\alpha}= \{z \in {\mathbb C}; \ {\rm Re}(z) < -\alpha; \ \alpha > 0\};$$
 $${\mathbb C}^{\beta}= \{z \in {\mathbb C}; \ {\rm Re}(z) > -\beta; \ \beta > 0\}.$$  Applying Property 3, we obtain that ${\mathbb C}^{-\alpha}_{-\beta}$ is a second order LMI region with the characteristic function (see \cite{CHG}, \cite{LIY}, \cite{EK})
\begin{equation}\label{sy} f_{\mathfrak D} = \begin{pmatrix} 2\alpha & 0 \\  0 & -2\beta \\ \end{pmatrix} + \begin{pmatrix} 1 & 0 \\  0 & -1 \\ \end{pmatrix}z + \begin{pmatrix} 1 & 0 \\  0 & -1 \\ \end{pmatrix}\overline{z}.\end{equation}

\begin{figure}[h]
\center{\includegraphics[scale=0.5]{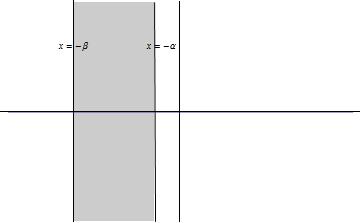}}
\caption{Real bounding}
\end{figure}

The localization of the eigenvalues of System \ref{systprop} inside ${\mathbb C}^{-\alpha}_{-\beta}$ measures the minimal $\alpha$ and the maximal $\beta$ decay rate of the system (see \cite{GUJU}).

Theorem \ref{lin} implies that any nonempty LMI region, defined by its characteristic function  $f_{\mathfrak D} = {\mathbf L} + {\mathbf M}z + {\mathbf M}^T\overline{z}$ of any order $n \geq 2$ with ${\mathbf M}$ being symmetric, is a vertical strip, thus it can be defined by the characteristic function of form \eqref{sy} of the lowest possible order $2$.

In this case, the main characteristics of $\mathfrak D$ are as follows:
\begin{enumerate}
\item[\rm 1.] ${\mathfrak D}_{rc} = {\mathbb I}$ (by Theorem \ref{lin});
\item[\rm 2.] $L_{{\mathfrak D}} = {\mathbb I}$ (by Theorem \ref{lin});
\item[\rm 3.] $\theta({\mathfrak D}) = 0$;
\item[\rm 4.] ${\mathfrak D}_{\mathbb R} = (-\beta,-\alpha)$;
\item[\rm 5.] $r({\mathfrak D}, x) = \min(|x + \beta|, |x + \alpha|)$, whenever $-\beta< x < -\alpha$.
\end{enumerate}
Note that in this case, we do not apply Formula \eqref{radius}, but calculate $r({\mathfrak D}, x)$ directly.

{\bf Example. Interval stability.} The following concept was introduced in \cite{ZHAX}, with a view to the applications to linear stochastic systems. An Ito-type stochastic differential system is called {\it $(-\beta, -\alpha)$-stable} with $0 \leq \alpha < \beta$ if the spectrum of the corresponding linear operator belongs to ${\mathbb C}^{-\alpha}_{-\beta}$. Thus the concept of ${\mathbf D}$-stability with respect to a region ${\mathfrak D} = {\mathbb C}^{-\alpha}_{-\beta}$ coincides with the concept of {\it interval stability} (see \cite{ZHA}, \cite{ZHAX}).

Note, that in \cite{LZH}, when studying an LMI region ${\mathfrak D}$, defined by its characteristic function $ f_{{\mathfrak D}} = {\mathbf L} + {\mathbf M}z+{\mathbf M}^T\overline{z}$, the authors assumed the matrices $\mathbf L$ and $\mathbf M$ to be diagonal (see \cite{LZH}, p. 292, Remark 1). By the above reasoning, this assumption reduces the region $\mathfrak D$ to the case of a vertical strip (halfplane), which can be defined by a characteristic function of order $\leq 2$.

\subsection{Horizontal strip (imaginary bounding)} The localization of the eigenvalues inside the stability region (see Figure 8)
$${\mathfrak D} = \{z \in {\mathbb C}: |y| < w_0, \ w_0 > 0\}$$ corresponds to such transient property of System \eqref{systprop} as {\it bounded frequency}, where $w_0$ measures the maximal damping frequency of the system (see \cite{DBE}, \cite{GUJU}). In this case, $\mathfrak D$ is defined by the characteristic function (see \cite{CHG}, also \cite{EK})
\begin{equation}\label{sk}f_{\mathfrak D} = \begin{pmatrix} -w_0 & 0 \\  0 & -w_0 \\ \end{pmatrix} + \begin{pmatrix} 0 & -1 \\  1 & 0 \\ \end{pmatrix}z + \begin{pmatrix} 0 & 1 \\  -1 & 0 \\ \end{pmatrix}\overline{z}.\end{equation}

\begin{figure}[h]
\center{\includegraphics[scale=0.5]{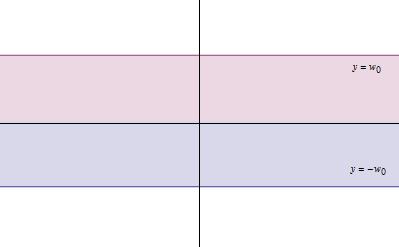}}
\caption{Imaginary bounding}
\end{figure}

Theorem \ref{lin} implies that any nonempty LMI region, defined by its characteristic function  $f_{\mathfrak D} = {\mathbf L} + {\mathbf M}z + {\mathbf M}^T\overline{z}$ of any order $n \geq 2$ with ${\mathbf M}$ being skew-symmetric, is a horizontal strip, thus it can be defined by the characteristic function of form \eqref{sk} of the lowest possible order $2$.

In this case, the main characteristics of $\mathfrak D$ are as follows:
\begin{enumerate}
\item[\rm 1.] ${\mathfrak D}_{rc} = {\mathbb R}$ (by Theorem \ref{semdef});
\item[\rm 2.] $L_{{\mathfrak D}} = {\mathbb R}$ (by Theorem \ref{lin});
\item[\rm 3.] $\theta({\mathfrak D}) = 0$;
\item[\rm 4.] ${\mathfrak D}_{R} = {\mathbb R}$;
\item[\rm 5.] $r({\mathfrak D},(x)) = w_0$.
\end{enumerate}
In this case, we also do not use Formula \ref{radius}.

\subsection{The set $S(\alpha,r,\theta)$} A particularly important for control purposes region $S(\alpha,r,\theta)$ (see \cite{CGA}, \cite{LIY}, \cite{SPS}, \cite{TMG} and many others) is defined as follows (see Figure 9):
$$S(\alpha,r,\theta) = \{z=x+iy \in {\mathbb C}: x < \alpha < 0, |z| < r, \tan(\theta x)< -|y|\}. $$
This composite region of order 5 represents the intersection of the conic sector with the inner angle $\theta$ around the negative direction of the real axis (see Subsection 7.1), the disk $D(0,r)$ of radius $r$ centered at the origin (see Subsection 7.3) and the shifted halfplane ${\mathbb C}_{\alpha}$, $\alpha < 0$ (see Subsection 7.4).
Placing all the eigenvalues of the system \eqref{systprop} in the region $S(\alpha,r,\theta)$ would guarantee a minimum decay rate $\alpha$, a minimum damping ratio $\xi = \cos(\theta)$ and a maximum undamped frequency $w_d = r\sin(\theta)$ (see, for example, \cite{CHG}, \cite{GUJU}, \cite{SZQ}).
\begin{figure}[h]
\center{\includegraphics[scale=0.5]{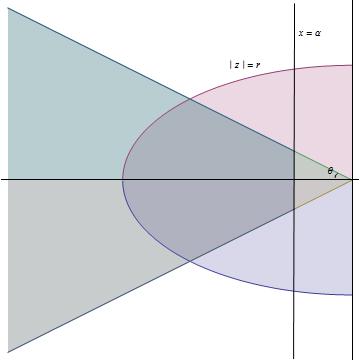}}
\caption{The region $S(\alpha,r,\theta)$}
\end{figure}
The region $S(\alpha,r,\theta)$ is defined by the characteristic function $f_{\mathfrak D}$ with $${\mathbf L} = \begin{pmatrix} -2\alpha & 0 & 0 & 0 & 0 \\ 0 & -r & 0 & 0 & 0 \\ 0 & 0 & -r & 0 & 0 \\ 0 & 0 & 0 & 0 & 0 \\ 0 & 0 & 0 & 0 & 0 \\ \end{pmatrix}$$
and
$${\mathbf M} = \begin{pmatrix} 1 & 0 & 0 & 0 & 0 \\ 0 & 0 & 1 & 0 & 0 \\ 0 & 0 & 0 & 0 & 0 \\0 & 0 & 0 & \sin(\theta) & \cos(\theta) \\0 & 0 & 0 & -\cos(\theta) & \sin(\theta) \\ \end{pmatrix},$$
see, for example, \cite{CHERGA}. In this case, $${\rm Sym}({\mathbf M}) = \begin{pmatrix} 1 & 0 & 0 & 0 & 0 \\ 0 & 0 & \frac{1}{2} & 0 & 0 \\ 0 & \frac{1}{2} & 0 & 0 & 0 \\0 & 0 & 0 & \sin(\theta) & 0 \\0 & 0 & 0 & 0 & \sin(\theta) \\ \end{pmatrix}$$ is obviously indefinite, and ${\rm Skew}({\mathbf M}) \neq 0$. Applying Theorem \ref{bounded}, we get that this region is bounded. Applying Lemma \ref{condempt}, we get that it is empty if $\alpha < -r$.

The main characteristics of $\mathfrak D$ are as follows:
\begin{enumerate}
\item[\rm 1.] ${\mathfrak D}_{rc} = \{0\}$ (by Theorem \ref{bounded});
\item[\rm 2.] $L_{{\mathfrak D}} = \{0\}$ (by Theorem \ref{lin});
\item[\rm 3.] $\theta({\mathfrak D}) = 0$;
\end{enumerate}
Calculating ${\mathfrak D}_{R}$, we apply Lemma \ref{App3}, representing the region $S(\alpha,r,\theta)$ as the intersection of three LMI regions:
$$S(\alpha,r,\theta) = {\mathfrak D}_1 \cap {\mathfrak D}_2 \cap {\mathfrak D}_3,$$
where ${\mathfrak D}_1 = D(0,r)$, ${\mathfrak D}_2 = {\mathbb C}(\theta)$ and ${\mathfrak D}_3 = {\mathbb C}_{\alpha}$. Hence we get $${\mathfrak D}_{R} = {\mathfrak D}^1_{R}\cap {\mathfrak D}^2_{R} \cap {\mathfrak D}^3_R,$$
where ${\mathfrak D}^i_R = {\mathfrak D}_i\cap {\mathbb R}$. By previous subsections, ${\mathfrak D}^1_{R} = (-r,r)$, ${\mathfrak D}^2_{R} = (- \infty, 0)$ and ${\mathfrak D}^3_{R} = (- \infty, \alpha)$ with $\alpha < 0$. Thus ${\mathfrak D}_{R} = (-r, \alpha)$.

Now, using the results of Subsection 6.5, we find $r({\mathfrak D},x)$.

 First, we choose $x_0 \in (-r, \alpha)$, for example, $x_0 = \frac{-r+\alpha}{2}$. Applying the shift $x':=x - x_0$ along the real axis, we obtain the shifted region $S'(\alpha,r,\theta)$, with $0 \in S'(\alpha,r,\theta)$. By Theorem \ref{Shift}, its generating matrix ${\mathbf L'} = {\mathbf L}+{\rm Sym}{\mathbf M}2x_0$, and by Lemma \ref{zero}, it is negative definite. Thus we obtain:
  $${\mathbf L}' = \begin{pmatrix} -r - \alpha & 0 & 0 & 0 & 0 \\ 0 & -r & \frac{-r+\alpha}{2} & 0 & 0 \\ 0 & \frac{-r+\alpha}{2} & -r & 0 & 0 \\ 0 & 0 & 0 & (-r+\alpha)\sin(\theta) & 0 \\ 0 & 0 & 0 & 0 & (-r+\alpha)\sin(\theta) \\ \end{pmatrix},$$
  and consequently,
   $$\tiny{\mathbf L}'(x) = \begin{pmatrix} -r - \alpha + 2x & 0 & 0 & 0 & 0 \\ 0 & -r & \frac{-r+\alpha}{2} + x & 0 & 0 \\ 0 & \frac{-r+\alpha}{2}+x & -r & 0 & 0 \\ 0 & 0 & 0 & (-r+\alpha+2x)\sin(\theta) & 0 \\ 0 & 0 & 0 & 0 & (-r+\alpha+2x)\sin(\theta) \\ \end{pmatrix}.$$
  By Cholesky decomposition, we get ${\mathbf L}'(x) = {\mathbf B}^T{\mathbf B}$, where
  $$\tiny{\mathbf B} = \begin{pmatrix}\sqrt{\alpha + r -2x} & 0 & 0 & 0 & 0 \\ 0 & \sqrt{r} & \frac{-\alpha + r - 2x}{2\sqrt{r}} & 0 & 0 \\
  0 & 0 & \sqrt{r - \frac{(-a+r-2x)^2}{4r}} & 0 & 0 \\ 0 & 0 & 0 & \sqrt{ (2x -r+\alpha)\sin(\theta)} & 0 \\ 0 & 0 & 0 & 0 &  \sqrt{(2x -r+\alpha)\sin(\theta)} \end{pmatrix}$$

 Calculating ${\rm Skew}({\mathbf M})_{{\mathbf L}'(x)} :=({\mathbf B}^T)^{-1}{\rm Skew}({\mathbf M}){\mathbf B}^{-1}$ and its eigenvalues, we get $$\tiny{\rm Skew}({\mathbf M})_{{\mathbf L}'(x)} = \begin{pmatrix} 0 & 0 & 0 & 0 & 0 \\
  0 & 0 & \frac{1}{2\sqrt{r^2 - (x + \frac{-r + \alpha}{2})^2}} & 0 & 0 \\
  0 & -\frac{1}{2\sqrt{r^2 - (x + \frac{-r + \alpha}{2})^2}} & 0 & 0 & 0 \\
  0 & 0 & 0 & 0 & \frac{\cos(\theta)}{(2x-r+\alpha)\sin(\theta)} \\
  0 & 0 & 0 & -\frac{\cos(\theta)}{(2x-r+\alpha)\sin(\theta)} & 0 \\
  \end{pmatrix}$$ with the eigenvalues $\pm \frac{i}{2\sqrt{r^2 - (x + \frac{-r + \alpha}{2})^2}}$, $\pm i\frac{\cos(\theta)}{(2x-r+\alpha)\sin(\theta)}$. Note, that symbolically computed eigenvalues after an easy transformation provide the boundary lines of the LMI region. Thus we can easily calculate $r({\mathfrak D}, x)$ for any $x \in {\mathfrak D}_{\mathbb R}$.

\subsection{Stability parabola} In the study of aeroelastic stability (see, for example, \cite{GUT2}), it is convenient to study the spectra localization in the region
$${\mathfrak D} = \{z = x+iy \in {\mathbb C}: \ y^2 < -\epsilon^2 x \}, $$
i.e. to the left of the stability parabola $y^2 = -\epsilon^2 x$, where $\epsilon$ is a damping parameter (see Figure 10).
\begin{figure}[h]
\center{\includegraphics[scale=0.5]{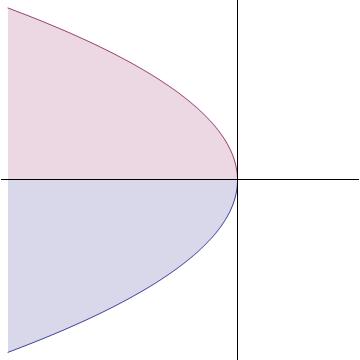}}
\caption{Stability parabola $y^2 = -\epsilon^2 x$}
\end{figure}

In this case, $\mathfrak D$ is a second-order LMI region defined by the characteristic function
$$f_{\mathfrak D} = \begin{pmatrix} -\epsilon^2 & 0 \\  0 & 0 \\ \end{pmatrix} + \begin{pmatrix} \frac{1}{2} & -1 \\  0 & \frac{1}{2} \\ \end{pmatrix}z + \begin{pmatrix} \frac{1}{2} & 0 \\  -1 & \frac{1}{2} \\ \end{pmatrix}\overline{z}=$$
$$\begin{pmatrix} -\epsilon^2 & 0 \\  0 & 0 \\ \end{pmatrix} + \begin{pmatrix} 1 & -1 \\  -1 & 1 \\ \end{pmatrix}x + \begin{pmatrix} 0 & -1 \\  1 & 0 \\ \end{pmatrix}iy, $$
with ${\rm Sym}({\mathbf M}) = \begin{pmatrix} \frac{1}{2}  & -\frac{1}{2}  \\  -\frac{1}{2}  & \frac{1}{2}  \\ \end{pmatrix}$ being positive semidefinite.

The main characteristics of $\mathfrak D$ are as follows:
\begin{enumerate}
\item[\rm 1.] ${\mathfrak D}_{rc} = {\mathbb R}_+$, i.e the positive direction of the real axis (by Corollary \ref{realcone});
\item[\rm 2.] $L_{{\mathfrak D}} = \{0\}$ (by Theorem \ref{lin});
\item[\rm 3.] $\theta({\mathfrak D}) = 0$;
\item[\rm 4.] ${\mathfrak D}_{R} = (-\infty, \ 0)$ (by Lemma \ref{zero});
\end{enumerate}

Now, using the results of Subsection 6.5, we find $r({\mathfrak D},x)$. First, fixing $x \in (-\infty, 0)$, we get ${\mathbf L}(x) = \begin{pmatrix} -\epsilon^2 + x & -x \\  -x & x \\ \end{pmatrix}$. By Cholesky decomposition, we get ${\mathbf L}(x) = {\mathbf B}^T{\mathbf B}$, where ${\mathbf B} = \begin{pmatrix}\sqrt{\epsilon^2 - x} & \frac{x}{\sqrt{\epsilon^2 - x}} \\ 0 & \sqrt{\frac{-x\epsilon^2}{\epsilon^2 - x}} \end{pmatrix}$. Calculating ${\rm Skew}({\mathbf M})_{{\mathbf L}(x)} :=({\mathbf B}^T)^{-1}{\rm Skew}({\mathbf M}){\mathbf B}^{-1}$ and its eigenvalues, we get $${\rm Skew}({\mathbf M})_{{\mathbf L}(x)} = \begin{pmatrix} 0 & \frac{1}{2\sqrt{-\epsilon^2x}} \\ -\frac{1}{2\sqrt{-\epsilon^2x}} & 0 \end{pmatrix}$$ with the eigenvalues $\pm \frac{i}{2\sqrt{-\epsilon^2x}}$ and the corresponding bounds for $y$ are $-\sqrt{-\epsilon^2x} < y < \sqrt{-\epsilon^2x}$ (note, that the exact substitution to the formula $y^2 = -\epsilon^2 x$ gives us the same result). Then, using Formula \eqref{radius}, we get $$r({\mathfrak D}, x) = \frac{-x \sqrt{-\epsilon^2x}}{\sqrt{x^2 -\epsilon^2x}} = \frac{-\epsilon x}{\sqrt{\epsilon^2 - x}}.$$

\end{document}